\documentclass[12pt]{amsart}
\usepackage[margin=2.0cm]{geometry}

\usepackage{amsfonts}
\usepackage{mathrsfs}
\usepackage{cite}
\usepackage{graphicx}

\newcommand{\R}{{\mathbb R}}
\newcommand{\N}{{\mathbb N}}

\newtheorem{theorem}{Theorem}[section]
\newtheorem{claim}[theorem]{Claim}
\newtheorem{lemma}[theorem]{Lemma}
\newtheorem{prop}[theorem]{Proposition}
\newtheorem{corollary}[theorem]{Corollary}

\DeclareMathOperator{\arccosh}{\mathrm{arccosh}}
\DeclareMathOperator{\diam}{\mathrm{diam}}

\title{Stability of the isodiametric problem on the sphere and in the hyperbolic space}
\author{K\'aroly J. B\"or\"oczky, \'Ad\'am Sagmeister }
\address{Alfr\'ed R\'enyi Institute of Mathematics,  Realtanoda u. 13-15, H-1053 Budapest, Hungary} \email{boroczky.karoly.j@renyi.hu}
\address{E\"otv\"os Lor\'and University, Institute of Mathematics, P\'azm\'any P\'eter s\'et\'any 1/c, Budapest, H-1117 Hungary} \email{sagmeister.adam@gmail.com }
\subjclass[2010]{Primary: }

\begin{document}

\maketitle

\begin{abstract}
We prove a stability version of the isodiametric inequality on the sphere and in the hyperbolic space.
\end{abstract}

\section{Introduction}

Let $\mathcal{M}^n$ be either the Euclidean space $\R^n$, hyperbolic space $H^n$ or spherical space $S^n$ for $n\geq 2$. We write $V_{\mathcal{M}^n}$ to denote the $n$-dimensional volume (Lebesgue measure) on $\mathcal{M}^n$,
and $d_{\mathcal{M}^n}(x,y)$ to denote the geodesic distance between $x,y\in\mathcal{M}^n$.
For $x,y\in \mathcal{M}^n$ where $x\neq - y$ if $\mathcal{M}^n=S^n$, we write $[x,y]_{\mathcal{M}^n}$ to denote the geodesic segment between $x$ and $y$ whose length is 
$d_{\mathcal{M}^n}(x,y)$. For pairwise different points $x,y,z\in \mathcal{M}^n$ (which are pairwise not antipodal in the spherical case), let $\angle(x,y,z)$ be the angle of the geodesic segments $[x,y]$ and $[y,z]$ at $y$.

For a bounded set $X\subset \mathcal{M}^n$,  its diameter ${\rm diam}_{\mathcal{M}^n} X$ is the supremum of the geodesic distances $d_{\mathcal{M}^n}(x,y)$ for $x,y\in X$.  
 For $D>0$ and $n\geq 2$, we are considering the maximal volume of a subset of  $\mathcal{M}^n$ of diameter at most $D$. For any $z\in \mathcal{M}^n$ and $r>0$, let
$$
B_{\mathcal{M}^n}(z,r)=\{x\in \mathcal{M}^n:\, d_{\mathcal{M}^n}(x,z)\leq r\}
$$ 
be the $n$-dimensional ball centered at $z$ where it is natural to assume $r<\pi$ if $\mathcal{M}^n=S^n$. When it is clear from the context what space we consider, we drop the subscript referring to the ambient space. In order to speak about the volume of a ball of radius $r$, we fix a reference point $z_0\in \mathcal{M}^n$ where $z_0=o$ the origin if $\mathcal{M}^n=\R^n$. The volume of the unit ball $B_{\R^n}\left(o,1\right)$ is denoted by $\kappa_n$.

The isoperimetric inequality requires the definition of surface area.
For any compact set  $X\subset \mathcal{M}^n$ and $\varrho\geq 0$, we consider the parallel domain
$$
X^{(\varrho)}=\left\{z\in \mathcal{M}^n:\, \exists x\in X\mbox{ \ with \ }
d_{\mathcal{M}^n}(x,z)\leq \varrho\right\}=\bigcup\left\{B_{\mathcal{M}^n}(x,\varrho):\,x\in X\right\}.
$$
If in addition, if ${\rm int}\,X\neq\emptyset$ and $\partial X$ is $(n-1)$-rectifiable; namely, it is the Lipschitz image of some compact subset of $\R^{n-1}$, or if $\partial X$ is a set of positive reach, then the surface area of $X$ can be interpreted
as the outer Minkowski content 
$$
S_{\mathcal{M}^n}(X)=\lim_{\varepsilon\to 0^+}
\frac{V_{\mathcal{M}^n}\left(X^{(\varepsilon)}\right)-V_{\mathcal{M}^n}(X)}\varepsilon
$$
(see Ambrosio, Colesanti, Villa \cite{ACV08}).
We note that if $\partial X$ is $(n-1)$-rectifiable (for example, $X$ is convex), then $S_{\mathcal{M}^n}(X)$ coincides with the $(n-1)$-dimensional Hausdorff measure of $\partial X$ ({\it cf.} \cite{ACV08}).

\begin{theorem}[Isoperimetric inequality]
\label{Isoperimetric}
If $\mathcal{M}^n$ is either $\R^n$, $S^n$ or $H^n$ and $X\subset \mathcal{M}^n$ is compact
and $V_{\mathcal{M}^n}(X)=V_{\mathcal{M}^n}(B(z_0,r))$ for
$r>0$, then
$$
V_{\mathcal{M}^n}\left(X^{(\varrho)}\right)\geq V_{\mathcal{M}^n}(B_{\mathcal{M}^n}(z_0,r+\varrho))
\mbox{ \ for $\varrho>0$}.
$$
\end{theorem}
\noindent {\bf Remark} It follows that if $X$ has outer Minkowski content, then $S_{\mathcal{M}^n}(X)\geq S_{\mathcal{M}^n}(B(z_0,r))$.\\

The isoperimetric inequality was known to the ancient Greeks in the Euclidean plane, and the Euclidean
case in any dimension was proved by the work of Steiner, Schwarz, Weierstrass and Minkowski in the 19th century
(see Gruber \cite{Gru07}). The isoperimetric inequality in the spherical and hyperbolic spaces is due to 
E. Schmidt \cite{Sch43}. We provide the elegant argument by Benyamini \cite{Ben84} because 
it works simultaneously in all spaces of constant curvature, and also the core ideas are essential ingredients for the isodiametric problem, which is our main focus.

Various stability versions of the isoperimetric inequality have been provided starting with Min\-kows\-ki.
In terms of the volume difference, see Fusco, Maggi, Pratelli \cite{FMP08} and Section~\ref{secEuclid}
in the Euclidean case,
 B\"ogelein, Duzaar, Scheven \cite{BDS15} in the hyperbolic case,
and B\"ogelein, Duzaar, N. Fusco \cite{BDF17} in the spherical case.
Essentially optimal stability version of the isoperimetric inequality in $\R^n$ in the terms of Hausdorff distance has been verified by 
 Fuglede \cite{Fug89}.

The main topic of the paper, the isodiametric inequality is proved by Bieberbach \cite{Bie15} in $\R^2$ and by
P. Urysohn \cite{Ury24} in $\R^n$, $n\geq 3$, by
W. Barthel, H. Pabel \cite{BaP87} in $n$-dimensional normed spaces
and by Schmidt \cite{Sch48,Sch49} and
B\"or\"oczky, Sagmeister \cite{BoS20} in the spherical space $S^n$ and the hyperbolic space $H^n$.
The isodiametric problem for bisections in the Euclidean plane is solved by 
Ca\~{n}ete, Merino \cite{CaM21}.

\begin{theorem}[Isodiametric inequality]
\label{Isodiametric}
If $\mathcal{M}^n$ is either $\R^n$, $S^n$ or $H^n$, $D>0$ (with $D<\pi$ if $\mathcal{M}^n=S^n$) and $X\subset \mathcal{M}^n$ is measurable and bounded with 
${\rm diam}X\leq D$, then
$$
V(X)\leq V(B(z_0,D/2)),
$$
and equality holds if and only if the closure of $X$ is a ball of radius $D/2$.
\end{theorem}

In this paper, we provide a new, conceptually more natural proof of the Isodiametric Inequality Theorem~\ref{Isodiametric}, and then we prove the following stability version.
The convex hull of an $X\subset \mathcal{M}^n$ is denoted by ${\rm conv}_{\mathcal{M}^n}X$ (see Section~\ref{secconvex}).

\begin{theorem}
\label{Isodiametricstab}
For $n\geq 2$, if $\mathcal{M}^n$ is either $\R^n$, $S^n$ or $H^n$, $D>0$ (where $D<\frac{\pi}{2}$ if $\mathcal{M}^n=S^n$) and $X\subset \mathcal{M}^n$ is measurable with ${\rm diam}X\leq D$ and
$$
V_{\mathcal{M}^n}\left(X\right)\geq (1-\varepsilon)
V_{\mathcal{M}^n}\left(B_{\mathcal{M}^n}\left(z_0,\frac{D}{2}\right)\right)
$$
for $\varepsilon\in\left[0,\varepsilon_{\mathcal{M}^n}\left(D\right)\right)$, then there exists a $c\in \mathcal{M}^n$ such that
$$
B\left(c,\mbox{$\frac{D}2$}-\gamma_{\mathcal{M}^n}\left(D\right)\cdot\varepsilon^{\frac2{3n+2}}\right)\subset {\rm conv}_{\mathcal{M}^n}X\subset B\left(c,\mbox{$\frac{D}2$}+\gamma_{\mathcal{M}^n}\left(D\right)\cdot\varepsilon^{\frac2{3n+2}}\right)
$$
where $\varepsilon_{\mathcal{M}^n}\left(D\right)>0$ depends on $D$ and $\mathcal{M}^n$ and
$$
\gamma_{\mathcal{M}^n}\left(D\right)=
\left\{
\begin{array}{ll}
  e^{21}n\cdot D&\mbox{ if $\mathcal{M}^n=H^n$ and $D\leq 2$, or $\mathcal{M}^n=\R^n$, or $\mathcal{M}^n=S^n$;}\\[1ex]
  n\cdot e^{7D+8}
 &\mbox{ if $\mathcal{M}^n=H^n$ and $D\geq 1$.}
\end{array} \right. 
$$
In addition, 
$V_{\mathcal{M}^n}\left(\left({\rm conv}_{\mathcal{M}^n}X\right)\backslash X\right)\leq \varepsilon\cdot 
V_{\mathcal{M}^n}\left(B_{\mathcal{M}^n}\left(z_0,\frac{D}{2}\right)\right)$.
\end{theorem}
\noindent {\bf Remark } The order of the error term in Theorem~\ref{Isodiametricstab} is not far from being optimal as if either $\mathcal{M}^n=S^n$ and $D\leq \frac{\pi}2$, or $\mathcal{M}^n=R^n$, or
$\mathcal{M}^n=H^n$, and $X$ is obtained from 
$B_{\mathcal{M}^n}\left(z_0,\frac{D}{2}\right)$ by cutting of a cap of volume $\varepsilon\cdot B_{\mathcal{M}^n}\left(z_0,\frac{D}{2}\right)$ for small $\varepsilon>0$, then
the radius of any ball contained in $X$ is at most
$\frac{D}2-\theta\cdot \varepsilon^{\frac2{n+1}}$
for $\theta>0$ depending on $\mathcal{M}^n$ and $D$. 

The exact value of $\varepsilon_{\mathcal{M}^n}\left(D\right)$ is stated in Theorem~\ref{Isodiametricstab0}.\\

The real importance of Theorem~\ref{Isodiametricstab} lies at the Spherical and Hyperbolic cases as
the known stability results about the Brunn-Minkowski inequality in the Euclidean space
(see Section~\ref{secEuclid} for a review)
 directly yield a stability version of the Isodiametric inequality with better error term  (see Theorem~\ref{IsodiametricstabEuclid} in Section~\ref{secEuclid}). Actually, for $n$-dimensional normed spaces (including the Euclidean space), the stability version of the isodiametric inequality is proved by
Diskant \cite{Dis97}. In addition, Hern\'andez Cifre, Mart\'{\i}nez Fern\'andez \cite{HCF15} proved 
Theorem~\ref{Isodiametricstab} in a more precise form for centrally symmetric subsets of $S^2$ of diameter less than $\pi/2$.

Let $\mathcal{M}^n$ be either $\R^n$,  $H^n$ or $S^n$, and let $D>0$ where
$D<\frac{\pi}2$ if $\mathcal{M}^n=S^n$.
We say that a subset $Z\subset \mathcal{M}^n$ of diameter $D$ is complete if $Z\subset Z'\subset \mathcal{M}^n$ and ${\rm diam}\,Z'=D$ imply $Z=Z'$.
Any $X\subset \mathcal{M}^n$ of diameter $D$ is contained in a complete set
$Z\subset \mathcal{M}^n$ of diameter $D$ (see Section~\ref{secconvex}); therefore, we discuss intensively properties of complete sets. We note that
 complete sets of diameter $D$ are also called convex bodies of constant width ({\it cf.} Section~\ref{secconvex}).

For related surveys on spherical convex bodies,
see Schramm \cite{Sch88b}, Lassak \cite{Las22} and Lassak, Musielak \cite{LaM18}. In addition, Groemer \cite{Gro86} surveys properties of complete sets in Minkowski spaces. For properties of convex bodies of constant width  in the hyperbolic space, see
B\"or\"oczky, Sagmeister \cite{BoS23}.

For some related results in the hyperbolic space, see 
Alfonseca, Cordier, Florentin \cite{ACF21},
Gallego, Reventos, Solanes, Teufel \cite{GRST08},
Jerónimo-Castro, Jimenez-Lopez \cite{JCJL17}.

It follows from the Isodiametric Inequality Theorem~\ref{Isodiametric} that among convex bodies of constant width $D$ in $\mathcal{M}^n$, balls have the maximum volume. However, the problem of the minimum volume of convex bodies of constant width $D$ in $\mathcal{M}^n$ is open if $n\geq 3$ even in the Euclidean case.

If $n=2$, then a Reuleaux triangle in a surface of constant curvature is the intersection of three circular discs of radius $D>0$ whose centers are vertices of a regular triangle of side length $D$ (where we assume $D<\frac{\pi}2$ in the spherical case). It is a convex domain of constant width $D$.
The Blaschke--Lebesgue Theorem, due to
Blaschke \cite{Bla15} and 
Lebesgue \cite{Leb14}
states that amongst bodies of
constant width in the Euclidean plane, the Reuleaux triangle has the minimal area (see Eggleston \cite{Egg52} for a particularly simple proof). The spherical version of the theorem was proved by Leichtweiss \cite{Lei05} based on some ideas of Blaschke. A new proof of the spherical case was recently published by K. Bezdek \cite{Bez21}.
After results by 
Ara\'ujo \cite{Ara97} and
Leichtweiss \cite{Lei05}  if the boundary is piecewise smooth,
B\"or\"oczky, Sagmeister \cite{BoS23} proved that  the Reuleaux triangle has the minimal area
amongst bodies of
constant width in the hyperbolic plane. We note that a stability version of the extremal property of the Reuleaux triangle is proved in \cite{BoS23} in all spaces of constant curvature.

\section{Spaces of constant curvature}
\label{secconstant}

Let $\mathcal{M}^n$ be either $\R^n$,  $H^n$ or $S^n$. Our focus is on the spherical- and hyperbolic space, and we assume that $S^n$ is embedded into $\R^{n+1}$ the standard way,
and $H^n$ is embedded into $\R^{n+1}$ using the hyperboloid model. We write $\langle \cdot,\cdot\rangle$ to denote
the standard scalar product in $\R^{n+1}$, 
and write 
$z^\bot=\{x\in\R^{n+1}:\,\langle x,z\rangle=0\}$ for a $z\in\R^{n+1}\backslash o$.
Fix an $e\in S^n$. In particular, we have
\begin{eqnarray*}
S^n&=& \{x\in\R^{n+1}:\,\langle x,x\rangle=1\}\\
H^n&=& \{x+te:\,x\in e^\bot \mbox{ and }t\geq 1\mbox{ and }t^2-\langle x,x\rangle=1\}.
\end{eqnarray*}

For $H^n$, we also consider the following symmetric bilinear form $\mathcal{B}$ on $\R^{n+1}$: If $x=x_0+te\in\R^{n+1}$ and $y=y_0+se\in\R^{n+1}$ for
$x_0,y_0\in e^\bot$ and $t,s\in\R$, then
$$
\mathcal{B}(x,y)=ts-\langle x_0,y_0\rangle.
$$
In particular,
\begin{equation}
\label{Hn}
\mathcal{B}(x,x)=1\mbox{ for $x\in H^n$}.
\end{equation}
We note that the geodesic distance of $x,y\in \mathcal{M}^n$ where $\mathcal{M}^n$ is either $H^n$ or $S^n$ is 
\begin{eqnarray*}
d_{S^n}\left(x,y\right)&=&\arccos{\langle x,y\rangle}
\mbox{ if $\mathcal{M}^n=S^n$;}\\
d_{H^n}\left(x,y\right)&=&\arccosh\left(\mathcal{B}\left(x,y\right)\right)\mbox{ if $\mathcal{M}^n=H^n$.}
\end{eqnarray*}
In particular, the isometries of $\R^n$ are the maps of the form $x\mapsto Ax+b$ where $A\in O(n)$ and $b\in \R$,
the isometries of $S^n\subset \R^{n+1}$ are the maps of the form $x\mapsto Ax$ where $A\in O(n+1)$,
and the  isometries of $H^n\subset \R^{n+1}$ are the maps of the form $x\mapsto Ax$ where $A\in {\rm GL}(n+1,\R)$
leaves $\mathcal{B}(\cdot,\cdot)$ invariant and $\langle Ae,e\rangle>0$. The isometry group of each
$\mathcal{M}^n$ acts transitively, and the subgroup fixing a $z\in \mathcal{M}^n$ is isomorphic to $O(n-1)$.

Again let $\mathcal{M}^n$ be either $\R^n$,  $H^n$ or $S^n$ using the models as above for $H^n$ and $S^n$. For
$z\in \mathcal{M}^n$, we define the tangent space $T_z$ as
\begin{eqnarray*}
T_z&=&\{x\in\R^{n+1}:\,\mathcal{B}(x,z)=0\} \mbox{ if $\mathcal{M}^n=H^n$}\\
T_z&=&z^\bot\subset \R^{n+1} \mbox{ if $\mathcal{M}^n=S^n$}\\
T_z&=&\R^n \mbox{ if $\mathcal{M}^n=\R^n$}.
\end{eqnarray*}  
We observe that $T_z$ is an $n$-dimensional real vector space equipped with the scalar product
$-\mathcal{B}(\cdot,\cdot)$ if $\mathcal{M}^n=H^n$, and with the scalar product 
$\langle \cdot,\cdot\rangle$ if $\mathcal{M}^n=S^n$ or $\mathcal{M}^n=\R^n$.

For $z\in \mathcal{M}^n$ and unit vector $u\in T_z$, the geodesic line $\ell$ passing through $z$ and determined by $u$ consists of the points
\begin{equation}
\label{geodesic}
p_t=\left\{
\begin{array}{rl}
z\,{\rm cosh}\,t+u\,{\rm sinh}\,t& \mbox{ \ if $\mathcal{M}^n=H^n$}\\
z\cos t+u\sin t& \mbox{ \ if $\mathcal{M}^n=S^n$}\\
z+tu& \mbox{ \ if $\mathcal{M}^n=\R^n$}
\end{array}\right.
\end{equation}
for $t\in \R$. Here the map $t\mapsto p_t$ is bijective onto $\ell$ and satisfies 
$d_{\mathcal{M}^n}(z,p_t)=|t|$ for $t\in \R$ if $\mathcal{M}^n=H^n$ or $\mathcal{M}^n=\R^n$,
and for $t\in(-\pi,\pi]$ if $\mathcal{M}^n=S^n$. If $t>0$  provided $\mathcal{M}^n=H^n$ or $\mathcal{M}^n=\R^n$,
or $0<t<\pi$ provided $\mathcal{M}^n=S^n$, then we say that $u$ points towards $p_t$ along the geodesic segment
$$
[z,p_t]_{\mathcal{M}^n}=\{p_s:\,0\leq s\leq t\}
$$
of length $t$.

A hyperplane $H$ in $\mathcal{M}^n$ passing through the point $z\in \mathcal{M}^n$ and having unit normal
$u\in T_z$, and
 the corresponding half-spaces $H^+$ and $H^-$ where
 $H^+$ has $u$ as the exterior unit normal are defined as follows: $H^-=\mathcal{M}^n\backslash{\rm int} H^+$, and
\begin{equation}
\label{hyperplane}
\begin{array}{lll}
H=\{x\in H^n: \mathcal{B}(x,u)=0\}& \mbox{ }H^+=\{x\in H^n: -\mathcal{B}(x,u)\geq 0\} &\mbox{ \ if $\mathcal{M}^n=H^n$}\\
H=\{x\in S^n: \langle x,u\rangle=0\}& \mbox{ }H^+=\{x\in S^n: \langle x,u\rangle\geq 0\}&\mbox{ \ if $\mathcal{M}^n=S^n$}\\
H=\{x\in \R^n: \langle x,u\rangle=\langle z,u\rangle\}& 
\mbox{ }H^+=\{x\in \R^n: \langle x,u\rangle\geq \langle z,u\rangle\} &\mbox{ \ if $\mathcal{M}^n=\R^n$}.
\end{array}
\end{equation}
The reflection $\sigma_H$ through $H$ is the unique isometry of $\mathcal{M}^n$ different from the identity fixing the points of $H$. In particular, for $x\in \mathcal{M}^n$ where $x\neq \pm u$ in the case $\mathcal{M}^n=S^n$, $H$ is the hyperplane perpendicularly bisecting the segment $[x,\sigma_Hx]_{\mathcal{M}^n}$ 
($H$ going through the midpoint of $[x,\sigma_Hx]_{\mathcal{M}^n}$).

 Let $\mathcal{M}^n$ be either $\R^n$,  $H^n$ or $S^n$.
An important tool to obtain convex bodies with extremal properties is the Blaschke Selection Theorem. 
First we impose a metric on compact subsets.
For a compact set $C\subset\mathcal{M}^n$ and $z\in \mathcal{M}^n$, 
we set $d_{\mathcal{M}^n}(z,C)=\min_{x\in C}d_{\mathcal{M}^n}(z,x)$. 
For any non-empty compact set $C_1,C_2\subset\mathcal{M}^n$, we define their Hausdorff distance
$$
\delta_{\mathcal{M}^n}(C_1,C_2)=\max\left\{
\max_{x\in C_2}d_{\mathcal{M}^n}(x,C_1),
\max_{y\in C_1}d_{\mathcal{M}^n}(y,C_2)\right\}.
$$
The Hausdorff distance is a metric on the space of compact subsets in $\mathcal{M}^n$.
We say that a sequence $\{C_m\}$ of compact subsets of $\mathcal{M}^n$ is bounded if there is a ball containing every $C_m$.
For compact sets $C_m,C\subset \mathcal{M}^n$, we write $C_m\to C$ to denote if
the sequence $\{C_m\}$ tends to $C$ in terms of the Hausdorff distance.

The following statement characterizing limits of compact sets with respect to the Hausdorff distance well-known  (see {\it e.g.} B\"or\"oczky, Sagmeister \cite{BoS20}).

\begin{lemma}
\label{Hausdorffconvergence}
For compact sets $C_m,C\subset \mathcal{M}^n$ where $\mathcal{M}^n$ is either $\R^n$,  $H^n$ or $S^n$, we have $C_m\to C$ if and only if
\begin{description}
\item[(i)] assuming $x_m\in C_m$, the sequence $\{x_m\}$ is bounded and any accumulation point of  $\{x_m\}$  
lies in $C$; and
\item[(ii)] for any $y\in C$, there exist $x_m\in C_m$ for each $m$ such that $\lim_{m\to \infty}x_m=y$.
\end{description}
\end{lemma}

The space of compact subsets of $\mathcal{M}^n$ is locally compact according to the Blaschke Selection Theorem
(see R. Schneider \cite{Sch14}).

\begin{theorem}[Blaschke]
\label{Blaschke}
If $\mathcal{M}^n$ is either $\R^n$,  $H^n$ or $S^n$,
then any bounded sequence of compact subsets of $\mathcal{M}^n$ has
a convergent subsequence.
\end{theorem}

For convergent sequences of compact subsets of $\mathcal{M}^n$, we have the following 
(see B\"or\"oczky, Sagmeister \cite{BoS20}).

\begin{lemma}
\label{limit}
Let $\mathcal{M}^n$ be either $\R^n$,  $H^n$ or $S^n$, and
let the sequence $\{C_m\}$ of compact subsets of $\mathcal{M}^n$ tend to $C$.
\begin{description}
\item[(i)] ${\rm diam}_{\mathcal{M}^n}\,C=\lim_{m\to \infty}{\rm diam}_{\mathcal{M}^n}\,C_m$
\item[(ii)] $V_{\mathcal{M}^n}(C)\geq \limsup_{m\to \infty} V_{\mathcal{M}^n}(C_m)$
\end{description}
\end{lemma}

Recall that for any compact set  $X\subset \mathcal{M}^n$ and $\varrho\geq 0$, the
parallel domain is
$$
X^{(\varrho)}=\{z\in \mathcal{M}^n:\, \exists x\in X\mbox{ \ with \ }
d_{\mathcal{M}^n}(x,z)\leq \varrho\}=\bigcup \{B(x,\varrho):\,x\in X\}.
$$
The triangle inequality and considering $x,y\in X$ with 
$d_{\mathcal{M}^n}(x,y)={\rm diam}_{\mathcal{M}^n}\, X$ lead to the
 following fundamental property of parallel domains.

\begin{lemma}
\label{paralleldiameter}
For $\varrho>0$ and a compact $X\subset \mathcal{M}^n$ where $\mathcal{M}^n$ is either $\R^n$,  $H^n$ or $S^n$, 
and $2\varrho+{\rm diam}\,X<\pi$ if $\mathcal{M}^n=S^n$,
we have 
$$
{\rm diam}\,X^{(\varrho)}=2\varrho+{\rm diam}\,X.
$$
\end{lemma}

We discuss further properties of parallel domains based on Benyamini \cite{Ben84}.

\begin{lemma}
\label{limit-parallel}
If $\varrho\geq 0$ and  $\mathcal{M}^n$ is either $\R^n$,  $H^n$ or $S^n$, and
the sequence $\{C_m\}$ of compact subsets of $\mathcal{M}^n$ tends to $C$, then
\begin{description}
\item[(i)] $\left\{C_m^{(\varrho)}\right\}$ tends to $C^{(\varrho)}$;
\item[(ii)] ${\rm diam}_{\mathcal{M}^n}\,C^{(\varrho)}=
\lim_{m\to \infty}{\rm diam}_{\mathcal{M}^n}\,C_m^{(\varrho)}$;
\item[(iii)] $V_{\mathcal{M}^n}\left(C^{(\varrho)}\right)\geq \limsup_{m\to \infty} V_{\mathcal{M}^n}\big(C_m^{(\varrho)}\big)$;
\item[(iv)] for any $\varepsilon>0$,
$V_{\mathcal{M}^n}\left(C^{(\varrho)}\right)\leq \liminf_{m\to \infty} V_{\mathcal{M}^n}\big(C_m^{(\varrho+\varepsilon)}\big)$.
\end{description}
\end{lemma}
\proof We deduce (i) from Lemma~\ref{Hausdorffconvergence}, and in turn (ii)
from Lemma~\ref{paralleldiameter} and (iii)
 from 
Lemma~\ref{limit} (ii).

For (iv), we only observe that $C^{(\varrho)}\subset C_m^{(\varrho+\varepsilon)}$ if $m$ is large.
\endproof

\section{Two-point symmetrization and a proof of Theorem~\ref{Isodiametric} without the equality case}
\label{sectwopoint}

Let $\mathcal{M}^n$ be either $\R^n$,  $H^n$ or $S^n$, let $H^+$ be a closed half-space bounded by the $(n-1)$-dimensional subspace $H$ in $\mathcal{M}^n$, and let
$X\subset \mathcal{M}^n$ be compact. We write $H^-$ to denote the other closed half-space of $\mathcal{M}^n$ determined by $H$ and
$\sigma_HX$ to denote the reflected image of $X$ through the $(n-1)$-subspace $H$. 

The two-point symmetrization $\tau_{H^+}X$ of $X$ with respect to $H^+$ is a rearrangement of $X$ by replacing
$(H^-\cap X)\backslash \sigma_HX$ by its reflected image through $H$ where readily this reflected image is disjoint from $X$. 
In particular, $\tau_{H^+}X$
can be defined by the properties Lemma~\ref{two-point-prop} (i) and (ii) below.
Naturally, interchanging the role of $H^+$ and $H^-$ results in taking the reflected image of $\tau_{H^+}X$ through $H$. Since this operation does not change any relevant property of the new set, we simply use the notation  
$\tau_HX$ (see Figure 1).

\begin{figure}
\begin{center}
\includegraphics[width=18em]{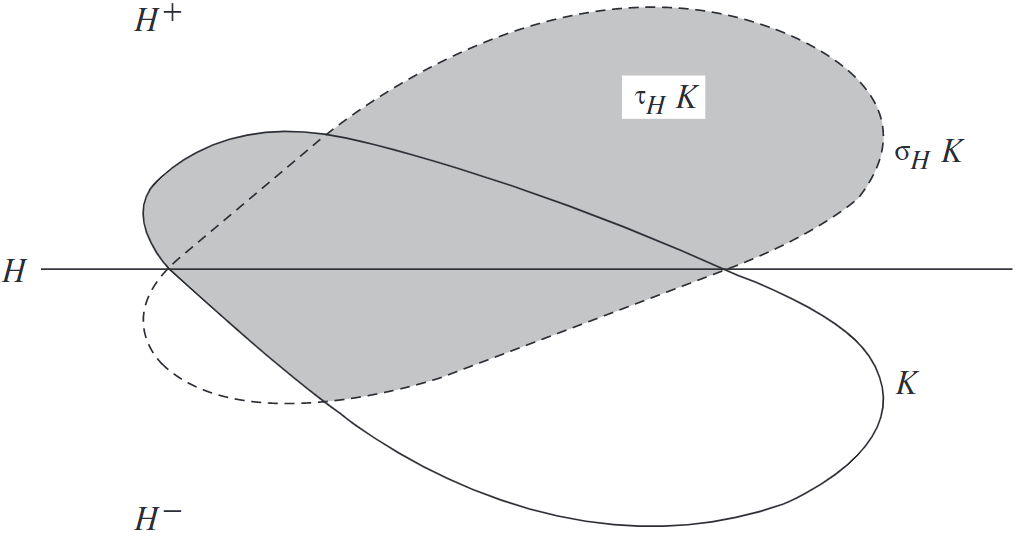}
\end{center}
\caption{ }
\end{figure}  

Two-point symmetrization appeared first in Wolontis \cite{Wol52}. It is applied to prove the
isoperimetric inequality in the spherical space by 
Benyamini \cite{Ben84}, and the spherical analogue of the Blaschke--Santal\'o inequality by
Gao, Hug, Schneider \cite{GHS03} where a crucial step is verified by Aubrun, Fradelizi \cite{AuF04}.

Two-point symmetrization does not lead to an object ``more symmetric'', however,
the definition directly yields that balls are invariant under two-point symmetrization (see Lemma~\ref{two-point-ball}). 

\begin{lemma}
\label{two-point-ball}
Let $\mathcal{M}^n$ be either $\R^n$,  $H^n$ or $S^n$ and let $H^+$ be a half-space of $\mathcal{M}^n$.
\begin{description}
\item[(i)] If $z\in H^+$ and $r>0$ where $r<\pi$ provided $\mathcal{M}^n=S^n$, then
$\tau_HB_{\mathcal{M}^n}(z,r)=B_{\mathcal{M}^n}(z,r)$;
\item[(ii)] if $Y\subset Z\subset \mathcal{M}^n$ compact, then $\tau_HY\subset \tau_HZ$.
\end{description}
\end{lemma}

The following are additional simple properties of two-point symmetrization 
(see B\"or\"oczky, Sagmeister \cite{BoS20}).

\begin{lemma}
\label{two-point-prop}
Let $\mathcal{M}^n$ be either $\R^n$,  $H^n$ or $S^n$, let $H^+$ be a half-space of $\mathcal{M}^n$, and let
$X\subset \mathcal{M}^n$ be compact such that ${\rm diam}_{\mathcal{M}^n}(X)<\pi$ if $\mathcal{M}^n=S^n$. Then $\tau_HX=\tau_{H^+}X$ is compact and satisfies
\begin{description}
\item[(i)] $(\tau_HX)\cap H^+=\big(X\cup\sigma_HX\big)\cap H^+$;
\item[(ii)] $(\tau_HX)\cap H^-=\big(X\cap\sigma_HX\big)\cap H^-$;
\item[(iii)] $V_{\mathcal{M}^n}(\tau_HX)=V_{\mathcal{M}^n}(X)$;
\item[(iv)] ${\rm diam}_{\mathcal{M}^n}(\tau_HX)\leq {\rm diam}_{\mathcal{M}^n}(X)$.
\end{description}
\end{lemma}

Benyamini \cite{Ben84} proved the following property of parallel sets. Since 
\cite{Ben84} is hard to access, we provide a simple argument.

\begin{lemma}
\label{two-point-parallel}
If $\varrho>0$,  $\mathcal{M}^n$ is either $\R^n$,  $H^n$ or $S^n$, and $H^+$ is a half-space 
of $\mathcal{M}^n$, then $V_{\mathcal{M}^n}((\tau_HX)^{(\varrho)})\leq V_{\mathcal{M}^n}(X^{(\varrho)})$
for any compact $X\subset \mathcal{M}^n$. 
\end{lemma}
\proof According to Lemma~\ref{two-point-prop} (iii) applied to $X^{(\varrho)}$, it is sufficient to prove that
\begin{equation}
\label{two-point-parallel-contain}
(\tau_HX)^{(\varrho)}\subset \tau_H\left(X^{(\varrho)}\right).
\end{equation}
Let $z\in (\tau_HX)^{(\varrho)}$, and hence there exists $y\in \tau_HX$
such that $d(y,z)\leq \varrho$. Since the role of $X$ and $\sigma_HX$ are symmetric in the definition of two-point symmetrization, we may assume that $y\in X$, and hence $z\in X^{(\varrho)}$.

If $z\in H^+$, then Lemma~\ref{two-point-prop} (i) applied to $X^{(\varrho)}$ yields that $z\in \tau_H(X^{(\varrho)})$.

If $z\not\in H^+$, then we use that
$$
X^{(\varrho)}\cap \sigma_H(X^{(\varrho)})\subset \tau_H(X^{(\varrho)})
$$
according to Lemma~\ref{two-point-prop} (i) and (ii) applied to $X^{(\varrho)}$. Therefore 
(\ref{two-point-parallel-contain}), and in turn the lemma follows if
\begin{equation}
\label{two-point-intersection}
z\in \sigma_H\left(X^{(\varrho)}\right)=(\sigma_HX)^{(\varrho)}.
\end{equation}

If $y\in \sigma_H X$, then (\ref{two-point-intersection}) readily holds. If
$y\not\in \sigma_H X$, then $y\in H^+$ by Lemma~\ref{two-point-prop} (i), thus
  $d(z,\sigma_Hy)\leq d(z,y)\leq\varrho$, proving (\ref{two-point-intersection}),
and in turn (\ref{two-point-parallel-contain}).
\endproof

For any compact $X\neq \mathcal{M}^n$, let $\mathcal{F}_X$ be the smallest closed subset of
the space of compact subsets of $\mathcal{M}^n$ equipped with the Hausdorff metric containing $X$ and 
being closed under two-point symmetrization. 
Benyamini \cite{Ben84} verified Lemma~\ref{two-point-family} whose proof we present for the convenience of reader.

\begin{lemma}
\label{two-point-family}
If $\mathcal{M}^n$ is either $\R^n$,  $H^n$ or $S^n$,  and 
$X\subset \mathcal{M}^n$, $X\neq \mathcal{M}^n$, is compact, then 
each $Y\in \mathcal{F}_X$ satisfies 
\begin{description}
\item[(i)] 
${\rm diam}_{\mathcal{M}^n}\, Y^{(\varrho)}\leq2\varrho+{\rm diam}_{\mathcal{M}^n}\, X$ 
for any $\varrho\geq 0$;
\item[(ii)] 
$V_{\mathcal{M}^n}(Y^{(\varrho)})\leq V_{\mathcal{M}^n}(X^{(\varrho)})$
for any $\varrho> 0$;
\item[(iii)] 
$V_{\mathcal{M}^n}(Y)= V_{\mathcal{M}^n}(X)$.
\end{description}
\end{lemma}
\proof
For any ordinal $\xi$, we define the subset $\mathcal{F}^{\xi}$ of the space of compact subsets of $\mathcal{M}^n$ by trasnfinite recursion. Let $\mathcal{F}^{0}=\{X\}$. For any ordinal $\xi$, we define
$$
\mathcal{F}^{\xi+1}=\bigcup\left\{\{C,\tau_{H+}C\}:\,
C\in \mathcal{F}^{\xi}\mbox{ and } H^+\subset  \mathcal{M}^n\mbox{ half-space}\right\}.
$$
Finally, if $\xi$ is a limit ordinal, then
$$
\mathcal{F}^{\xi}=\bigcup\left\{\lim_{m\to\infty}C_m:\,C_m\in\mathcal{F}^{\xi_m}\mbox{ where }\xi_m<\xi
\mbox{ and }\lim_{m\to\infty}C_m\mbox{ exists}\right\}.
$$
In particular, $\mathcal{F}^{\alpha}\subset \mathcal{F}^{\xi}$ if $\alpha<\xi$.

It follows from the definition above that if $C\in \mathcal{F}^{\xi}$ for a countable ordinal $\xi$ and $C\neq X$, then
\begin{equation}
\label{family-recursion}
\begin{array}{l}
\mbox{either there exists $\alpha<\xi$, $Z\in\mathcal{F}^{\alpha}$
and a half-space $H^+$ such that $C=\tau_{H+}Z$};\\[1ex]
\mbox{or there exist $\xi_m<\xi$ for $m\in\N$ and $C_m\in\mathcal{F}^{\xi_m}$
such that $\lim_{m\to\infty}C_m=C$}.
\end{array}
\end{equation}

We claim that
\begin{equation}
\label{family-omega1}
\mathcal{F}_{X}=\bigcup\{\mathcal{F}^{\xi}:\,\xi \mbox{ countable ordinal}\}.
\end{equation}
Readily, $\mathcal{F}=\bigcup\{\mathcal{F}^{\xi}:\,\xi \mbox{ countable ordinal}\}\subset \mathcal{F}_{X}$, and $\mathcal{F}$ is closed under two-point symmetrization. To show that $\mathcal{F}$ is a closed subset of
the space of compact subsets of $\mathcal{M}^n$, let $\lim_{m\to\infty}C_m=C$ for 
$C_m\in\mathcal{F}$ and compact $C\subset \mathcal{M}^n$.
Then $C_m\in\mathcal{F}^{\xi_m}$ for some countable ordinals $\xi_m$, and hence $\bigcup_{m\in\N}\xi_m$
is also countable. Let $\xi$ be the smallest (and hence countable) ordinal at least 
$\bigcup_{m\in\N}\xi_m$, thus $C\in \mathcal{F}^{\xi}$, proving (\ref{family-omega1}).

We deduce from (\ref{family-recursion}), (\ref{family-omega1}) and transfinite induction
that each $Y\in \mathcal{F}_X$ satisfies (i) and 
\begin{equation}
\label{VYVX}
V_{\mathcal{M}^n}(Y)\geq  V_{\mathcal{M}^n}(X) 
\end{equation}
where (i) follows from Lemma~\ref{limit-parallel} (ii),  Lemma~\ref{two-point-prop} (iv) 
and Lemma~\ref{paralleldiameter}, and $V_{\mathcal{M}^n}(Y)\geq  V_{\mathcal{M}^n}(X)$ 
follows from Lemma~\ref{limit} (ii)
and Lemma~\ref{two-point-prop} (iii)

For (ii), we again use transfinite induction, so let 
$Y\in \mathcal{F}^{\xi}$ with $Y\neq X$ for a countable ordinal $\xi$. If 
$Y=\tau_{H+}Z$ where $Z\in\mathcal{F}^{\alpha}$ for $\alpha<\xi$ 
and $H^+$ is a half-space, then 
$V_{\mathcal{M}^n}(Y^{(\varrho)})\leq V_{\mathcal{M}^n}(Z^{(\varrho)})\leq V_{\mathcal{M}^n}(X^{(\varrho)})$
follows from Lemma~\ref{two-point-parallel}. Otherwise
there exist $\xi_m<\xi$  and $C_m\in\mathcal{F}^{\xi_m}$ for $m\in\N$
such that $\lim_{m\to\infty}C_m=Y$. For a $\varepsilon>0$,
it follows from transfinite induction that
$V_{\mathcal{M}^n}(C_m^{(\varrho+\varepsilon)})\leq V_{\mathcal{M}^n}(X^{(\varrho+\varepsilon)})$,
thus Lemma~\ref{limit-parallel} (iv) yields
$$
V_{\mathcal{M}^n}(Y^{(\varrho)})\leq \liminf_{m\to \infty} V_{\mathcal{M}^n}(C_m^{(\varrho+\varepsilon)})
\leq V_{\mathcal{M}^n}(X^{(\varrho+\varepsilon)}).
$$
Letting $\varepsilon$ tending to zero, we deduce (ii).

For (iii), it follows from (ii) that
$$
V_{\mathcal{M}^n}(Y)= \lim_{\varrho\to 0^+}  V_{\mathcal{M}^n}(Y^{(\varrho)})
\leq \lim_{\varrho\to 0^+} V_{\mathcal{M}^n}(X^{(\varrho)})= V_{\mathcal{M}^n}(X).
$$
Therefore (\ref{VYVX}) implies (iii).
\endproof

\begin{lemma}
\label{two-point-ball-in-family}
If $\mathcal{M}^n$ is either $\R^n$,  $H^n$ or $S^n$,  and 
$X\subset \mathcal{M}^n$  is compact with $X\neq \mathcal{M}^n$ and $V_{\mathcal{M}^n}(X)>0$, then 
 there exists $B_{\mathcal{M}^n}(z,r)\in \mathcal{F}_X$ 
with $V_{\mathcal{M}^n}(B_{\mathcal{M}^n}(z,r))=V_{\mathcal{M}^n}(X)$
for some $r> 0$ and 
$z\in \mathcal{M}^n$.
\end{lemma}
\proof We choose $z\in \mathcal{M}^n$ and $r> 0$ such that $V(X\cap B)>0$
and $V(B)=V(X)$ for $B=B(z,r)$, and let
$$
\widetilde{F}=\{C\in \mathcal{F}_X:\, C\cap B\neq \emptyset\},
$$
which is a closed subset of $\mathcal{F}_X$.

Setting $D={\rm diam}(X)$, we have $Y\subset B(z,r+D)$ for $Y\in \widetilde{F}$.
Let 
$$
v=\sup\{V(Y\cap B):\,Y\in \widetilde{F}\},
$$
and let $C_m\in \widetilde{F}$, $m\in\N$ such that $\lim_{m\to\infty}V(C_m\cap B)=v$.
Since $C_m\subset B(z,r+D)$ by Lemma~\ref{two-point-prop} (iv), we may assume according to the
Blaschke Selection Theorem~\ref{Blaschke} that $C_m$ tends to a $C\in \widetilde{F}$, and
$C_m\cap B$ tends to a compact set $Z$. It follows from Lemma~\ref{Hausdorffconvergence}
that $Z\subset C$, thus Lemma~\ref{limit} (ii) implies
$$
V(C\cap B)\geq V(Z)\geq
\limsup_{m\to\infty}V(C_m\cap B)=v,
$$
thus $V(C\cap B)=v$.

We suppose that $C\neq B$, and hence $V(W)<V(B)$ holds for
$W=C\cap B$,
and seek a contradiction. We choose an $x\in ({\rm int}\,B)\backslash  C$. Since $V(C)=V(B)$, 
we have $V(C\backslash  B)>0$, therefore there exists a density point $y\in C\backslash  B$; namely,
$$
\lim_{t\to 0^+}\frac{V[(C\backslash  B)\cap B(y,t)]}{V[B(y,t)]}=1.
$$
Let $H$ be the hyperplane in $\mathcal{M}^n$ such that $\sigma_Hx=y$,
and let $H^+$ be the half-space determined by $H$ containing $x$.
Since $x\in B$ and $y\not\in B$, we have $z\in H^+$.
We choose $\varrho>0$ such that $B(x,\varrho)\subset (H^+\cap {\rm int}\,B)\backslash  C$
and $B(y,\varrho)\cap B=\emptyset$, and in particular, we have $V(B(y,\varrho)\cap C)>0$ by the choice of $y$.

It follows from Lemma~\ref{two-point-ball} that $\tau_{H^+}B=B$ and
$\tau_{H^+}W\subset B\cap\tau_{H^+}C$, and Lemma~\ref{two-point-prop} (iii) yields
$V(\tau_{H^+}W)=v$.
In addition, $B(y,\varrho)\cap B=\emptyset$ and $V(B(y,\varrho)\cap C)>0$ imply that 
$B(x,\varrho)\cap \tau_{H^+}W=\emptyset$
and $V(B(x,\varrho)\cap \tau_{H^+}C)>0$, therefore 
$V(B\cap \tau_{H^+}C)>v$. Since $\tau_{H^+}C\in \widetilde{F}$, we have arrived at a contradiction, proving 
that $B=C\in \mathcal{F}_X$. Finally, $V(B)=V(X)$ by Lemma~\ref{two-point-family} (iii).
\endproof

\noindent{\bf Proofs of Theorem~\ref{Isoperimetric} and 
Theorem~\ref{Isodiametric} without the characterization of equality: }
For any compact set $X\subset \mathcal{M}^n$, $X\neq \mathcal{M}^n$,
the family $\mathcal{F}_X$ contains a ball $B_{\mathcal{M}^n}(z,r)$
of the same volume as $X$
by Lemma~\ref{two-point-ball-in-family}. In addition,
Lemma~\ref{two-point-family} implies that
${\rm diam}_{\mathcal{M}^n}\,B_{\mathcal{M}^n}(z,r)\leq {\rm diam}_{\mathcal{M}^n}\, X$
and $V_{\mathcal{M}^n}(B_{\mathcal{M}^n}(z,r+\varrho))\leq V_{\mathcal{M}^n}(X^{(\varrho)})$
for any $\varrho>0$. \hfill $\Box$

\section{$D$-hull, $D$-maximal sets and the uniqueness in Theorem~\ref{Isodiametric}}
\label{secDhull}

Let $\mathcal{M}^n$ be either $\R^n$,  $H^n$ or $S^n$, and let $D>0$ where
$D<\pi$ if $\mathcal{M}^n=S^n$. If ${\rm diam}_{\mathcal{M}^n}X\leq D$ holds for $X\subset \mathcal{M}^n$,  then we define its $D$-hull to be
\begin{equation}
\label{Dhulldef}    
D\mbox{-hull}\,X=\bigcap\{B_{\mathcal{M}^n}(z,D):\,z\in \mathcal{M}^n\mbox{ and }
X\subset B_{\mathcal{M}^n}(z,D)\}.
\end{equation}

\begin{lemma}
\label{Dhull}
Let $\mathcal{M}^n$ be either $\R^n$,  $H^n$ or $S^n$, and let $D>0$ where
$D<\pi$ if $\mathcal{M}^n=S^n$. If ${\rm diam}_{\mathcal{M}^n}\,X\leq D$ 
 holds for $X\subset \mathcal{M}^n$,  then ${\rm diam}_{\mathcal{M}^n}D\mbox{-hull}\,X\leq D$.
\end{lemma} 
\proof First we claim that if $p\in D\mbox{-hull}\,Z$ for a set 
$Z\subset \mathcal{M}^n$ with ${\rm diam}\,Z\leq D$, then
\begin{equation}
\label{Dhullp}
d(x,p)\leq D\mbox{ \ for \ }x\in Z.
\end{equation}
Here (\ref{Dhullp}) follows from the fact that $p\in B(x,D)$ for any $x\in Z$.

Next let $y,z\in D\mbox{-hull}\,X$. We observe that ${\rm diam}\,X'\leq D$ for $X'=X\cup\{y\}$
by (\ref{Dhullp}), and  $D\mbox{-hull}\,X'=D\mbox{-hull}\,X$ by the choice of $y$. In particular, 
$z\in D\mbox{-hull}\,X'$, thus (\ref{Dhullp}) yields that $d(y,z)\leq D$.
\endproof

For $D>0$ where $D<\pi$ if $\mathcal{M}^n=S^n$, we say that a compact set $C\subset \mathcal{M}^n$ 
is $D$-maximal if ${\rm diam}_{\mathcal{M}^n}\,C\leq D$ and
$$
V_{\mathcal{M}^n}(C)=\sup\{V_{\mathcal{M}^n}(X):\, X\subset \mathcal{M}^n\mbox{ compact and }
{\rm diam}_{\mathcal{M}^n} X\leq D \}.
$$
We observe that a ball in $\mathcal{M}^n$ of diameter $D$ is a $D$-maximal set. 

\begin{lemma}
\label{maximal-sets}
Let $\mathcal{M}^n$ be either $\R^n$,  $H^n$ or $S^n$, let
 $D>0$ where $D<\pi$ if $\mathcal{M}^n=S^n$ and let $X\subset \mathcal{M}^n$ with
${\rm diam}_{\mathcal{M}^n}\,X\leq D$.
\begin{description}
\item[(i)] There exists a $D$-maximal set $Z$ in $\mathcal{M}^n$ containing $X$.
\item[(ii)] For any $D$-maximal set $C$ in $\mathcal{M}^n$ and $z\in\partial_{\mathcal{M}^n}C$,
there exists $y\in\partial_{\mathcal{M}^n}C$ such that $d_{\mathcal{M}^n}(z,y)=D$.
\end{description}
\end{lemma}
\proof Let $\{C_m\}$ be a sequence of compact subsets of $\mathcal{M}^n$ with $X\subset C_m$,
${\rm diam}\,C_m\leq D$ and
$$
\lim_{m\to \infty}V(C_m)
=\sup\{V(Y):\, Y\subset \mathcal{M}^n\mbox{ compact and }
X\subset Y\mbox{ and }
{\rm diam} \,Y\leq D \}.
$$
According to the Blaschke Selection Theorem ({\it cf.} Theorem~\ref{Blaschke}), we may asume that
the sequence $\{C_m\}$ tends to a  compact subset $Z\subset \mathcal{M}^n$. Here $Z$ is a $D$-maximal set by Lemma~\ref{limit}.

Next let $C$ be any $D$-maximal set in $\mathcal{M}^n$, and let $z\in\partial_{\mathcal{M}^n}C$.
We suppose that \\
$\widetilde{D}=\max_{x\in C}d(z,x)<D$,
and seek a contradiction. As $z\in\partial C$, there exists some 
$y\in B_{\mathcal{M}^n}(z,\frac12(D-\widetilde{D}))\backslash C$, and hence
$B(y,r)\cap C=\emptyset$ for some $r\in(0,\frac12(D-\widetilde{D}))$.
Therefore $C_0=C\cup  B(z,r)$ satisfies that
$V(C_0)>V(C)$ and ${\rm diam} \,C_0\leq D$,
which is a contradiction verifying (ii).
\endproof

For a compact set $C\subset \mathcal{M}^n$ with $C\neq S^n$, we say that $x\in\partial_{\mathcal{M}^n}C$ is strongly regular if there exist
$r>0$ and  $y,z\in \mathcal{M}^n$ such that 
$x=B_{\mathcal{M}^n}(y,r)\cap B_{\mathcal{M}^n}(z,r)=
C\cap B_{\mathcal{M}^n}(z,r)$ and $B_{\mathcal{M}^n}(y,r)\subset C$. In this case, 
the  exterior unit normal $N_C(x)\in T_x$ to $C$ at $x$ is the unit vector such that 
$-N_C(x)$
points towards $y$ along $[x,y]_{\mathcal{M}^n}$. We deduce from Lemma~\ref{maximal-sets} that boundary points of parallel domains
of $D$-maximal sets are strongly regular.

\begin{lemma}
\label{maximal-sets-strong}
If $\varrho>0$,  $\mathcal{M}^n$ is either $\R^n$,  $H^n$ or $S^n$, 
 $D>0$ with $D+\varrho<\pi$ if $\mathcal{M}^n=S^n$ and $X\subset \mathcal{M}^n$ is a $D$-maximal set,
 then any $x\in\partial_{\mathcal{M}^n} X^{(\varrho)}$ is strongly regular.
\end{lemma}

The upcoming Claim~\ref{ballboundary0} provides a tool to distinguish between points on the boundary of a ball or outside  of the boundary ball, and this tool  is useful to understand boundary structure of a two-point symmetrization.

\begin{claim}
\label{ballboundary0} 
For $B=B_{\mathcal{M}^n}(y_0,R)$ where $\mathcal{M}^n$ is either $\R^n$,  $H^n$ or $S^n$ and $R<\frac{\pi}2$ if $\mathcal{M}^n=S^n$, let $y_1,y_2\in \partial_{\mathcal{M}^n} B$ and $p\in \mathcal{M}^n$ such that 
$y_0\in [y_1,y_2]_{\mathcal{M}^n}$ (and hence $[y_1,y_2]_{\mathcal{M}^n}$ is a diameter of $B$)
and the geodesic line of $[y_1,y_2]_{\mathcal{M}^n}$ does not contain $p$. For 
the hyperplane $H_i$ in $\mathcal{M}^n$ perpendicularly bisecting 
the geodesic segments $[p,y_i]_{\mathcal{M}^n}$, $i=1,2$, we have
$\sigma_{H_1}N_B(y_1)=\sigma_{H_2}N_B(y_2)$ if and only if
$p\in \partial B$. 
\end{claim}

\proof We note that $p,y_1,y_2$ span a two dimensional subspace $\Pi$. Since both $H_1$ and $H_2$ are orthogonal to $\Pi$ if $n\geq 3$, we may actually assume that $n=2$ and $\Pi=\mathcal{M}^2$. In particular, $H_1$ and $H_2$ are lines in this case.

We observe that $\sigma_{H_2}\sigma_{H_1}y_1=\sigma_{H_2}p=y_2$.
Since $\sigma_{H_2}\sigma_{H_2}$ is the identity, $\sigma_{H_1}N_B(y_1)=\sigma_{H_2}N_B(y_2)$ is equivalent to $\sigma_{H_2}\sigma_{H_1}N_B(y_1)=N_B(y_2)$.  As $-N_B(y_1)\in T_{y_1}$ points towards $y_2$ along the segment $[y_1,y_2]_{\mathcal{M}^2}$, and $-N_B(y_2)\in T_{y_2}$ points towards $y_1$ along the segment $[y_2,y_1]_{\mathcal{M}^2}$,
we deduce that $\sigma_{H_1}N_B(y_1)=\sigma_{H_2}N_B(y_2)$ is equivalent to
$\sigma_{H_2}\sigma_{H_1}\left([y_1,y_2]_{\mathcal{M}^2}\right)=\left([y_2,y_1]_{\mathcal{M}^2}\right)$.

Now if $p\in \partial B$, then $d_{\mathcal{M}^2}(y_1,y_0)=d_{\mathcal{M}^2}(p,y_0)=d_{\mathcal{M}^2}(y_2,y_0)$; therefore,
$\left\{y_0\right\}=H_1\cap H_2$ and $H_1$ and $H_2$ are orthogonal. It follows that $\sigma_{H_2}\sigma_{H_1}$ is a rotation around $y_0$ of angle $\pi$, and hence $\sigma_{H_2}\sigma_{H_1}\left([y_1,y_2]_{\mathcal{M}^2}\right)=\left([y_2,y_1]_{\mathcal{M}^2}\right)$.

Finally, let $\sigma_{H_1}N_B(y_1)=\sigma_{H_2}N_B(y_2)$, thus $\sigma_{H_2}\sigma_{H_1}\left([y_1,y_2]_{\mathcal{M}^2}\right)=\left([y_2,y_1]_{\mathcal{M}^2}\right)$ yields
$\sigma_{H_2}\sigma_{H_1}y_0=y_0$.
If $H_1$ and $H_2$ do not intersect, then $\mathcal{M}^2=H^2$, and $\sigma_{H_2}\sigma_{H_1}$ has no fixed points; therefore, this case can't occur. In particular, $H_1$ and $H_2$ intersect in a point $q$.
In this case, $\sigma_{H_2}\sigma_{H_1}$ is a non-trivial rotation around $q$, and hence $\sigma_{H_2}\sigma_{H_1}y_0=y_0$ yields that $\left\{y_0\right\}=\left\{q\right\}=H_1\cap H_2$. We deduce that $d_{\mathcal{M}^2}(y_1,y_0)=d_{\mathcal{M}^2}(p,y_0)=d_{\mathcal{M}^2}(y_2,y_0)$,
thus $p\in \partial B$.  
\endproof

If $\mathcal{M}^n$ is either $\R^n$,  $H^n$ or $S^n$, and $x_1,x_2,x_3\in \mathcal{M}^n$ are not contained in a
geodesic line, then the triangle $[x_1,x_2,x_3]_{\mathcal{M}^n}$ with vertices $x_1,x_2,x_3$ is the union of all geodesic segments starting from $x_1$ and ending at a point of $[x_2,x_3]_{\mathcal{M}^n}$. Actually, it does not matter which vertex is $x_1$, we always obtain the same object. Provided $\{i,j,k\}=\{1,2,3\}$, the angle of the triangle at $x_i$ is the angle of the geodesic segments $[x_i,x_j]_{\mathcal{M}^n}$ and $[x_i,x_k]_{\mathcal{M}^n}$.\\

\noindent {\bf Proof of the uniqueness result in Theorem~\ref{Isodiametric}: } 
Let $\mathcal{M}^n$ be either $\R^n$,  $H^n$ or $S^n$, and let $y_0\in M^n$.
We drop the index notation $\mathcal{M}^n$ as the we always have a fixed ambient space.

For $D>0$ where $D<\pi$ if $\mathcal{M}^n=S^n$, we choose
$\varrho>0$ such that $D+2\varrho<\pi$ if $\mathcal{M}^n=S^n$.
For a $D$-maximal set $X$, which exists by
Lemma~\ref{maximal-sets} (i), we have $V(X)\geq V(B(y_0,\frac{D}2))$
according to Theorem~\ref{Isodiametric}.
It follows from the Isoperimetric Inequality Theorem~\ref{Isoperimetric} that
$V(X^{(\varrho)})\geq V(B(y_0,\frac{D}2+\varrho))$, thus ${\rm diam}\,X^{(\varrho)}\leq D+2\varrho$ implies that $Y=X^{(\varrho)}$
is $(D+2\varrho)$-maximal.

Let $x_1,x_2\in\partial X$ satisfy that $d(x_1,x_2)=D$. It follows that there exist
$y_1,y_2\in\partial Y$ such that $x_1,x_2\in[y_1,y_2]$, and $d(x_i,y_i)=\varrho$ for $i=1,2$.
We may assume that $y_0$ is the midpoint of $[y_1,y_2]$, and we claim first that
\begin{equation}
\label{XrhoB}
B\left(y_0,\frac{D}2+\varrho\right)\subset Y.
\end{equation}
We suppose that (\ref{XrhoB}) does not hold, and seek a contradiction.
Since $B(x_i,\varrho)\subset X^{(\varrho)}$ for $i=1,2$, we deduce the existence of a
$$
p\in \left(\partial Y\cap {\rm int}B\left(y_0,\frac{D}2+\varrho\right)\right)\bigg\backslash[y_1,y_2].
$$

For $i=1,2$, let $H_i$ be the hyperplane perpendicularly bisecting the segment $[y_i,p]$, and let $H_i^+$ be the corresponding half-space with $p\in {\rm int}\,H_i^+$.  We choose $r\in(0,\varrho)$ such that
$$
B(p,r)\subset H_1^+\cap H_2^+.
$$
For $i=1,2$, $-N_Y(y_i)$ points along the segment $[y_1,y_2]$. Since  $d(p,y_0)<d(y_i,y_0)$ yields $\angle(y_0,p,y_i)>\angle(y_0,y_i,p)$ in
the triangle $[y_0,p,y_i]$, we have that $-\sigma_{H_i}N_Y(y_i)\in T_p$ points along a segment connecting $p$ and
a point of the segment $[y_0,y_i]\subset[y_1,y_2]$. Therefore
$$
\sigma_{H_1}N_Y(y_1)\neq -\sigma_{H_2}N_Y(y_2).
$$
In addition, $\sigma_{H_1}N_Y(y_1)\neq \sigma_{H_2}N_Y(y_2)$ according to Claim~\ref{ballboundary0},
therefore after possibly interchanging $y_1$ and $y_2$, we may assume that
\begin{equation}
\label{Npnotpmsigmay1}
N_Y(p)\neq \pm \sigma_{H_1}N_Y(y_1).
\end{equation}

We claim that
\begin{equation}
\label{pontauboundary}
p\in\partial \tau_{H_1}Y.
\end{equation}
Let $p'\in X$ such that $p\in\partial B(p',\varrho)\subset Y$.
According to Lemma~\ref{maximal-sets} (ii), there exists a $q_0\in Y$ such that
$d(p,q)=D+2\varrho$, and hence $Y\subset B(q_0,D+2\varrho)$. It follows that
$-N_Y(p)$ points towards the segment $[p,q_0]$ and $p'\in [p,q_0]$.
On the other hand, $p\in\partial B(x'_1,\varrho)\subset \sigma_{H_1}Y$ for 
$x'_1=\sigma_{H_1}x_1\in[p,\sigma_{H_1}y_2]$, and
$\sigma_{H_1}Y\subset B(\sigma_{H_1}y_2,D+2\varrho)$. Therefore,
\begin{equation}
\label{Nvectors}
\begin{array}{rcl}
N_Y(p)&=& N_{B(q_0,D+2\varrho)}(p)=N_{B(p',\varrho)}(p)\\
N_{\sigma_{H_1}Y}(p)&=&\sigma_{H_1}N_Y(y_1)=N_{B(\sigma_{H_1}y_2,D+2\varrho)}(p)
=N_{B(x'_1,\varrho)}(p).
\end{array}
\end{equation}
Since  $Y\subset B(q_0,D+2\varrho)$ and $\sigma_{H_1}Y\subset B(\sigma_{H_1}y_2,D+2\varrho)$,
$N_Y(p)\neq - \sigma_{H_1}N_Y(y_1)$ ({\it cf.} \eqref{Npnotpmsigmay1}) yields that $p\in\partial (Y\cup \sigma_{H_1}Y)$, which in turn implies 
\eqref{pontauboundary}.

It follows from Lemma~\ref{two-point-prop} that $\tau_{H_1}Y$ is also $D$-maximal,
and hence \eqref{pontauboundary} and  Lemma~\ref{maximal-sets} (ii) yield the existence of a 
$q\in \tau_H Y$ such that
$d(p,q)=D+2\varrho$ and $\tau_HY\subset B(q,D+2\varrho)$. However,
both $B(x'_1,\varrho)\cap H^+\subset \tau_HY$ and $B(p',\varrho)\cap H^+\subset \tau_HY$, therefore
(see also \eqref{Nvectors})
$$
N_Y(p)=N_{B(p',\varrho)}(p)=N_{B(q,D+2\varrho)}(p)=N_{B(x'_1,\varrho)}(p)=\sigma_{H_1}N_Y(y_1).
$$
This contradicts \eqref{Npnotpmsigmay1}, and in turn proves \eqref{XrhoB}.

Since ${\rm diam}\,Y\leq D+2\varrho$, \eqref{XrhoB} implies $X^{(\varrho)}=Y=B\left(y_0,\frac{D}2+\varrho\right)$,
therefore $X=B\left(y_0,\frac{D}2\right)$.
\hfill$\Box$

\bigskip

\noindent{\bf Remark } The condition $N_Y(p)\neq - \sigma_{H_1}N_Y(y_1)$ in \eqref{Npnotpmsigmay1}
is only needed to prove $p\in\partial \tau_{H_1}Y$ ({\it cf.} \eqref{pontauboundary}) when $\mathcal{M}^n=S^n$
and $D\geq \frac{\pi}2$. Otherwise,  \eqref{pontauboundary} simply follows from
$Y\subset B(q_0,D+2\varrho)$ and  $\sigma_{H_1}Y\subset B(\sigma_{H_1}y_2,D+2\varrho)$.

\section{Convex sets and complete sets}
\label{secconvex}

In this section, we start our preparation for the proof of the stability version of the Isodiametric Inequality.
Let $\mathcal{M}^n$ be either $\R^n$,  $H^n$ or $S^n$.
We call $X\subset \mathcal{M}^n$ convex if $[x,y]_{\mathcal{M}^n}\subset X$
for any $x,y\in X$, and in addition, we also assume that $X$ is contained in an open hemisphere if $\mathcal{M}^n=S^n$.
 For $Z\subset \mathcal{M}^n$ where we assume that 
$Z$ is contained in an open hemisphere if $\mathcal{M}^n=S^n$, the convex hull ${\rm conv}_{\mathcal{M}^n} Z$ is the intersection of all convex sets containing $Z$. 

\begin{lemma}
\label{ball-convex}
If either $\mathcal{M}^n=H^n$ or $\mathcal{M}^n=\R^n$ and $r>0$, or $\mathcal{M}^n=S^n$
and $r\in\left(0,\frac{\pi}2\right)$, then $B_{\mathcal{M}^n}(z,r)$ is convex for any $z\in \mathcal{M}^n$. 
\end{lemma}

We remark that 
 $B_{S^n}(z,r)$ is not convex if $r\in\left[\frac{\pi}2,\pi\right)$.

\begin{lemma} 
\label{convexhull}
Let $\mathcal{M}^n$ be either $\R^n$,  $H^n$ or $S^n$, and let
$X\subset \mathcal{M}^n$ be compact and satisfy
${\rm diam}\,X< \frac{\pi}2$ in the case of  $\mathcal{M}^n=S^n$. Then
\begin{description}
\item[(i)] ${\rm diam}_{\mathcal{M}^n}{\rm conv}_{\mathcal{M}^n}\,X= {\rm diam}_{\mathcal{M}^n}\,X$;
\item[(ii)] ${\rm conv}_{\mathcal{M}^n}\,X$ is compact;
\item[(iii)] $V_{\mathcal{M}^n}\left({\rm conv}_{\mathcal{M}^n}\,X\right)>
V_{\mathcal{M}^n}(X)$ if $V_{\mathcal{M}^n}\left({\rm conv}_{\mathcal{M}^n}\,X\right)>0$
and ${\rm conv}_{\mathcal{M}^n}\,X\neq X$.
\item[(iv)] For $D\geq {\rm diam}_{\mathcal{M}^n}\,X$
where $D< \frac{\pi}2$ if  $\mathcal{M}^n=S^n$, the $D$-hull of $X$ is convex.
\end{description}
\end{lemma}
\proof 
For (i), let ${\rm diam}\,X=D$, and let $x_1,x_2\in {\rm conv}_{\mathcal{M}^n} X$.
Since ${\rm conv}_{\mathcal{M}^n}\,X\neq X$ is the intersection of all convex sets containing $X$
and $X\subset B_{\mathcal{M}^n}(x,D)$ for $x\in X$ where $B_{\mathcal{M}^n}(x_1,D)$ is convex
 by Lemma~\ref{ball-convex}, we have $X\subset B_{\mathcal{M}^n}(x_1,D)$, and in turn even 
 ${\rm conv}_{\mathcal{M}^n}\,X\subset B_{\mathcal{M}^n}(x_1,D)$. Therefore $d_{\mathcal{M}^n}(x_1,x_2)\leq D$,
proving (i).

Turning to (ii), if $\mathcal{M}^n=\R^n$, then (ii) follows from Carath\'eodory's theorem (see Schneider \cite{Sch14}); namely, that ${\rm conv}_{\R^n} X$ is the union of the convex hulls of subsets of $X$ of cardinality at most $n+1$. Using suitable radial projection yields the statement in the spherical and in the hyperbolic case.
 
For (iii), we have $V(Z)>0$ for $Z={\rm conv}_{\mathcal{M}^n}\,X$ and $Z\neq X$.
As $Z$ is convex, compact and $V(Z)>0$, it follows that
${\rm int} Z\neq \emptyset$ and the closure of ${\rm int} Z$ is $Z$.
Since  $X$ is compact and $X\neq Z$, there exists some $z\in ({\rm int} Z)\backslash X$. Therefore
$B_{\mathcal{M}^n}(z,r)\subset  ({\rm int} Z)\backslash X$ for some $r>0$, verifying that
$V(Z)>V(X)$.

Finally, (iv) follows from the definition
\eqref{Dhulldef} of the D-hull. 
\endproof

Let $X\subset \mathcal{M}^n$ be compact, and let  $z\in \partial_{\mathcal{M}^n} X$.
We say that $H$ is a supporting hyperplane at $z$ to $X$ if $z\in H$
and $X$ lies in one of the closed half-spaces determined by $H$.
In addition, a ball  $B_{\mathcal{M}^n}(y,r)$ is a supporting ball at $z$
if $z\in B_{\mathcal{M}^n}(y,r)$ but $X\cap{\rm int}_{\mathcal{M}^n}B_{\mathcal{M}^n}(y,r)=\emptyset$.

\begin{lemma}
\label{supporting-hyperplane}
Let $\mathcal{M}^n$ be either $\R^n$,  $H^n$ or $S^n$, and let
$K\subset \mathcal{M}^n$ be convex compact.
\begin{description}
\item[(i)] If $x\not\in K$, and $z\in K$ is a closest point to $x$, then  the hyperplane $H$ passing through $z$ and orthogonal to $[x,z]_{\mathcal{M}^n}$ is a supporting hyperplane to $K$.
\item[(ii)] If $z\in \partial_{\mathcal{M}^n}K$ then there exists a supporting hyperplane $H$ at $z$.
\end{description}
\end{lemma}
\noindent{\bf Remarks } The closest point of a compact convex set to an external point is unique in the Euclidean and the hyperbolic cases, but may not be unique in the spherical case. \\

We also note that if $\varrho>0$ and $K$ is convex in the Euclidean space $\R^n$ or in the hyperbolic space $H^n$, then
$K^{(\varrho)}$ is convex. On the other hand, in the spherical space $S^n$, the parallel domain of a spherical segment is not convex.
\proof\leavevmode
For (i), it is sufficient to prove that if $y\in K$ with $y\neq z$, and $\gamma=\angle(x,z,y)$, then
\begin{equation}
\label{angleclosestpoint}
\gamma\geq\frac{\pi}2.
\end{equation}
Let $d_{\mathcal{M}^n}(x,z)=a$ and $d_{\mathcal{M}^n}(y,z)=b$. For any $t\in[0,b]$, let $y_t\in [z,y]_{\mathcal{M}^n}$
be the point with $d_{\mathcal{M}^n}(y_t,z)=t$, and hence $d_{\mathcal{M}^n}(y_t,x)\geq d_{\mathcal{M}^n}(z,x)$.
Using the Law of Cosines in the triangle $[x,y,z]_{\mathcal{M}^n}$, we deduce that
if $\mathcal{M}^n=\R^n$, then
$$
0\leq \left.\frac{d}{dt}\,d_{\R^n}(y_t,x)^2\right|_{t=0^+}=
\left.\frac{d}{dt}\left(a^2+t^2-2at\cos\gamma\right)\right|_{t=0^+}=-2a\cos\gamma,
$$
if $\mathcal{M}^n=H^n$, then
$$
0\leq \left.\frac{d}{dt}\,\cosh d_{H^n}(y_t,x)\right|_{t=0^+}=
\left.\frac{d}{dt}\left(\cosh a \cdot\cosh t-\sinh a\cdot\sinh t\cdot\cos\gamma\right)\right|_{t=0^+}=-\sinh a\cdot\cos\gamma,
$$
and if $\mathcal{M}^n=S^n$, then
$$
0\geq \left.\frac{d}{dt}\,\cos d_{S^n}(y_t,x)\right|_{t=0^+}=
\left.\frac{d}{dt}\left(\cos a \cdot\cos t+\sin a\cdot\sin t\cdot\cos\gamma\right)\right|_{t=0^+}=\sin a\cdot\cos\gamma.
$$
Therefore, $\cos\gamma\leq 0$, proving \eqref{angleclosestpoint}, and in turn verifying (i).

To prove (ii), we consider a sequence $\{x_m\}$ with $x_m\not\in K$ tending to $z$, and a closest point $z_m\in K$ to $x_m$ for each $z_m$. It follows that $\lim_{m\to\infty}z_m=z$. For each $z_m$, there exists a supporting hyperplane $H_m$ to $K$
passing through $z_m$ according to (i), and it follows from applying the Blaschke Selection Theorem to $H_m\cap B_{\mathcal{M}^n}(z,1)$
that some subsequence $\{H_{m'}\}$ of the sequence $\{H_{m}\}$ tend to a supporting hyperplane $H$ to $K$ at $z$.
\endproof

Let $K\subset\mathcal{M}^n$ be a compact convex set, and let $p\in \partial_{\mathcal{M}^n} K$. We say that
a non-zero $v\in T_p$ is an exterior normal to $K$ at $p$ if there exists
a supporting hyperplane $H$ of $K$ passing through $p$ and orthogonal to $v$, and $-v$ points towards the closed half-space bounded by $H$ and containing $K$. According to Lemma~\ref{supporting-hyperplane}, there exists some exterior unit normal vector $v\in T_p$ at each $p\in \partial_{\mathcal{M}^n} K$. If $p\in \partial_{\mathcal{M}^n} K$ is strongly regular
in the sense of Lemma~\ref{maximal-sets-strong}; namely, there exist
$r>0$ and  $y\in \mathcal{M}^n$ such that $B_{\mathcal{M}^n}(y,r)\subset K$ and $x\in B_{\mathcal{M}^n}(y,r)$, then 
 $N_K(x)\in T_x$ is an exterior normal also in this new sense.

We call a bounded set $K\subset\mathcal{M}^n$  \emph{complete} if $\diam_{\mathcal{M}^n}\left(K\cup\left\{y\right\}\right)>\diam_{\mathcal{M}^n}\left(K\right)$ holds for each $y\in\mathcal{M}^n\backslash K$. Readily, any complete set is closed.

If $\mathcal{M}^n=S^n$, then there exist some unexpected examples of complete sets with obtuse diameter. For example, the $n+2$ vertices of an regular Euclidean $(n+1)$-dimensional simplex inscribed into $S^n$ form a complete set of diameter $\arccos \frac{-1}{n+1}$. Therefore, we will consider only complete sets of acute diameter in the spherical case.

Let $\mathcal{M}^n$ be either $\R^n$, $S^n$ or $H^n$, and let $D>0$
 where $D<\frac{\pi}{2}$ if $\mathcal{M}^n=S^n$.
 We call a compact convex set $K\subset \mathcal{M}^n$ with non-empty interior a convex body.
 We say that a convex body $K\subset\mathcal{M}^n$ is a \emph{convex body of constant width} $D$  if $\diam_{\mathcal{M}^n}\left(K\right)=D$, and for any $p\in\partial_{\mathcal{M}^n}K$ and any outer unit normal $v\in T_p$ to $K$ at $p$, 
  the geodesic segment $\left[p,q\right]_{\mathcal{M}^n}$ of length $D$ and pointing in the direction of $-v$ is contained in $K$. The following statement is well-known (see e.g. Dekster \cite{Dek95} or B\"or\"oczky, Sagmeister \cite{BoS23}).
  
\begin{lemma}
\label{completeness}
Let $\mathcal{M}^n$ be either $\R^n$, $S^n$ or $H^n$, and let
$K\subset\mathcal{M}^n$ be a compact set of diameter $D>0$ where $D<\frac{\pi}{2}$ if $\mathcal{M}^n=S^n$. The following are equivalent:
\begin{description}
\item[(i)] $K$ is a convex body of constant width $D$;
\item[(ii)] $K$ is complete;
\item[(iii)] $K=\bigcap_{x\in K}{B_{\mathcal{M}^n}\left(x,D\right)}$.
\end{description}
\end{lemma}

Theorem~\ref{maximal-sets} guarantees the existence of $D$-maximal sets, and Lemmas~\ref{convexhull} and \ref{completeness} yield

\begin{corollary}
\label{maximal-convex}
If $\mathcal{M}^n$ is either $\R^n$,  $H^n$ or $S^n$, and 
 $D>0$ where $D<\frac{\pi}2$ if $\mathcal{M}^n=S^n$, then
 any $D$-maximal set in $\mathcal{M}^n$ is complete, and hence convex.
\end{corollary}

For a convex body $K$ in $\mathcal{M}^n$ where $\mathcal{M}^n$ be either $\R^n$, $S^n$ or $H^n$, 
the circumradius $R(K)$ is the minimal radius of any ball of $\mathcal{M}^n$ containing $K$. It is well-known (see
Dekster \cite{Dek95} or B\"or\"oczky, Sagmeister \cite{BoS23}) that there exists a unique so-called circumcenter $p\in\mathcal{M}^n$
with $K\subset B_{\mathcal{M}^n}(p,R(K))$, and actually $p\in K$.
In addition,
the inradius 
$r(K)$ is the maximal radius of any $n$-dimensional ball in $K$. There might be various balls of maximal radius contained in $K$ in the Euclidean and the hyperbolic case, as it is shown by the example when $K$ is the parallel domain of a segment.
The following nice connection between the inradius and the circumradius of complete sets in any space of constant curvature is proved by Dekster \cite{Dek95} (see also B\"or\"oczky, Sagmeister \cite{BoS23}).
 
\begin{prop}\label{r_plus_R}
Let $\mathcal{M}^n$ be either $\R^n$, $S^n$ or $H^n$. If 
$K\subset\mathcal{M}^n$ is a complete set of diameter $D>0$ where $D<\frac{\pi}{2}$ in the case $\mathcal{M}^n=S^n$, then
$$
R\left(K\right)+r\left(K\right)=D.
$$
Furthermore, $K$ contains a unique inscribed ball of radius $r(K)$ whose center is the circumcenter $c$ of $K$, and there exists some diameter $[x_1,x_2]_{\mathcal{M}^n}$ of $K$ that contains $c$ where $d_{\mathcal{M}^n}(x_1,c)=R(K)$ and $d_{\mathcal{M}^n}(x_2,c)=r(K)$.
\end{prop}

Now we show that the parallel domains of 
convex bodies of constant width are also
convex bodies of constant width.

\begin{lemma}
\label{completeparallel}
If $\varrho>0$, $\mathcal{M}^n$ is either $\R^n$, $S^n$ or $H^n$, and
$K\subset\mathcal{M}^n$ is a convex body of constant width $D>0$ where $D+2\varrho<\frac{\pi}{2}$ if $\mathcal{M}^n=S^n$, then
$K^{\left(\varrho\right)}$
is a convex body of constant width $D+2\varrho$ satisfying $R\left(K^{\left(\varrho\right)}\right)=R(K)+\varrho$,
$r\left(K^{\left(\varrho\right)}\right)=r(K)+\varrho$
and
\begin{equation}
\label{completeparallel-eq}
K^{\left(\varrho\right)}=\bigcap_{x\in K}{B_{\mathcal{M}^n}\left(x,D+\varrho\right)}.
\end{equation}
\end{lemma}
\proof We drop the reference to $\mathcal{M}^n$ in the formulas to simplify notation.  First we verify \eqref{completeparallel-eq}. For any $y\in K$, we have
$B\left(y,\varrho\right)\subset B\left(x,D+\varrho\right)$ for any $x\in K$; therefore,
$K^{\left(\varrho\right)}\subset\bigcap_{x\in K}{B\left(x,D+\varrho\right)}$.

On the other hand, let $z\in \mathcal{M}^n\backslash K^{\left(\varrho\right)}$, and let
$y\in K$ be closest to $z$. In particular, $d(z,y)>\varrho$. We deduce from Lemma~\ref{supporting-hyperplane} (i)
that the unit vector $v\in T_y$ pointing towards $z$ along $[y,z]$ is an exterior unit normal to $K$ at $y$,
and hence there exists $x\in K$ such that $-v\in T_y$ points towards $x$ along $[y,x]$ and $d(x,y)=D$. It follows that
$d(x,z)=d(x,y)+d(y,z)>D+\varrho$, completing the proof of \eqref{completeparallel-eq}.

Now \eqref{completeparallel-eq} and $D+2\varrho<\frac{\pi}2$ yield that $K^{\left(\varrho\right)}$ is convex, and
Lemma~\ref{paralleldiameter} implies that ${\rm diam}\,K^{\left(\varrho\right)}=D+2\varrho<\frac{\pi}2$. 

To verify that $K^{\left(\varrho\right)}$ has constant width $D+2\varrho$, let
$z\in \partial K^{\left(\varrho\right)}$, and let $v\in T_z$ be an exterior unit normal to $K^{\left(\varrho\right)}$ at $z$. Let
$y\in K$ be closest to $z$, and hence  $d(z,y)=\varrho$. Since $B(y,\varrho)\subset K^{\left(\varrho\right)}$,
$-v\in T_z$ points towards $y$ along $[z,y]$. For the unit vector $w\in T_y$ pointing towards $z$ along $[y,z]$,
$w$ is an exterior unit normal to $K$ at $y$ by Lemma~\ref{supporting-hyperplane} (i), thus
  there exists $x_0\in K$ such that $-v\in T_y$ points towards $x_0$ along $[y,x_0]$ and $d(x_0,y)=D$. 
For the $x\in \partial B(x_0,\varrho)\subset K^{\left(\varrho\right)}$ satisfying that $y,x_0\in[z,x]$, we have 
$d(x,z)=D+2\varrho$, verifying that $K^{\left(\varrho\right)}$ has constant width $D+2\varrho$.
\endproof

\section{Angles and balls}
\label{secangleball}

The main tool to prove Proposition~\ref{Improve} is Lemma~\ref{ballboundary}.  However, we first point out the following two simple observations.

\begin{claim}
\label{sumtwoangles}
The sum of any two angles of  a triangle $T$ is less than $\pi$ provided that either $T$ is a triangle in $\R^n$ or $H^n$,
or $T$ lies in $S^n$ and each side is less than $\pi/2$.
\end{claim}
\proof The only non-trivial case is when $T$ lies in $S^n$. In this case, we may assume that $n=2$. Let $\gamma$ be the third angle of $T$, and let $T'$ be the spherical triangle having an angle $\gamma$ enclosed by two sides of lengths $\pi/2$. Since
$V(T)<V(T')=\gamma$, and $\pi+V(T)$ is the sum of the angles of $T$, we conclude Claim~\ref{sumtwoangles}.
\hfill $\Box$

\begin{claim}
\label{sumtwoanglesmonotone}
Given  $y_0,y_1,z\in \mathcal{M}^n$ not lying on a geodesic arc
and $w\in [z,y_0]_{\mathcal{M}^n}\backslash\{z,y_0\}$
where $\mathcal{M}^n$ is either $\R^n$,  $H^n$ or $S^n$, 
we have 
$$
\angle(y_0,y_1,z)-\angle(y_0,z,y_1)>\angle(y_0,y_1,w)-\angle(y_0,w,y_1).
$$
\end{claim}
\proof We may assume that $n=2$. If $\mathcal{M}^2=H^2$ or $\mathcal{M}^2=\R^2$, then
  applying Claim~\ref{sumtwoangles}
to the triangle $[z,w,y_1]_{\mathcal{M}^2}$ shows that 
$\angle(y_0,z,y_1)<\angle(y_0,w,y_1)$. Since readily $\angle(y_0,y_1,z)>\angle(y_0,y_1,w)$,
we deduce Claim~\ref{sumtwoanglesmonotone} in this case.

Therefore let $\mathcal{M}^2=S^2$. Writing $Y,Z,W$ to denote the angles of the triangle
$T=[y_1,z,w]_{\mathcal{M}^2}$ at $y_1,z,w$, repectively, we deduce that
$$
\angle(y_0,y_1,z)-\angle(y_0,z,y_1)-\left(\angle(y_0,y_1,w)-\angle(y_0,w,y_1)\right)
=Y-Z+\pi-W=2Y-V(T).
$$
However, the ``two-gon'' $M$ containing $T$ and
 bounded by the two half-circles connecting $y_1$ and its antipodal $-y_1$ and passing through $z$ and $w$, respectively,
has area $2Y$, thus
$$
2Y-V(T)=V(M)-V(T)>0,
$$
completing the proof of Claim~\ref{sumtwoanglesmonotone}.
\hfill $\Box$

\begin{lemma}
\label{ballboundary} 
For $B=B_{\mathcal{M}^n}(y_0,R)$ where $\mathcal{M}^n$ is either $\R^n$,  $H^n$ or $S^n$ and 
$R<\frac{\pi}2$ if $\mathcal{M}^n=S^n$, let $y_1\in \partial B$ and $p\in \mathcal{M}^n$ such that $0<\angle(y_1,y_0,p)\leq \frac{\pi}2$,
let  $u\in T_p$  denote the unit vector such that $-u$ points along the segment
$[p,y_0]_{\mathcal{M}^n}$
 and
let the hyperplane $H$ in $\mathcal{M}^n$ be the perpendicular bisector
of the geodesic segment $[p,y_1]_{\mathcal{M}^n}$. 

If $d(p,y_0)\geq R+\eta$  where
$\eta\in\left(0,\frac{R}{2}\right)$ provided $\mathcal{M}^n=H^n$ or $\mathcal{M}^n=\R^n$,
and $\eta\in\left(0,\min\left\{\frac{R}3,\frac{\pi}6-\frac{R}3\right\}\right)$  provided $\mathcal{M}^n=S^n$,
then the angle of
 $\sigma_{H}N_B(y_1)$ and $u$ in $T_p$
is at least 
\begin{eqnarray*}
\frac{2\cdot\eta}{\sqrt{26}\cdot R}>\frac{\eta}{3R}&\mbox{ \ if $\mathcal{M}^n=\R^n$,}\\
\frac{{\rm sinh}\,\eta}{{\rm sinh}\left(\frac{5}2\,R\right)\cdot \sqrt{2}\,{\rm cosh}\,R}>\frac{\eta}{\sinh{5R}}&\mbox{ \ if $\mathcal{M}^n=H^n$,}\\
\frac{3\cdot\eta}{8\sqrt{2}\cdot R}>\frac{\eta}{4R}&\mbox{ \ if $\mathcal{M}^n=S^n$}.
\end{eqnarray*}
\end{lemma}
\proof We observe that the angle of 
$\sigma_{H}N_B(y_1)$ and $u$ in $T_p$ is
$|\angle(y_0,y_1,p)-\angle(y_0,p,y_1)|$. 
Let $p_0\in [p,y_0]_{\mathcal{M}^n}$ such that
$d(p_0,y_0)=R+\eta$, and let $q= \partial B\cap [p,y_0]_{\mathcal{M}^n}$.
Readily
\begin{equation}
\label{anglepymonotone2}
\angle(y_0,y_1,q)=\angle(y_0,q,y_1), 
\end{equation}
and it follows from Claim~\ref{sumtwoanglesmonotone} and (\ref{anglepymonotone2}) that
\begin{equation}
\label{anglepymonotone1}
\angle(y_0,y_1,p)-\angle(y_0,p,y_1)>\angle(y_0,y_1,p_0)-\angle(y_0,p_0,y_1)>
\angle(y_0,y_1,q)-\angle(y_0,q,y_1)=0.
\end{equation}
Therefore Lemma~\ref{ballboundary} holds in Case~1 if
\begin{equation}
\label{ballboundaryp0y1y0} 
\angle(y_0,y_1,p_0)-\angle(y_0,p_0,y_1)\geq \widetilde{\gamma}\eta
\end{equation}
for the suitable $\widetilde{\gamma}>0$ depending on  $\mathcal{M}^n$ and $R$.

We write $P,Q,Y$ to denote the angles of the triangle $[p_0,q,y_1]_{\mathcal{M}^n}$
at $p_0,q,y_1$. Since $\angle(q,y_0,y_1)\leq\frac{\pi}2$ in the triangle 
$[q,y_0,y_1]_{\mathcal{M}^n}$ where the other two angles are $\pi-Q$, we have
$$ 
\tan Q\leq\left\{
\begin{array}{ll}
 \frac{-1}{{\rm cosh}\, R}&\mbox{ \ if $\mathcal{M}^n$ is $H^n$,}\\[1ex]
 -1&\mbox{ \ if $\mathcal{M}^n$ is $\R^n$ or $S^n$}.
\end{array} \right.
$$
In particular, as ${\rm cosh}\, R\geq 1$, we have
\begin{equation}
\label{sinQ} 
\sin Q\geq\left\{
\begin{array}{ll}
 \frac{1}{\sqrt{2}\,{\rm cosh}\, R}&\mbox{ \ if $\mathcal{M}^n$ is $H^n$,}\\[1ex]
 \frac1{\sqrt{2}}&\mbox{ \ if $\mathcal{M}^n$ is $\R^n$ or $S^n$}.
\end{array} \right.
\end{equation}
On the one hand, first applying (\ref{anglepymonotone2}), and after that
Claim~\ref{sumtwoangles}
imply that
\begin{equation}
\label{anglediffY} 
\angle(y_0,y_1,p_0)-\angle(y_0,p_0,y_1)=(Y+\pi-Q)-P>Y.
\end{equation}
On the other hand, we define
$$ 
f_{\mathcal{M}^n}(t)=\left\{
\begin{array}{ll}
t&\mbox{ \ if $\mathcal{M}^n=\R^n$,}\\[1ex]
 {\rm sinh}\, t&\mbox{ \ if $\mathcal{M}^n=H^n$,}\\[1ex]
 \sin t&\mbox{ \ if $\mathcal{M}^n=S^n$}.
\end{array} \right.
$$
We observe that $d_{\mathcal{M}^n}(p_0,y_1)<\frac{5}2\,R$ if $\mathcal{M}^n=H^n$, and $d_{\mathcal{M}^n}(q,y_1)<R$ if $\mathcal{M}^n=S^n$. For $\mathcal{M}^n=\R^n$ since $Q\in\left(\frac{\pi}{2};\frac{3\cdot\pi}{4}\right]$, using the Law of Cosines for the triangle $\left[p_0,q,y_1\right]_{\R^n}$ we deduce $d_{\R^n}(p_0,y_1)<\frac{\sqrt{13}}{2}\,R$. For $\mathcal{M}^n=S^n$ using the Law of Cosines in the triangle $\left[p_0,q,y_1\right]_{S^n}$ we have
\begin{gather*}
\cos\left(d_{S^n}\left(p_0;y_1\right)\right)=\cos\left(\eta\right)\cdot\cos\left(d_{S^n}\left(q;y_1\right)\right)+\sin\left(\eta\right)\cdot\sin\left(d_{S^n}\left(q;y_1\right)\right)\cdot\cos\left(Q\right)\geq\\
\geq\cos\left(\eta\right)\cdot\cos\left(R\right)-\frac{1}{\sqrt{2}}\cdot\sin\left(\eta\right)\cdot\sin\left(R\right)>\cos\left(R+\eta\right)
\end{gather*}
and since by the choice of $\eta$ we know that $R+\eta<\frac{\pi}{2}$, from this we get
$$
\frac{1}{f_{S^n}\left(d_{S^n}\left(p_0,y_1\right)\right)}>\frac{1}{\sin\left(R+\eta\right)}>\frac{3}{4\cdot R}
$$
also using the upper bound $\eta<\frac{R}{3}$. Substituting  $d_{\mathcal{M}^n}(p_0,q)=\eta$ in the Law of Sines for the triangle $[p,q,y_1]_{\mathcal{M}^n}$ and using (\ref{sinQ}) and $\frac{t}2<\sin t<t$ for $t\in\left(0,\frac{\pi}2\right)$, we deduce that
$$
Y>\sin Y=\frac{f_{\mathcal{M}^n}(d_{\mathcal{M}^n}(p_0,q))}{f_{\mathcal{M}^n}(d_{\mathcal{M}^n}(p_0,y_1))}\cdot \sin Q
> \left\{
\begin{array}{ll}
\frac{2\cdot\eta}{\sqrt{26}\cdot R}&\mbox{ \ if $\mathcal{M}^n=\R^n$,}\\[2ex]
\frac{{\rm sinh}\,\eta}{{\rm sinh}\left(\frac{5}2\,R\right)\cdot \sqrt{2}\,{\rm cosh}\,R}&\mbox{ \ if $\mathcal{M}^n=H^n$,}\\[2ex]
\frac{3\cdot\eta}{8\sqrt{2}\cdot R}&\mbox{ \ if $\mathcal{M}^n=S^n$}.
\end{array} \right.
$$
Therefore if $d(p,y_0)\geq R+\eta$, then Lemma~\ref{ballboundary} follows from
(\ref{ballboundaryp0y1y0}) and (\ref{anglediffY}). 
\hfill $\Box$

\section{The gap between two-point symmetrization and its convex hull}
\label{sechulltwopoint}

We write $\angle(x,y)$ to denote the angle of $x,y\in\R^n\backslash\{o\}$.
In Proposition~\ref{Improve}, we use the following observation to estimate the extra volume from below.

\begin{claim}
\label{ballconvexhull}
For $r>0$ and $w_1,w_2\in\R^n$ with $0<\|w_1-w_2\|<2r$, let $p\in\partial B(w_1,r)\cap\partial B(w_2,r)$,
and for $i=1,2$, let $H_i^+= \{x\in\R^n:\,\langle x,p-w_i\rangle\geq \langle p,p-w_i\rangle\}$.
If $\angle (w_1,p,w_2)=\alpha$ satisfies $0<\alpha_0\leq\alpha\leq\alpha_1<\pi$
and $c_n=\frac{1}{2^{4n}n^n}$, then
$$
V\left(
{\rm conv}\left(B(w_1,r)\cup B(w_2,r)\right)\cap H_1^+\cap H_2^+ \right)\geq 
c_nr^n\alpha_0^{n+1}\cos\frac{\alpha_1}2.
$$
\end{claim}
\proof Let $\omega_{-1}=1$ and $\omega_k$ be the $k$ dimensional Lebesgue measure of $S^k$ for $k\geq 0$. 
First we verify that
\begin{equation}
\label{ballconvexhull0}
V\left({\rm conv}\left(B(w_1,r)\cup B(w_2,r)\right)\cap H_1^+\cap H_2^+ \right)\geq \widetilde{c}_n\cdot 
\frac{\left(1-\cos\frac{\alpha}2\right)^{\frac{n}2+1}\cos\frac{\alpha}2}{\alpha/2}\cdot r^n
\end{equation}
where $\widetilde{c}_2=1$ and $\widetilde{c}_n=\frac{2\omega_{n-3}}{n(n-1)(n-2)}$.
Let $w=\frac12(w_1+w_2)$, and let $u=\frac{p-w}{\|p-w\|}\in S^{n-1}$.
For $q=w+ru$, we observe that $p\in[w,q]$ and $\|q-p\|=r(1-\cos\frac{\alpha}2)$.
For $i=1,2$, we consider $v_i=w_i+ru\in \partial B(w_i,r)$ , and the point
$\widetilde{v}_i\in[v_1,v_2]$ with $\langle \widetilde{v}_i,p-w_i\rangle= \langle p,p-w_i\rangle$, and hence
$$
\|\widetilde{v}_i-q\|=\frac{\|q-p\|}{\tan\frac{\alpha}2}=\frac{1-\cos\frac{\alpha}2}{\tan\frac{\alpha}2}\cdot r.
$$
It follows that the area of the triangle $[p,\widetilde{v}_1,\widetilde{v}_2]$ is
$$
\frac{(1-\cos\frac{\alpha}2)^2}{\tan\frac{\alpha}2}\cdot r^2\geq 
\frac{(1-\cos\frac{\alpha}2)^2\cos\frac{\alpha}2}{\alpha/2}\cdot r^2.
$$
If $n=2$, then we have verified \eqref{ballconvexhull0}. 
If $n\geq 3$, then let $L$ be the affine $(n-2)$-space passing through $p$
and orthogonal to the linear two-space ${\rm lin}\{v_1,v_2\}$ containing all of our points defined so far.
Then $L$ intersects 
${\rm conv}\left(B(w_1,r)\cup B(w_2,r)\right)$ in an $(n-2)$-dimensional ball centered at $p$ and of
radius
$$
\sqrt{\|q-p\|\cdot(2r-\|q-p\|)}\geq \sqrt{r^2\left(1-\cos\frac{\alpha}2\right)},
$$
completing the proof of \eqref{ballconvexhull0}. 

Since $\left(\frac{1-\cos t}t\right)'=\frac{t\sin t+\cos t-1}{t^2}>0$ for $t\in\left(0,\frac{\pi}2\right)$, we deduce that
$\frac{\left(1-\cos\frac{\alpha}2\right)^{\frac{n}2+1}}{\alpha/2}$ is monotonically
increasing in $\alpha\in(0,\pi)$. However, $1-\cos t=\frac{(\sin t)^2}{1+\cos t}\geq \frac{t^2}8$
for $t\in\left(0,\frac{\pi}2\right)$, therefore \eqref{ballconvexhull0} yields
Claim~\ref{ballconvexhull} with $c'_n=\widetilde{c}_n/2^{\frac{5n}2+4}$ instead of $c_n$.

If $n=2$, then 
$c'_n>c_n$ and
Claim~\ref{ballconvexhull} readily holds.
If $n\geq 3$, then we use that
$\omega_{k-1}=\frac{k\pi^{k/2}}{\Gamma\left(\frac{k}2+1\right)}$
holds for $k\geq 1$ and that
 $\Gamma(t+1)<\left(\frac{t}{e}\right)^{t}\sqrt{2\pi(t+1)}$ holds for $t\geq \frac12$ according to Artin \cite{Art15} to conclude
\begin{eqnarray*}
c'_n&=&
\frac{2\omega_{n-3}}{2^{\frac{5n}2+4}n(n-1)(n-2)}=
\frac{\pi^{\frac{n}2-1}}{2^{\frac{5n}2+3}n(n-1)\Gamma\left(\frac{n}2\right)}>\frac{(2e\pi)^{\frac{n}2-1}}{2^{\frac{5n}2+3}n(n-1)(n-2)^{\frac{n}2-1}\sqrt{\pi n}}\\
&>&\left(\frac{e\pi}{2^4}\right)^{\frac{n}2}\frac1{2^4e\pi^{\frac32}n^n}>\frac1{2^{\frac{n}2+8}n^n}>\frac1{2^{4n}n^n}=c_n
\end{eqnarray*}
where $\frac{e\pi}{2^4}>\frac12$ and $e\pi^{\frac32}<2^4$.
 \hfill$\Box$\\

To compare a bounded portion of the hyperbolic or the spherical geometry to Euclidean ge\-o\-met\-ry,
we consider the map $\varphi:\R^{n+1}\backslash e^\bot\to e+e^\bot$ defined by
$$
\varphi(x)=\frac{x}{\langle x,e\rangle}.
$$
Writing $B^{n+1}$ to denote the unit Euclidean ball centred at the origin in $\R^{n+1}$
and $S^n_+=\{x\in S^n:\,\langle x,e\rangle>0\}$, we consider
the diffeomorphisms
\begin{eqnarray*}
\varphi_{H^n}&:\,&H^n\to e+(e^\bot\cap{\rm int}\, B^{n+1})\\
\varphi_{S^n}&:\,&S^n_+\to e+e^\bot
\end{eqnarray*}
that are the restrictions of $\varphi$
to $H^n$ and $S^n$, respectively. We observe that $\varphi_{H^n}$ and $\varphi_{S^n}$
map convex subsets of $H^n$ and $S^n_+$, respectively into convex subsets of $e+e^\bot$.
We observe that if $\mathcal{M}^n$ is either $H^n$ or $S^n$, then for any ball $B_{\mathcal{M}^n}(z,r)$ in $\mathcal{M}^n$ , there exists a hyperplane $\Pi$ in $\R^{n+1}$ such that
$\partial_{\mathcal{M}^n}B_{\mathcal{M}^n}(z,r)=\Pi\cap \mathcal{M}^n$. Since $\varphi|_\Pi$ maps lines into lines, 
it is a projective map $\Pi\to e+e^\bot$, therefore we have
\begin{eqnarray}
\label{ellipsoidH}
\varphi_{H^n}B_{H^n}(z,r)&&\mbox{ is an ellipsoid for $z\in H^n$ and $r>0$,}\\
\label{ellipsoidS}
\varphi_{S^n}B_{S^n}(z,r)&&\mbox{ is an ellipsoid for $z\in S^n_+$ and $r>0$ with 
$B_{S^n}(z,r)\subset S^n_+$.}
\end{eqnarray}
It follows from (\ref{geodesic}) that
for $z\in H^n$ with $d_{H^n}(e,z)=r$ or $z\in S^n_+$ with $d_{S^n}(e,z)=r<\frac{\pi}2$, if $u\in T_z$, then
\begin{eqnarray}
\label{differentialH}
\frac{\sqrt{-\mathcal{B}(u,u)}}{({\rm cosh}\,r)^2}&\leq\|D\varphi_{H^n}(z)u\|\leq &
\frac{\sqrt{-\mathcal{B}(u,u)}}{{\rm cosh}\,r}
\mbox{ \ \ and $|\det D\varphi_{H^n}(z)|=({\rm cosh}\,r)^{-(n+1)}$} \\
\label{differentialS}
\frac{\|u\|}{(\cos r)^2}&\geq\|D\varphi_{S^n}(z)u\|\geq &
\frac{\|u\|}{\cos r}\mbox{ \ \ and $|\det D\varphi_{S^n}(z)|=(\cos r)^{-(n+1)}$} 
\end{eqnarray}
where $D\varphi_{H^n}(z)$ and $D\varphi_{S^n}(z)$ stand for the differentials at $z$, which are linear maps
$T_z\to e^\bot$ (here $T_z$ is equipped with the scalar product $-\mathcal{B}(\cdot,\cdot)$ if $\mathcal{M}^n=H^n$).

If  $\mathcal{M}^n$ is either $H^n$ or $S^n$, then we observe that
$\partial_{\mathcal{M}^n} B_{\mathcal{M}^n}(z,r)$
is contained in a hyperplane of $\R^{n+1}$ for $z\in \mathcal{M}^n$
and $r>0$ (even $r\in(0,\pi)$ provided $\mathcal{M}^n=S^n$); namely
\begin{equation}
\label{ball-boundary}
\begin{array}{rcl}
\partial_{H^n} B_{H^n}(z,r)&=&\{x\in H^n: \mathcal{B}(x,z)= {\rm cosh}\,r\}\\
\partial_{S^n} B_{S^n}(z,r)&=&\{x\in S^n: \langle x,z\rangle= \cos r\}.
\end{array}
\end{equation}

\begin{lemma}
\label{ball-ellipsoid}
Let $u\in T_e=e^\bot$ with $\|u\|=1$ where $\mathcal{M}^n$ is either $H^n$ or $S^n$ (in particular, $\mathcal{B}(u,u)=-1$ provided $\mathcal{M}^n=H^n$).
\begin{description}
\item[(i)] If $\mathcal{M}^n=H^n$, $r>0$ and 
$z=e\,{\rm cosh}\,r+u\,{\rm sinh}\,r\in H^n$ (and hence $e\in\partial B_{H^n}(z,r)$), then $\varphi_{H^n}B_{H^n}(z,r)\subset e+T_e$ is an ellipsoid having $e$ on its boundary with exterior normal $-u$ such that one of its principal
axes is parallel to $u$ and is of length $a=\frac{{\rm sinh}\,2r}{{\rm cosh}\,2r}$, and all other principal axes are of length $b=\frac{2{\rm sinh}\,r}{\sqrt{{\rm cosh}\,2r}}$ where
$\frac{{\rm sinh}\,2r}{{\rm cosh}\,2r}<\frac{2{\rm sinh}\,r}{\sqrt{{\rm cosh}\,2r}}$.

\item[(ii)] If $\mathcal{M}^n=S^n$, $r\in\left(0,\frac{\pi}4\right)$ and 
$z_0=e\cos r+u \sin r\in S^n$ (and hence $e\in\partial B_{S^n}(z,r)$), then $\varphi_{S^n}B_{S^n}(z,r)\subset e+T_e$ is an ellipsoid having $e$ on its boundary with exterior normal $-u$ such that one of its principal
axes is parallel to $u$ and is of length $a=\frac{\sin 2r}{\cos 2r}$, and all other principal axes are of length 
$b=\frac{2\sin r}{\sqrt{\cos 2r}}$ where
$\frac{\sin 2r}{\cos 2r}>\frac{2\sin r}{\sqrt{\cos 2r}}$.
\end{description}
\end{lemma}
\proof In both cases, we consider an orthonormal basis $e_0,\ldots,e_n$ of $\R^{n+1}$ where $e_0=e$ and $e_1=u$, and for $x\in\R^{n+1}$, we write its coordinates with respect to $e_0,\ldots,e_n$. 
In particular, for $x=(x_0,x_1,\ldots,x_n)\in\R^{n+1}$ with $x_0\neq 0$, we have
$\varphi(x)=\left(1,\frac{x_1}{x_0},\ldots,\frac{x_n}{x_0}\right)$.
We  observe that the ellipsoid $\varphi_{\mathcal{M}^n}B_{\mathcal{M}^n}(z,r)$
in $e+T_e$ has an axial rotational symmetry around the line $e+\R u$.
It follows that $-u$ is an exterior normal to the ellipsoid at the boundary point $e$, and one of the principal axes of the ellipsoid is parallel to $u$.  

If $\mathcal{M}^n=H^n$, then $e\,{\rm cosh}\,2r+u\,{\rm sinh}\,2r\in H^n$ is the diametrically opposite point of $B_{H^n}(z,r)$ to $e$; therefore,
the principal axis of the ellipsoid $\varphi_{H^n}B_{H^n}(z,r)$ parallel to $u$ is of length 
$a=\frac{{\rm sinh}\,2r}{{\rm cosh}\,2r}$.

To determine the common length $b$ of the principal axes
of the ellipsoid $\varphi_{H^n}B_{H^n}(z,r)$ orthogonal to $u$, we note that according to \eqref{ball-boundary}, 
$x=(x_0,x_1,\ldots,x_n)\in\R^{n+1}$ satisfies $x\in\partial_{H^n}B_{H^n}(z,r)$ if and only if
\begin{equation}
\label{ball-boundaryH}
\begin{array}{rcl}
1+x_1^2+\ldots+x_n^2&=&x_0^2\\
x_0\cdot {\rm cosh}\,r-x_1\cdot {\rm sinh}\,r&=&{\rm cosh}\,r.
\end{array}
\end{equation}
We deduce from \eqref{ball-boundaryH} and the symmetries of
$\varphi_{H^n}B_{H^n}(z,r)$ that
\begin{eqnarray}
\nonumber
\left(\frac{b}2\right)^2&=&
\max\left\{\frac{x_2^2}{x_0^2}:\,
1+x_1^2+x_2^2=x_0^2\mbox{ \ and \ }
x_0\cdot {\rm cosh}\,r-x_1\cdot {\rm sinh}\,r={\rm cosh}\,r\right\}\\
\nonumber
&=&
\max\left\{\frac{x_0^2-(x_0-1)^2\,\frac{({\rm cosh}\,r)^2}{({\rm sinh}\,r)^2}-1}{x_0^2}:\,
x_0\geq 1\right\}\\
\label{bhalfH}
&=&
\max\left\{1-\frac{({\rm cosh}\,r)^2}{({\rm sinh}\,r)^2}+
\frac{2}{x_0}\cdot \frac{({\rm cosh}\,r)^2}{({\rm sinh}\,r)^2}
-\frac{1}{x_0^2}\cdot \frac{({\rm cosh}\,r)^2+({\rm sinh}\,r)^2}{({\rm sinh}\,r)^2}:\,
 x_0\geq 1\right\}.
\end{eqnarray}
For $f(t)=\frac{2}{t}\cdot \frac{({\rm cosh}\,r)^2}{({\rm sinh}\,r)^2}
-\frac{1}{t^2}\cdot \frac{({\rm cosh}\,r)^2+({\rm sinh}\,r)^2}{({\rm sinh}\,r)^2}$,
we have
$f'(t)=\frac{2}{t^3}\cdot \frac{1}{({\rm sinh}\,r)^2}
\left(({\rm cosh}\,r)^2+({\rm sinh}\,r)^2-t({\rm cosh}\,r)^2\right)$, and hence
$f$ has its maximum at $t=\frac{({\rm cosh}\,r)^2+({\rm sinh}\,r)^2}{({\rm cosh}\,r)^2}$.
Therefore, \eqref{bhalfH} and $({\rm cosh}\,r)^2-({\rm sinh}\,r)^2=1$ imply
\begin{eqnarray*}
\left(\frac{b}2\right)^2&=&\frac{-1}{({\rm sinh}\,r)^2}+
\frac{2({\rm cosh}\,r)^2}{({\rm cosh}\,r)^2+({\rm sinh}\,r)^2}\cdot 
\frac{({\rm cosh}\,r)^2}{({\rm sinh}\,r)^2}
-\frac{({\rm cosh}\,r)^4}{(({\rm cosh}\,r)^2+({\rm sinh}\,r)^2)({\rm sinh}\,r)^2}\\
&=&
\frac{({\rm cosh}\,r)^2(({\rm cosh}\,r)^2-1)-({\rm sinh}\,r)^2}
{(({\rm cosh}\,r)^2+({\rm sinh}\,r)^2)({\rm sinh}\,r)^2}=
\frac{({\rm sinh}\,r)^2}
{({\rm cosh}\,r)^2+({\rm sinh}\,r)^2}.
\end{eqnarray*}
Since ${\rm cosh}\,2r=({\rm cosh}\,r)^2+({\rm sinh}\,r)^2$, we have
$b=\frac{2\,{\rm sinh}\,r}{\sqrt{{\rm cosh}\,2r}}$.

To compare $a$ and $b$, using ${\rm sinh}\,2r=2({\rm sinh}\,r)({\rm cosh}\,r)$ yields
$\frac{a}{b}=\frac{{\rm cosh}\,r}{\sqrt{{\rm cosh}\,2r}}=
\frac{{\rm cosh}\,r}{\sqrt{({\rm cosh}\,r)^2+({\rm sinh}\,r)^2}}<1$.

Let us sketch the argument, analogously to the one above, in the spherical case $\mathcal{M}^n=H^n$ and $0<r<\frac{\pi}4$.
In particular, $e\cos 2r+u\sin 2r\in S^n$ is the diametrically opposite point of $B_{S^n}(z,r)$ to $e$, and hence 
the principal axis of the ellipsoid $\varphi_{S^n}B_{S^n}(z,r)$ parallel to $u$ is of length $a=\frac{\sin 2r}{\cos 2r}$.

To determine the common length $b$ of the principal axes
of the ellipsoid $\varphi_{S^n}B_{S^n}(z,r)$ orthogonal to $u$, we note that according to \eqref{ball-boundary}, 
$x=(x_0,x_1,\ldots,x_n)\in\R^{n+1}$ satisfies $x\in\partial_{S^n}B_{S^n}(z,r)$ if and only if
\begin{equation}
\label{ball-boundaryS}
\begin{array}{rcl}
x_0^2+x_1^2+\ldots+x_n^2&=&1\\
x_0\cdot \cos r+x_1\cdot \sin r&=&\cos r.
\end{array}
\end{equation}
Using $(\cos r)^2+(\sin r)^2=1$ and similar argument as in the hyperbolic case, we deduce that
$$
b=\frac{2\sin r}{\sqrt{(\cos r)^2-(\sin r)^2}}=\frac{2\sin r}{\sqrt{\cos 2r}}.
$$

To compare $a$ and $b$, using $\sin 2r=2\sin r\cos r$ yields
$\frac{a}{b}=\frac{\cos r}{\sqrt{\cos 2r}}=
\frac{\cos r}{\sqrt{(\cos r)^2-(\sin r)^2}}>1$.
\hfill\makebox{$\Box$}\\

\begin{lemma}
\label{pythegorean}
If $[y_0,y_1,p]$ is a triangle in $\mathcal{M}^n$
where
$\mathcal{M}^n$ is either $\R^n$,  $H^n$ or $S^n$, and 
$\angle(y_0,y_1,p)\leq\frac{\pi}2$,
$d_{\mathcal{M}^n}(y_0,y_1)=R$
and $d_{\mathcal{M}^n}(y_0,p)\geq R+\eta$
for $R>0$ and $\eta\geq 0$, then
$$
d_{\mathcal{M}^n}(y_1,p)\geq
\left\{
\begin{array}{ll}
\sqrt{2R\cdot\eta}&\mbox{ \ if $\mathcal{M}^n=\R^n$,}\\[1ex]
\sqrt{\tanh R\cdot\eta}&\mbox{ \ if $\mathcal{M}^n=H^n$,}\\[1ex]
\sqrt{\tan R\cdot\eta}&\mbox{ \ if $\mathcal{M}^n=S^n$ and $R,\eta\leq\frac{\pi}4$}.
\end{array} \right.
$$
\end{lemma}
\proof
We use $\angle(y_0,y_1,p)\leq\frac{\pi}2$ and
the Law of Cosines in $\mathcal{M}^n$.

If $\mathcal{M}^n=\R^n$, then
$$
d_{\R^n}(y_1,p)\geq \sqrt{(R+\eta)^2-R^2}\geq \sqrt{2R\eta}.
$$

If $\mathcal{M}^n=H^n$, then
$$
\cosh{d_{H^n}(y_1,p)}\geq \frac{\cosh(R+\eta)}{\cosh{R}}
=\cosh{\eta}+\tanh{R}\cdot \sinh{\eta}.
$$
Therefore, Lemma~\ref{pythegorean} follows 
in the case  $\mathcal{M}^n=H^n$ if
\begin{equation}
\label{Hpyth}
\cosh{\left(\sqrt{\tanh{R}}\cdot\sqrt{\eta}\right)}\leq \cosh{\eta}+\tanh{R}\cdot \sinh{\eta}.
\end{equation}
The equation \eqref{Hpyth} readily holds if $\eta=0$. Thus using differentiation,
it is sufficient to verify
\begin{equation}
\label{Hpythder}
\sinh{\left(\sqrt{\tanh{R}}\cdot\sqrt{\eta}\right)}\cdot\frac{\sqrt{\tanh{R}}}{2\sqrt{\eta}}\leq \sinh{\eta}+\tanh{R}\cdot\cosh{\eta}
\end{equation}
for $\eta>0$.
If $\eta\in(0,\tanh{R}]$, then we use that
$\sinh{t}/t$ is increasing for $t>0$, $\tanh{R}<1$
and $\sinh{1}<2$, and hence
$$
\sinh{\left(\sqrt{\tanh{R}}\cdot\sqrt{\eta}\right)}\cdot\frac{\sqrt{\tanh{R}}}{2\sqrt{\eta}}=\frac{\sinh{\left(\sqrt{\tanh{R}}\cdot\sqrt{\eta}\right)}}{\sqrt{\tanh{R}}\cdot\sqrt{\eta}}\cdot\frac{\tanh{R}}2
<\tanh{R}<\tanh{R}\cdot\cosh{\eta}.
$$
On the other hand, if $\eta\geq\tanh{R}$, then
$$
\sinh{\left(\sqrt{\tanh{R}}\cdot\sqrt{\eta}\right)}\cdot\frac{\sqrt{\tanh{R}}}{2\sqrt{\eta}}< \sinh{\eta},
$$
proving \eqref{Hpythder}, and in turn \eqref{Hpyth}.

If $\mathcal{M}^n=S^n$, then we use a similar argument as in the hyperbolic case; in particular, we have
$$
\cos{d_{S^n}(y_1,p)}\leq \frac{\cos(R+\eta)}{\cos{R}}
=\cos{\eta}-\tan{R}\cdot\sin{\eta}.
$$
Therefore, Lemma~\ref{pythegorean} follows 
in the case  $\mathcal{M}^n=S^n$ if
\begin{equation}
\label{Spyth}
\cos{\left(\sqrt{\tan{R}}\cdot\sqrt{\eta}\right)}\geq \cos{\eta}-\tan{R}\cdot \sin{\eta}.
\end{equation}
The equation \eqref{Spyth} readily holds if $\eta=0$. Thus using differentiation,
it is sufficient to verify
\begin{equation}
\label{Spythder}
\sin{\left(\sqrt{\tan{R}}\cdot\sqrt{\eta}\right)}\cdot
\frac{\sqrt{\tan{R}}}{2\sqrt{\eta}}\leq \sin{\eta}+\tan{R}\cdot\cos{\eta}
\end{equation}
for $\eta>0$.
If $\eta\leq\tan{R}$ where $\tan{R}\leq 1$ because
of $R\leq\frac{\pi}4$, then we use that
$\sin{t}< t$ for $t>0$
and $\cos{1}>\cos\frac{\pi}3=\frac12$, and hence
$$
\sin{\left(\sqrt{\tan{R}}\cdot\sqrt{\eta}\right)}\cdot
\frac{\sqrt{\tan{R}}}{2\sqrt{\eta}}\leq\frac{\tan{R}}2
<\tan{R}\cdot\cos{\eta}.
$$
On the other hand, if $\eta\geq\tan{R}$, then
$$
\sin{\left(\sqrt{\tan{R}}\cdot\sqrt{\eta}\right)}\cdot
\frac{\sqrt{\tan{R}}}{2\sqrt{\eta}}< \sin{\eta},
$$
proving \eqref{Spythder}, and in turn \eqref{Spyth}.
\hfill\makebox{$\Box$}\\

\begin{prop}
\label{Improve}
Let $\mathcal{M}^n$ be either $\R^n$,  $H^n$ or $S^n$, and let
$X\subset \mathcal{M}^n$ be a convex body of constant width $D$ for $D>0$ where $D<\frac{\pi}2$  if $\mathcal{M}^n=S^n$. In addition, let $x_1,x_2\in X$ satisfy $d_{\mathcal{M}^n}(x_1,x_2)=D$, let $y_0$
be the midpoint of $[x_1,x_2]_{\mathcal{M}^n}$, and let 
$$
\varrho=\left\{
\begin{array}{ll}
\frac{D}2&\mbox{ \ if $\mathcal{M}^n=\R^n$ or $\mathcal{M}^n=H^n$,}\\[1ex]
\min\left\{\frac{D}{2},\frac{\pi}{8}-\frac{D}{4}\right\}&\mbox{ \ if $\mathcal{M}^n=S^n$ and $D<\frac{\pi}{2}$}.
\end{array} \right.
$$

If there exists $z\in\partial_{\mathcal{M}^n} X$ such that $d_{\mathcal{M}^n}(z,y_0)\geq \frac{D}2+\eta$ for $\eta\in(0,\eta_0)$,
 then
there exists a hyperplane $H$ such that
$$
V_{\mathcal{M}^n}\left(\mbox{\rm conv}\,\tau_H X^{(\varrho)}\right)
-V_{\mathcal{M}^n}(\tau_H X^{(\varrho)})\geq \widetilde{\gamma}_1 \eta^{\frac{3n+2}2}
$$
where
$$
\widetilde{\gamma}_1=
\left\{
\begin{array}{ll}
\frac{1}{2^{12n}n^n}\cdot
\frac{1}{D^{\frac{n+2}2}}&\mbox{ \ if $\mathcal{M}^n=\R^n$ or $\mathcal{M}^n=S^n$ and $D<\frac{\pi}2$;}\\[1ex]
\frac{1}{2^{8n}n^n}\cdot
\frac{1}{(\sinh{10D})^{\frac{n+2}2}}&\mbox{ \ if $\mathcal{M}^n=H^n$;}
\end{array} \right.
$$
$$
\eta_0=
\left\{
\begin{array}{ll}
\frac{D}2&\mbox{ \ if $\mathcal{M}^n=\R^n$ or $\mathcal{M}^n=S^n$ and $D\leq\frac{\pi}6$;}\\[1ex]
\frac{3}{4\pi}\left(\frac{\pi}2-D\right)^2&\mbox{ \ if $\mathcal{M}^n=S^n$ and $\frac{\pi}6\leq D<\frac{\pi}2$;}
\\[1ex]
\min\left\{1,\frac{D}2\right\}&\mbox{ \ if $\mathcal{M}^n=H^n$.}
\end{array} \right.
$$
\end{prop}
\proof We may  assume that $z$ is the farthest point of $X$ from $y_0$, and hence
\begin{equation}
\label{XinBally0}
X\subset B_{\mathcal{M}^n}(y_0, d_{\mathcal{M}^n}(z,y_0)).
\end{equation}

Let $Y=X^{(\varrho)}$ that is a convex body of diameter $D+2\varrho$ according to Lemma~\ref{completeparallel}. There exist $y_1,y_2\in \partial_{\mathcal{M}^n} Y$ such that
$x_1,x_2\in[y_1,y_2]_{\mathcal{M}^n}$ and $d_{\mathcal{M}^n}(x_i,y_i)=\varrho$, $i=1,2$.
It follows from (\ref{XinBally0}) that there exists
 $p\in \partial_{\mathcal{M}^n} Y$ such that
$z\in[p,y_0]_{\mathcal{M}^n}$, and hence 
$-N_Y(p)\in T_p$ points towards $y_0$ along $[p,y_0]_{\mathcal{M}^n}$ and
$d_{\mathcal{M}^n}(p,y_0)\geq \frac{D}2+\varrho+\eta$. We may assume that $\angle(p,y_0,y_1)\leq\frac{\pi}2$ after possibly interchanching $y_1$ and $y_2$.

Let $H$ be the perpendicular bisector hyperplane  in $\mathcal{M}^n$ of the segment $[p,y_1]_{\mathcal{M}^n}$,
let $H^+$ be the corresponding half-space containing $p$, and
let $\widetilde{x}=\sigma_Hx_1$. We observe that $p\in B_{\mathcal{M}^n}(\widetilde{x},\varrho)\subset \sigma_HY$
and $p\in B_{\mathcal{M}^n}(z,\varrho)\subset Y$.

It follows from applying Lemma~\ref{ballboundary}
with $R=\frac{D}2+\varrho\leq D$ that 
for $N_Y(p)\in T_p$ and $\sigma_HN_Y(y_1)=N_{\sigma_HY}(p)\in T_p$, we have
\begin{equation}
\label{angleatplow}
\angle(p,y_1,y_0)-\angle(y_1,p,y_0)=
\angle(N_{\sigma_HY}(p),N_Y(p))\geq \widetilde{\gamma}\cdot \eta
\end{equation}
where 
\begin{equation}
\label{gamma*73}
\widetilde{\gamma}=\left\{
\begin{array}{cl}
  \frac1{4D}   &\mbox{ if $\mathcal{M}^n=\R^n$ or $\mathcal{M}^n=S^n$;}  \\[1ex]
   \frac1{\sinh{5D}}  & \mbox{ if $\mathcal{M}^n=H^n$.}
\end{array}\right.
\end{equation}

Since $d_{\mathcal{M}^n}(p,y_0)\geq \frac{D}2+\varrho+\eta=R+\eta$ and 
$\angle(y_0,y_1,p)\leq\frac{\pi}2$, we deduce from
Lemma~\ref{pythegorean}
applied with $R=\frac{D}2+\varrho=D$ if
$\mathcal{M}^n=\R^n$ or $\mathcal{M}^n=H^n$,
and with $R=\frac{D}2+\varrho>\frac{D}2$ if
$\mathcal{M}^n=S^n$
that
$$
d_{\mathcal{M}^n}(y_1,p)\geq 4\widetilde{\gamma}_0\cdot \sqrt{\eta}
$$
where 
\begin{equation}
\label{tilde-gamma0}    
\widetilde{\gamma}_0=\left\{
\begin{array}{ll}
\frac{\sqrt{D}}4&\mbox{ \ if $\mathcal{M}^n=\R^n$;}\\[1ex]
\frac{\sqrt{\tanh D}}4&\mbox{ \ if $\mathcal{M}^n=H^n$;}\\[1ex]
\frac{\sqrt{\tan{\frac{D}2}}}4&\mbox{ \ if $\mathcal{M}^n=S^n$;}
\end{array} \right.
\end{equation}
and hence $\eta$  satisfies
$$
\widetilde{\gamma}_0\sqrt{\eta}< \varrho
\mbox{ \ and \ }\eta<\varrho.
$$
In turn, we deduce that
\begin{equation}
\label{ballinH+}
B_{\mathcal{M}^n}\left(p,2\widetilde{\gamma}_0\sqrt{\eta}\right)\subset H^+.
\end{equation}
Let $\widetilde{x}_\eta\in[\widetilde{x},p]_{\mathcal{M}^n}$ and $z_\eta\in[z,p]_{\mathcal{M}^n}$ such that
$d_{\mathcal{M}^n}(p,\widetilde{x}_\eta)=\widetilde{\gamma}_0\sqrt{\eta}$ and $d_{\mathcal{M}^n}(p,z_\eta)=\widetilde{\gamma}_0\sqrt{\eta}$.
It follows that $N_{\sigma_HY}(p)$ is the exterior unit normal to 
$B_{\mathcal{M}^n}\left(\widetilde{x}_\eta,\widetilde{\gamma}_0\sqrt{\eta}\right)$ at $p$ and
$N_{Y}(p)$ is the exterior unit normal to 
$B_{\mathcal{M}^n}\left(z_\eta,\widetilde{\gamma}_0\sqrt{\eta}\right)$ at $p$; moreover,  \eqref{ballinH+} yields  that
\begin{eqnarray}
\nonumber
B_{\mathcal{M}^n}\left(\widetilde{x}_\eta,\widetilde{\gamma}_0\sqrt{\eta}\right)\cup
B_{\mathcal{M}^n}\left(z_\eta,\widetilde{\gamma}_0\sqrt{\eta}\right)
&\subset&\left(B_{\mathcal{M}^n}(\widetilde{x},\varrho)\cap B_{\mathcal{M}^n}\left(p,2\widetilde{\gamma}_0\sqrt{\eta}\right)\right)\cup
\left(B_{\mathcal{M}^n}(z,\varrho)\cap B_{\mathcal{M}^n}\left(p,2\widetilde{\gamma}_0\sqrt{\eta}\right)\right)\\
\nonumber
&\subset&(B_{\mathcal{M}^n}(\widetilde{x},\varrho)\cap H^+)\cup(B_{\mathcal{M}^n}(z,\varrho)\cap H^+)\\
\label{ballsintauY}
&\subset&  (\sigma_H Y\cap H^+)\cup (Y\cap H^+)\subset\tau_H Y.
\end{eqnarray}

Let $H^+_x$ be the half-space of $\mathcal{M}^n$ touching 
$B_{\mathcal{M}^n}\left(\widetilde{x}_\eta,\widetilde{\gamma}_0\sqrt{\eta}\right)$ at $p$ with $H^+_x\cap B_{\mathcal{M}^n}\left(\widetilde{x}_\eta,\widetilde{\gamma}_0\sqrt{\eta}\right)=\{p\}$,
and hence  ${\rm int}\,H^+_x\cap \sigma_HY=\emptyset$. Similarly, let $H^+_z$ be the half-space of $\mathcal{M}^n$ touching 
$B_{\mathcal{M}^n}\left(z_\eta,\widetilde{\gamma}_0\sqrt{\eta}\right)$ at $p$ with $H^+_z\cap B_{\mathcal{M}^n}\left(z_\eta,\widetilde{\gamma}_0\sqrt{\eta}\right)=\{p\}$,
and hence  ${\rm int}\,H^+_z\cap Y=\emptyset$.
Therefore
if either $\mathcal{M}^n=S^n$ and $D<\frac{\pi}2$, or $\mathcal{M}^n=\R^n$ or $\mathcal{M}^n=H^n$, then 
\begin{equation}
\label{difference-convex}
{\rm conv}\left\{B_{\mathcal{M}^n}\left(\widetilde{x}_\eta,\widetilde{\gamma}_0\sqrt{\eta}\right),
B_{\mathcal{M}^n}\left(z_\eta,\widetilde{\gamma}_0\sqrt{\eta}\right)\right\}
\cap ({\rm int}\,H^+_x)\cap ({\rm int}\,H^+_z)\subset
\left(\mbox{\rm conv}\tau_H X^{(\varrho)}\right)
\backslash\tau_H X^{(\varrho)}.
\end{equation}
In the triangle $[y_1,p,y_0]_{\mathcal{M}^n}$, we have $d_{\mathcal{M}^n}(p,y_0)\geq D+\varrho=d_{\mathcal{M}^n}(y_1,y_0)$, thus $\angle(y_1,p,y_0)<\angle(p,y_1,y_0)$, and hence $\angle\left(y_1,p,y_0\right)<\frac{\pi}{2}$ because either $\mathcal{M}^n=S^n$ and all sides of $[y_1,p,y_0]_{S^n}$ are acute, or $\mathcal{M}^n=\R^n$ or $\mathcal{M}^n=H^n$. We deduce from \eqref{angleatplow} and $\angle(y_1,p,y_0)<\frac{\pi}2$ that
\begin{equation}
\label{angleatpup}
\angle(N_{\sigma_HY}(p),N_Y(p))<\frac{\pi}2.
\end{equation}
The rest of the argument is divided into three cases depending on $\mathcal{M}^n$.\\

\noindent{\bf Case 1} $\mathcal{M}^n=\R^n$

In this case, \eqref{angleatplow} reads as
$$
\angle(\widetilde{x}_\eta,p,z_\eta)=\angle(N_{\sigma_HY}(p),N_Y(p))\geq \frac{\eta}{4D}.
$$ 
Since $\angle(N_{\sigma_HY}(p),N_Y(p))<\frac{\pi}2$ according to \eqref{angleatpup},
we deduce from \eqref{tilde-gamma0}, \eqref{difference-convex}  and
Claim~\ref{ballconvexhull} (where $c_n=\frac1{2^{4n}n^n}$) that
$$
V\left(\left(\mbox{\rm conv}\,\tau_H X^{(\varrho)}\right)
\backslash\tau_H X^{(\varrho)}\right)\geq
\frac{c_n}{\sqrt{2}}\cdot\left(\frac{\sqrt{D}\cdot \sqrt{\eta}}4\right)^n\cdot\left(\frac{\eta}{4D}\right)^{n+1}=
\frac{c_n}{4^{2n+1}\sqrt{2}}\cdot
\frac{\eta^{\frac{3n}2+1}}{D^{\frac{n}2+1}}>
\frac{1}{2^{9n}n^n}\cdot
\frac{\eta^{\frac{3n}2+1}}{D^{\frac{n}2+1}}.
$$

\noindent{\bf Case 2} $\mathcal{M}^n=H^n$

In this case, we may assume that
$$
p=e.
$$
We consider the ellipsoids $E_x=\varphi_{H^n}B_{H^n}(\widetilde{x}_\eta,\widetilde{\gamma}_0\sqrt{\eta})\subset e+T_e$
and $E_z=\varphi_{H^n}B_{H^n}(z_\eta,\widetilde{\gamma}_0\sqrt{\eta})\subset e+T_e$ satisfying
that $N_{\sigma_HY}(p)$ is an exterior normal to $E_x$ at $p=e$,
and $N_{Y}(p)$ is an exterior normal to $E_z$ at $p=e$. In addition,
$N_{\sigma_HY}(p)$ and $N_{Y}(p)$ are interior normals for the half-spaces
$$
\widetilde{H}_x^+=\varphi_{H^n}H_x^+\mbox{ \ and \ }\widetilde{H}_z^+=\varphi_{H^n}H_z^+
$$
in $e+T_e$, respectively. We deduce from \eqref{difference-convex} that
\begin{equation}
\label{difference-convexHE}
{\rm conv}_{e+T_e}(E_x\cup E_z)
\cap ({\rm int}\,\widetilde{H}^+_x)\cap ({\rm int}\,\widetilde{H}^+_z)\subset
\varphi_{H^n}\left(
\left(\mbox{\rm conv}_{H^n}\,\tau_H X^{(\varrho)}\right)
\backslash\tau_H X^{(\varrho)}\right).
\end{equation}

According to Lemma~\ref{ball-ellipsoid}, $E_z$ has axial rotational symmetry around
$p+\R\,N_{Y}(p)$, and the axis of $E_z$ contained in $p+\R\,N_{Y}(p)$ is of length
$\tanh 2\widetilde{\gamma}_0\sqrt{\eta}=\frac{\sinh 2\widetilde{\gamma}_0\sqrt{\eta}}{\cosh 2\widetilde{\gamma}_0\sqrt{\eta}}<
\frac{2\sinh\widetilde{\gamma}_0\sqrt{\eta}}{\sqrt{\cosh 2\widetilde{\gamma}_0\sqrt{\eta}}}$
where $\frac{2\sinh\widetilde{\gamma}_0\sqrt{\eta}}{\sqrt{\cosh 2\widetilde{\gamma}_0\sqrt{\eta}}}$ is the length of the orthogonal axes.
Therefore setting
$$
z^*=p-N_{Y}(p)\cdot \frac{\tanh{2\widetilde{\gamma}_0\sqrt{\eta}}}2,
$$  
using $\widetilde{B}(\cdot,\cdot)$ to denote $n$-dimensional balls in $e+T_e$,
we have
$$
\widetilde{B}\left(z^*, \frac{\tanh{2\widetilde{\gamma}_0\sqrt{\eta}}}2\right)
\subset E_z,
$$
and has $N_{Y}(p)$ as exterior unit normal at $p$. Similarly,
$$
x^*=p-N_{\sigma_HY}(p)\cdot \frac{\tanh{2\widetilde{\gamma}_0\sqrt{\eta}}}2,
$$  
satisfies that
$$
\widetilde{B}\left(x^*, \frac{\tanh{2\widetilde{\gamma}_0\sqrt{\eta}}}2
\right)\subset E_x,
$$
and has $N_{\sigma_HY}(p)$ as exterior unit normal at $p$. We deduce from \eqref{difference-convexHE} that
\begin{equation}
\label{difference-convexH}
{\rm conv}_{e+T_e}\left\{\widetilde{B}\left(x^*, \frac{\tanh{2\widetilde{\gamma}_0\sqrt{\eta}}}2\right),
\widetilde{B}\left(z^*, \frac{\tanh{2\widetilde{\gamma}_0\sqrt{\eta}}}2\right) \right\}
\cap\left({\rm int}\,\widetilde{H}^+_x\right)\cap\left({\rm int}\,\widetilde{H}^+_z\right)\subset
\end{equation}
$$
\subset \varphi_{H^n}\left(
\left(\mbox{\rm conv}_{H^n}\,\tau_H X^{(\varrho)}\right)
\backslash\tau_H X^{(\varrho)}\right).
$$
It follows from \eqref{angleatplow} and Lemma~\ref{ballboundary} that for $R=\frac{D}2+\varrho\leq D$, we have
\begin{equation}
\label{anglecentersH}
\angle(\widetilde{x}^*,p,z^*)=\angle(N_{\sigma_HY}(p),N_Y(p))>\frac{\eta}{\sinh 5D}.
\end{equation} 
Since $\angle(N_{\sigma_HY}(p),N_Y(p))<\frac{\pi}2$ according to \eqref{angleatpup},
and \eqref{differentialH} yields that
$$
V_{H^n}(X)\geq V_{\R^n}(\varphi_{H^n}X)
$$
for any bounded Borel set $X\subset H^n$,
we deduce from \eqref{difference-convexH}, \eqref{anglecentersH}  and
Claim~\ref{ballconvexhull} that
\begin{equation}
\label{Hbasicest}    
V_{H^n}\left(\left(\mbox{\rm conv}_{H^n}\,\tau_H X^{(\varrho)}\right)
\backslash\tau_H X^{(\varrho)}\right)\geq
\frac{c_n}{\sqrt{2}}\cdot
\left(\frac{\tanh{2\widetilde{\gamma}_0\sqrt{\eta}}}2\right)^n\cdot\left(\frac{\eta}{\sinh 5D}\right)^{n+1}.
\end{equation}
Here $\widetilde{\gamma}_0=\frac{\sqrt{\tanh{D}}}4<\frac14$ and $\eta\in(0,1)$ yield that $$
\frac{\tanh{2\widetilde{\gamma}_0\sqrt{\eta}}}2>0.9\cdot \frac{\sqrt{\eta\tanh{D}}}4.
$$
On the other hand, $\tanh{D}>\frac15\tanh{5D}$ and
$\sinh{10D}=2\sinh{5D}\cdot \cosh{5D}$ imply
\begin{eqnarray*}
\frac{(\tanh{D})^{\frac{n}2}}{(\sinh{5D})^{n+1}}&>& \frac{\left(\frac15\tanh{5D}\right)^{n/2}}{(\sinh{5D})^{n+1}}=
\left(\frac25\right)^{\frac{n}2}\frac1{(\sinh{10D})^{\frac{n}2}\sinh{5D}}\\
&>&
\left(\frac25\right)^{\frac{n}2}\frac1{(\sinh{10D})^{\frac{n}2+1}}.
\end{eqnarray*}
Using these last two estimates and $0.9\sqrt{\frac25}>\frac12$, we deduce from
\eqref{Hbasicest} that
$$
V\left(\left(\mbox{\rm conv}\,\tau_H X^{(\varrho)}\right)
\backslash\tau_H X^{(\varrho)}\right)\geq
\frac{c_n}{4^n\sqrt{2}}\left(0.9\sqrt{\frac25}\right)^n\cdot
\frac{\eta^{\frac{3n}2+1}}{(\sinh{10D})^{\frac{n}2+1}}
>\frac1{2^{8n}n^n}
\cdot
\frac{\eta^{\frac{3n}2+1}}{(\sinh{10D})^{\frac{n}2+1}}.
$$

\noindent{\bf Case 3} $\mathcal{M}^n=S^n$

We may assume once more that $p=e$. Considering the ellipsoids $E_x=\varphi_{S^n}B_{S^n}\left(\widetilde{x}_{\eta},\widetilde{\gamma}_0\sqrt{\eta}\right)\subset e+T_e$ and $E_z=\varphi_{S^n}B_{S^n}\left(z_{\eta},\widetilde{\gamma}_0\sqrt{\eta}\right)\subset e+T_e$, we have $N_{\sigma_HY}\left(p\right)$ as an exterior normal to $E_x$ at $p$ and $N_Y\left(p\right)$ an exterior normal to $E_z$ at $p$. Then again $N_{\sigma_HY}(p)$ and $N_{Y}(p)$ are interior normals for the half-spaces
$$
\widetilde{H}_x^+=\varphi_{S^n}H_x^+\mbox{ \ and \ }\widetilde{H}_z^+=\varphi_{S^n}H_z^+
$$
in $e+T_e$, respectively.  We deduce from \eqref{difference-convex} that
\begin{equation}
\label{difference-convexSE}
{\rm conv}_{e+T_e}(E_x\cup E_z)
\cap ({\rm int}\,\widetilde{H}^+_x)\cap ({\rm int}\,\widetilde{H}^+_z)\subset
\varphi_{S^n}\left(
\left(\mbox{\rm conv}_{H^n}\,\tau_H X^{(\varrho)}\right)
\backslash\tau_H X^{(\varrho)}\right).
\end{equation}

Lemma~\ref{ball-ellipsoid} implies that $E_z$ has axial rotational symmetry around
$p+\R\,N_{Y}(p)$, and the axis of $E_z$ contained in $p+\R\,N_{Y}(p)$ is of length
$\tan\left(2\widetilde{\gamma}_0\sqrt{\eta}\right)>\frac{2\sin\left(\widetilde{\gamma}_0\sqrt{\eta}\right)}{\sqrt{\cos\left(2\widetilde{\gamma}_0\sqrt{\eta}\right)}}$
where $\frac{2\sin\left(\widetilde{\gamma}_0\sqrt{\eta}\right)}{\sqrt{\cos\left(2\widetilde{\gamma}_0\sqrt{\eta}\right)}}$ is the length of the orthogonal axes.
Since if $a=\tan\left(2\widetilde{\gamma}_0\sqrt{\eta}\right)$ is the major axis and $b=\frac{2\sin\left(\widetilde{\gamma}_0\sqrt{\eta}\right)}{\sqrt{\cos\left(2\widetilde{\gamma}_0\sqrt{\eta}\right)}}$
is the minor axis of an ellipse, then its minimal radius of curvature is
$\frac{b^2}{2a}=\tan{\widetilde{\gamma}_0\sqrt{\eta}}$, we deduce
for
\begin{eqnarray*}
z^*&=&p-N_{Y}(p)\cdot\tan{\widetilde{\gamma}_0\sqrt{\eta}},\\
x^*&=&p-N_{\sigma_HY}(p)\cdot\tan{\widetilde{\gamma}_0\sqrt{\eta}},
\end{eqnarray*} 
that
$$
\widetilde{B}\left(z^*,\tan{\widetilde{\gamma}_0\sqrt{\eta}}\right)\subset E_z,
$$
and has $N_{Y}(p)$ as exterior unit normal at $p$. Similarly,
$$
\widetilde{B}\left(x^*,\tan{\widetilde{\gamma}_0\sqrt{\eta}}\right)\subset E_x,
$$
and has $N_{\sigma_HY}(p)$ as exterior unit normal at $p$. We deduce from \eqref{difference-convexSE} that
\begin{equation}
\label{difference-convexS}
{\rm conv}_{e+T_e}\left\{\widetilde{B}\left(x^*,\tan{\widetilde{\gamma}_0\sqrt{\eta}}\right),
\widetilde{B}\left(z^*,\tan{\widetilde{\gamma}_0\sqrt{\eta}}\right) \right\}
\cap ({\rm int}\,\widetilde{H}^+_x)\cap ({\rm int}\,\widetilde{H}^+_z)\subset
\end{equation}
$$
\subset \varphi_{S^n}\left(\left(\mbox{\rm conv}_{H^n}\,\tau_H X^{(\varrho)}\right)\backslash\tau_H X^{(\varrho)}\right).
$$
From \eqref{angleatplow} and Lemma~\ref{ballboundary} setting $R=\frac{D}2+\varrho\leq D$, we have
\begin{equation}
\label{anglecentersS}
\angle(\widetilde{x}^*,p,z^*)=\angle(N_{\sigma_HY}(p),N_Y(p))>\frac{\eta}{4D}.
\end{equation} 
Since \eqref{differentialS} yields that
$$
V_{S^n}(X)\geq (\cos r)^{n+1}\cdot
V_{\R^n}(\varphi_{S^n}X)
$$
for any bounded Borel set $X\subset B_{S^n}(e,r)$,
and
\begin{equation}
\label{withinpi4}    
{\rm conv}\Big\{B_{S^n}\left(\widetilde{x}_\eta,\widetilde{\gamma}_0\sqrt{\eta}\right),
B_{S^n}\left(z_\eta,\widetilde{\gamma}_0\sqrt{\eta}\right)\Big\}\subset
B_{S^n}(p,2\varrho)\subset B_{S^n}\left(e,\frac{\pi}4\right),
\end{equation}
we deduce from \eqref{difference-convexS}, \eqref{anglecentersS}  and
Claim~\ref{ballconvexhull} that
$$
V_{\mathbb{R}^n}\left(\varphi_{S^n}\left(\left(\mbox{\rm conv}_{S^n}\,\tau_H X^{(\varrho)}\right)
\backslash\tau_H X^{(\varrho)}\right)\right)\geq
\frac1{2^{\frac{n+1}2}}\cdot\frac{c_n}{\sqrt{2}}\cdot
\left(\tan\left(\widetilde{\gamma}_0\sqrt{\eta}\right)\right)^n\cdot\left(\frac{\eta}{4D}\right)^{n+1}.
$$
Since $\eta<\frac{D}{2}<\frac{\pi}{4}$, we have
$$
0<\widetilde{\gamma}_0\sqrt{\eta}<\frac{\sqrt{\eta}}{4}<\frac{\pi}{2}
$$
implying
$$
\tan\left(\widetilde{\gamma}_0\sqrt{\eta}\right)\geq\frac{\sqrt{\tan\left(\frac{D}{2}\right)\eta}}{4}\geq\frac{\sqrt{\frac{D}{2}\cdot\eta}}{4}.
$$
This gives us
$$
V_{\mathbb{R}^n}\left(\varphi_{S^n}\left(\left(\mbox{\rm conv}_{S^n}\,\tau_H X^{(\varrho)}\right)
\backslash\tau_H X^{(\varrho)}\right)\right)\geq
\frac{\eta^{\frac{3n}{2}+1}}{D^{\frac{n+2}{2}}\cdot 2^{9n+3}\cdot n^n},
$$
concluding via \eqref{withinpi4} that
$$
V_{S^n}\left(\left(\mbox{\rm conv}_{S^n}\,\tau_H X^{(\varrho)}\right)
\backslash\tau_H X^{(\varrho)}\right)\geq\frac{\eta^{\frac{3n}{2}+1}\cdot\left(\cos\left(\frac{\pi}{4}\right)\right)^{n+1}}{D^{\frac{n+2}{2}}\cdot 2^{9n+3}\cdot n^n}=\frac{\eta^{\frac{3n}{2}+1}}{D^{\frac{n+2}{2}}\cdot 2^{\frac{19n+7}{2}}\cdot n^n}>
\frac{\eta^{\frac{3n}{2}+1}}{D^{\frac{n+2}{2}}\cdot 2^{12 n}\cdot n^n}.
$$
\endproof

\section{Estimates about volumes of balls and spherical caps}

We prepare the proof of Theorem~\ref{Isodiametricstab} with a series of lemmas about balls.

\begin{lemma}
\label{ballvolumelow}
Let $\mathcal{M}^n$ be either $\R^n$, $H^n$ or $S^n$. Let $r$ be a positive number, for $\mathcal{M}^n=S^n$ we also assume $r\leq\frac{\pi}{2}$. For $0<s<r$ we can give the following lower bound for the volume of a ball of radius $r-s$:
$$
V_{\mathcal{M}^n}\left(B_{\mathcal{M}^n}\left(x_0,r-s\right)\right)\geq V_{\mathcal{M}^n}\left(B_{\mathcal{M}^n}\left(x_0,r\right)\right)-s\cdot\left(f_{\mathcal{M}^n}\left(r\right)\right)^{n-1}\cdot n\cdot\kappa_n
$$
where
$$
f_{\mathcal{M}^n}\left(t\right)=\left\{\begin{array}{lll}
    t & \rm{if} & \mathcal{M}^n=\R^n \\
    \sinh t & \rm{if} & \mathcal{M}^n=H^n\\
    \sin t & \rm{if} & \mathcal{M}^n=S^n
\end{array}\right..
$$
\end{lemma}
\proof
For $\varrho>0$ the Lebesgue measure of the ball of radius $\varrho$ in $\mathcal{M}^n$ is
\begin{equation}
\label{ballvolume}
V_{\mathcal{M}^n}\left(B_{\mathcal{M}^n}\left(x_0,\varrho\right)\right)=\int_0^{\varrho}{\left(f_{\mathcal{M}^n}\left(t\right)\right)^{n-1}\cdot n\cdot\kappa_n}\,\mathrm{d}t.
\end{equation}
We observe that $f_{\mathcal{M}^n}$ is monotonically increasing (for $\mathcal{M}^n=S^n$ we assume $\varrho\leq\frac{\pi}{2}$), so
\begin{gather*}
V_{\mathcal{M}^n}\left(B_{\mathcal{M}^n}\left(x_0,r\right)\right)=\int_0^{r-s}{\left(f_{\mathcal{M}^n}\left(t\right)\right)^{n-1}\cdot n\cdot\kappa_n}\,\mathrm{d}t+\int_{r-s}^r{\left(f_{\mathcal{M}^n}\left(t\right)\right)^{n-1}\cdot n\cdot\kappa_n}\,\mathrm{d}t\leq\\
\leq\int_0^{r-s}{\left(f_{\mathcal{M}^n}\left(t\right)\right)^{n-1}\cdot n\cdot\kappa_n}\,\mathrm{d}t+s\cdot\left(f_{\mathcal{M}^n}\left(r\right)\right)^{n-1}\cdot n\cdot\kappa_n=\\
=V_{\mathcal{M}^n}\left(B_{\mathcal{M}^n}\left(x_0,r-s\right)\right)+s\cdot\left(f_{\mathcal{M}^n}\left(r\right)\right)^{n-1}\cdot n\cdot\kappa_n.
\end{gather*}
\hfill\makebox{$\Box$}\\

\begin{lemma}
\label{ballvolumeup}
Let $\mathcal{M}^n$ be either $\R^n$, $H^n$ or $S^n$, and let $r>0$ where we also assume $r\leq\frac{\pi}{3}$ if $\mathcal{M}^n=S^n$. For $0<s<\frac{r}{2}$ we can give the following upper bounds for the volume of a ball of radii $r-s$ and $r+s$:
$$
V_{\mathcal{M}^n}\left(B_{\mathcal{M}^n}\left(x_0,r-s\right)\right)\leq V_{\mathcal{M}^n}\left(B_{\mathcal{M}^n}\left(x_0,r\right)\right)-s\cdot\left(f_{\mathcal{M}^n}\left(\frac{r}{2}\right)\right)^{n-1}\cdot n\cdot\kappa_n
$$
and
$$
V_{\mathcal{M}^n}\left(B_{\mathcal{M}^n}\left(x_0,r+s\right)\right)\leq V_{\mathcal{M}^n}\left(B_{\mathcal{M}^n}\left(x_0,r\right)\right)+s\cdot\left(f_{\mathcal{M}^n}\left(\frac{3r}{2}\right)\right)^{n-1}\cdot n\cdot\kappa_n
$$
where
$$
f_{\mathcal{M}^n}\left(t\right)=\left\{\begin{array}{lll}
    t & \rm{if} & \mathcal{M}^n=\R^n \\
    \sinh t & \rm{if} & \mathcal{M}^n=H^n\\
    \sin t & \rm{if} & \mathcal{M}^n=S^n
\end{array}\right..
$$
\end{lemma}
\proof
We use a similar argument as in Lemma~\ref{ballvolumelow}. Using \eqref{ballvolume} for $\varrho=r$ we have
\begin{gather*}
V_{\mathcal{M}^n}\left(B_{\mathcal{M}^n}\left(x_0,r\right)\right)=\int_0^{r-s}{\left(f_{\mathcal{M}^n}\left(t\right)\right)^{n-1}\cdot n\cdot\kappa_n}\,\mathrm{d}t+\int_{r-s}^r{\left(f_{\mathcal{M}^n}\left(t\right)\right)^{n-1}\cdot n\cdot\kappa_n}\,\mathrm{d}t\geq\\
\geq\int_0^{r-s}{\left(f_{\mathcal{M}^n}\left(t\right)\right)^{n-1}\cdot n\cdot\kappa_n}\,\mathrm{d}t+s\cdot\left(f_{\mathcal{M}^n}\left(\frac{r}{2}\right)\right)^{n-1}\cdot n\cdot\kappa_n=\\
=V_{\mathcal{M}^n}\left(B_{\mathcal{M}^n}\left(x_0,r-s\right)\right)+s\cdot\left(f_{\mathcal{M}^n}\left(\frac{r}{2}\right)\right)^{n-1}\cdot n\cdot\kappa_n
\end{gather*}
by the choice of $s$ and the monotonicity of $f_{\mathcal{M}^n}$.

By the choice $\varrho=r+s$ we can obtain the other inequality:
\begin{gather*}
V_{\mathcal{M}^n}\left(B_{\mathcal{M}^n}\left(x_0,r+s\right)\right)=\int_0^r{\left(f_{\mathcal{M}^n}\left(t\right)\right)^{n-1}\cdot n\cdot\kappa_n}\,\mathrm{d}t+\int_r^{r+s}{\left(f_{\mathcal{M}^n}\left(t\right)\right)^{n-1}\cdot n\cdot\kappa_n}\,\mathrm{d}t\leq\\
\leq\int_0^r{\left(f_{\mathcal{M}^n}\left(t\right)\right)^{n-1}\cdot n\cdot\kappa_n}\,\mathrm{d}t+s\cdot\left(f_{\mathcal{M}^n}\left(\frac{3r}{2}\right)\right)^{n-1}\cdot n\cdot\kappa_n=\\
=V_{\mathcal{M}^n}\left(B_{\mathcal{M}^n}\left(x_0,r\right)\right)+s\cdot\left(f_{\mathcal{M}^n}\left(\frac{3r}{2}\right)\right)^{n-1}\cdot n\cdot\kappa_n.
\end{gather*}
\hfill\makebox{$\Box$}\\
\noindent {\bf Remark} Note that we only used $r\leq\frac{\pi}{3}$ for the second inequality, the first upper estimate holds for $r\leq\frac{\pi}{2}$.

\begin{prop}
\label{error}
If $X\subset\mathcal{M}^n$ is compact with diameter at most $D>0$ (also $D<\frac{\pi}{2}$ if $\mathcal{M}^n=S^n$), $0<\varepsilon<V_{\mathcal{M}^n}\left(B_{\mathcal{M}^n}\left(x_0,\frac{D}{2}\right)\right)$ and $V_{\mathcal{M}^n}\left(X\right)\geq V_{\mathcal{M}^n}\left(B_{\mathcal{M}^n}\left(x_0,\frac{D}{2}\right)\right)-\varepsilon$, then for $0<\varrho\leq\frac{D}{2}$ satisfying also that $D+2\varrho<\frac{\pi}{2}$ if $\mathcal{M}^n=S^n$, we have
$$
V_{\mathcal{M}^n}\left(X^{\left(\varrho\right)}\right)\geq V_{\mathcal{M}^n}\left(B_{\mathcal{M}^n}\left(x_0,\frac{D}{2}+\varrho\right)\right)-E_{\mathcal{M}^n}\left(\varepsilon,D\right)
$$
where
$$
E_{\mathcal{M}^n}\left(\varepsilon,D\right)=\left\{\begin{array}{lll}
    \varepsilon\cdot 4^{n-1} & \rm{if} & \mathcal{M}^n=\R^n\mbox{ or }\mathcal{M}^n=S^n \\
    \varepsilon\cdot\left(4\cosh\frac{3D}{4}\right)^{n-1} & \rm{if} & \mathcal{M}^n=H^n
\end{array}\right..
$$
\end{prop}
\proof
We choose $\widetilde{E}_{\mathcal{M}^n}\left(\varepsilon,D\right)$ so that it satisfies
$$
V_{\mathcal{M}^n}\left(B_{\mathcal{M}^n}\left(x_0,\frac{D}{2}-\widetilde{E}_{\mathcal{M}^n}\left(\varepsilon,D\right)\right)\right)\leq V_{\mathcal{M}^n}\left(B_{\mathcal{M}^n}\left(x_0,\frac{D}{2}\right)\right)-\varepsilon.
$$
Let $R=\frac{D}{2}+\varrho$. Applying Lemma~\ref{ballvolumeup} we can set
$$
\widetilde{E}_{\mathcal{M}^n}\left(\varepsilon,D\right)=\frac{\varepsilon}{\left(f_{\mathcal{M}^n}\left(\frac{D}{4}\right)\right)^{n-1}\cdot n\cdot\kappa_n}
$$
where
$$
f_{\mathcal{M}^n}\left(t\right)=\left\{\begin{array}{lll}
    t & \rm{if} & \mathcal{M}^n=\R^n \\
    \sinh t & \rm{if} & \mathcal{M}^n=H^n\\
    \sin t & \rm{if} & \mathcal{M}^n=S^n
\end{array}\right..
$$
Hence the Isoperimetric Inequality Theorem~\ref{Isoperimetric} yields
$$
V_{\mathcal{M}^n}\left(X^{\left(\varrho\right)}\right)\geq V_{\mathcal{M}^n}\left(B_{\mathcal{M}^n}\left(x_0,R-\widetilde{E}_{\mathcal{M}^n}\left(\varepsilon,D\right)\right)\right).
$$
Now using Lemma~\ref{ballvolumelow} we have
\begin{gather*}
V_{\mathcal{M}^n}\left(B_{\mathcal{M}^n}\left(x_0,R-\widetilde{E}_{\mathcal{M}^n}\left(\varepsilon,D\right)\right)\right)\geq V_{\mathcal{M}^n}\left(B_{\mathcal{M}^n}\left(x_0,R\right)\right)-\widetilde{E}_{\mathcal{M}^n}\left(\varepsilon,D\right)\cdot f_{\mathcal{M}^n}\left(R\right)^{n-1}\cdot n\cdot\kappa_n=\\
=V_{\mathcal{M}^n}\left(B_{\mathcal{M}^n}\left(x_0,R\right)\right)-\varepsilon\cdot\left(\frac{f_{\mathcal{M}^n}\left(R\right)}{f_{\mathcal{M}^n}\left(\frac{R}{2}\right)}\right)^{n-1}.
\end{gather*}
For $\mathcal{M}^n=\R^n$ we are already done. It is trivial that
$$
\frac{f_{S^n}\left(R\right)}{f_{S^n}\left(\frac{D}{4}\right)}\leq\frac{f_{S^n}\left(R\right)}{f_{S^n}\left(\frac{R}{4}\right)}\leq 4.
$$
Finally, in the hyperbolic case
$$
\frac{f_{H^n}\left(R\right)}{f_{H^n}\left(\frac{D}{4}\right)}\leq\frac{f_{H^n}\left(R\right)}{f_{H^n}\left(\frac{R}{4}\right)}=4\cosh\left(\frac{R}{2}\right)\cosh\left(\frac{R}{4}\right)\leq4\cosh\left(\frac{3R}{4}\right),
$$
which finishes the proof.
\hfill$\Box$

\mbox{ }

If $\mathcal{M}^n$ is either $\R^n$, $H^n$ or $S^n$, then we need some lower bound on the volume of spherical caps.
If $H^+$ is a half-space in $\mathcal{M}^n$ such that 
$\partial_{\mathcal{M}^n}H^+\cap {\rm int}\,B_{\mathcal{M}^n}\left(x_0,t\right)\neq\emptyset$ for 
$x_0\in\mathcal{M}^n$ and $t>0$ where $t<\frac{\pi}2$ if $\mathcal{M}^n=S^n$, then
$C=H^+\cap B_{\mathcal{M}^n}\left(x_0,t\right)$ is called a spherical cap, and the depth $\delta$ of $C$ is the maximal distance of points of $C$ from $\partial_{\mathcal{M}^n}H^+$. We observe that if $\delta<t$, then 
the distance of $x_0$ from $H^+$ is $t-\delta$.

We need some estimate about the $(n-1)$-volume $\kappa_{n-1}$ of an $(n-1)$-dimensional Euclidean unit ball. It follows from $\frac{\kappa_{n-1}}{\kappa_n}>\sqrt{\frac{n}{2\pi}}$ and $\Gamma\left(x+1\right)<\left(\frac{x}{e}\right)^x\cdot\sqrt{2\pi\left(x+1\right)}$ that
\begin{equation}
\label{kappanest}   
\frac{\kappa_{n-1}}{n}>\frac{\kappa_n}{\sqrt{2n\pi}}=\frac{\pi^{\frac{n}{2}}}{\sqrt{2n\pi}\cdot\Gamma\left(\frac{n}{2}+1\right)}>\frac{\pi^{\frac{n}{2}}}{\sqrt{2n\pi}\cdot\left(\frac{n}{2e}\right)^{\frac{n}{2}}\cdot\sqrt{\pi\left(n+2\right)}}>\frac{1}{n^{\frac{n}{2}}}.
\end{equation}

In the Euclidean case, if $H\cap {\rm int}\,B_{\R^n}\left(x_0,t\right)\neq \emptyset$ for a hyperplane $H$ of $\R^n$
where the distance of $x_0$ from $H$ is $t-\delta$ for $\delta\in(0,t]$ (and hence the small cap cut off by $H$ is of depth $\delta$), then $H\cap B_{\R^n}\left(x_0,t\right)$ is an $(n-1)$-dimensional Euclidean ball of radius $a$ where
\begin{equation}
\label{Euclidballsectionest}
a=\sqrt{t^2-(t-\delta)^2}\geq \sqrt{t\delta}.
\end{equation}

\begin{lemma}
\label{Euclideancap}
For $x_0\in \R^n$, $t>0$ and $0<\delta\leq t$, if 
$H^+$ is a half-space in $\R^n$ such that 
$\partial_{\R^n}H^+\cap {\rm int}\,B_{\R^n}\left(x_0,t\right)\neq\emptyset$,
and $H^+\cap B_{\R^n}\left(x_0,t\right)$ is a spherical cap of depth at least $\delta$, then
$$
V_{\R^n}\left(H^+\cap B_{\R^n}\left(x_0,t\right)\right)\geq
\frac{2\kappa_{n-1}}{n+1}\cdot t^{\frac{n-1}{2}}\cdot\delta^{\frac{n+1}{2}}
\geq 
\frac{1}{n^{\frac{n}{2}}}\cdot t^{\frac{n-1}{2}}\cdot\delta^{\frac{n+1}{2}}.
$$
\end{lemma}
\proof It follows from applying first \eqref{Euclidballsectionest} and then \eqref{kappanest} that
$$
V_{\R^n}\left(H^+\cap B_{S^n}\left(x_0,t\right)\right)\geq
\int_0^\delta\kappa_{n-1}t^{\frac{n-1}{2}}\cdot s^{\frac{n-1}{2}}\,ds=
\frac{2\kappa_{n-1}}{n+1}\cdot t^{\frac{n-1}{2}}\cdot\delta^{\frac{n+1}{2}}
\geq 
\frac{1}{n^{\frac{n}{2}}}\cdot t^{\frac{n-1}{2}}\cdot\delta^{\frac{n+1}{2}}.
$$
\endproof

In Lemma~\ref{hyperboliccap}, we use two estimates for hyperbolic volume. Concerning balls, if $r>0$ and 
$x_0\in H^n$, then
\begin{equation}
\label{hypballvolest}
V_{H^n}\left(B_{H^n}\left(x_0,r\right)\right)=
\int_0^rn\kappa_n(\sinh s)^{n-1}\,ds\geq
\int_0^rn\kappa_n s^{n-1}\,ds=\kappa_nr^n.
\end{equation}
It is a reasonable lower  bound if $r\leq 2$.

Next let $\ell$ be any line in $H^n$, let $\delta>0$, and let $x_s\in \ell$ for $s\in [0,\delta]$ be
a parametrization of the segment $[x_0,x_\delta]_{H^n}\subset\ell$
where $d_{H^n}(x_0,x_s)=s$. In addition, let $H_s$ be the hyperplane in $H^n$ passing through $x_s$ and orthogonal to $\ell$. Now if $X$ is any compact set lying between $H_0$ and $H_\delta$, then the fact that the distance between any point of $H_s$ and any point of $H_t$ for $0\leq s<t\leq\delta$ is at least $t-s$ yields that
\begin{equation}
\label{hypsectionsvolest}
V_{H^n}\left(X\right)\geq
\int_0^\delta V_{H^{n-1}}\left(H_s\cap X\right) \,ds.
\end{equation}

\begin{lemma}
\label{hyperboliccap}
For $x_0\in H^n$, $t>0$ and $0<\delta\leq\min\left\{\frac{t}{2},1\right\}$, if $H^+$ is a half-space in $H^n$ 
such that $H^+\cap B_{H^n}(x_0,t)$ is a spherical cap of depth at least $\delta$, then
$$
V_{H^n}\left(H^+\cap B_{H^n}\left(x_0,t\right)\right)\geq n^{\frac{-(n-1)}{2}}\cdot \left(\tanh (t-\delta)\right)^{\frac{n-1}2}\cdot\delta^{\frac{n+1}{2}}.
$$
\end{lemma}
\proof We may assume that $H^+=H_\delta^+$ is a half-space in $H^n$ 
such that $H_\delta^+\cap B_{H^n}$ is a spherical cap of depth $\delta$.

First, let $w_0\in \partial_{H^n}B_{H^n}\left(x_0,t\right)$ 
be the point such that the segment $[x_0,w_0]_{H^n}$ intersects $\partial_{H^n}H_\delta^+$ in some
$w_\delta$.
For any $s\in[0,\delta]$, let $w_s$ be the point of $[x_0,w_0]_{H^n}$ with $d_{H^n}(x_0,w_s)=s$, and let  $H_s$ be the hyperplane of $H^n$ passing through $w_s$ and orthogonal to $[x_0,w_0]_{H^n}$. The Law of Cosines yields that $H_s$ intersects
$B_{H^n}\left(x_0,t\right)$ in an $(n-1)$-ball of radius $a_s$ where
\begin{eqnarray*}
\cosh a_s&=&\frac{\cosh t}{\cosh (t-s)}=\frac{\cosh (t-s)\cosh s+\sinh (t-s)\cdot\sinh s}{\cosh (t-s)}\\
&=&
\cosh s+\tanh (t-s)\cdot\sinh s\geq 1+s\cdot \tanh (t-\delta).
\end{eqnarray*}
Since $\cosh 1+\sinh 1=e\leq \cosh 2$ and $\cosh z\leq 1+z^2$ if $|z|\leq 2$, we deduce that
$$
a_s\geq \sqrt{s\cdot \tanh (t-\delta)}
$$
It follows from \eqref{hypballvolest} that if
$s\in[0,\delta]$, then
$$
V_{H^{n-1}}\left(H_s\cap B_{H^n}\left(x_0,t\right)\right)\geq \kappa_{n-1}a_s^{n-1}\geq \kappa_{n-1}\cdot
s^{\frac{n-1}2}\cdot \left(\tanh (t-\delta)\right)^{\frac{n-1}2};
$$
therefore, \eqref{hypsectionsvolest}
yields that
$$
V_{H^n}\left(H_\delta^+\cap B_{H^n}\left(x_0,t\right)\right)\geq
\int_0^\delta \kappa_{n-1}\cdot
s^{\frac{n-1}2}\cdot \left(\tanh (t-\delta)\right)^{\frac{n-1}2}\,ds=
\kappa_{n-1}\cdot \left(\tanh (t-\delta)\right)^{\frac{n-1}2}
\cdot \delta^{\frac{n+1}2}.
$$
Since \eqref{kappanest} yields 
$\kappa_{n-1}\geq n^{\frac{-(n-1)}{2}}$,
we conclude Lemma~\ref{hyperboliccap}.
\endproof

In the case of spherical caps $C$, the analogue of \eqref{hypsectionsvolest} does not hold in the spherical space; therefore, we estimate the volume of $C$ by projecting $C$ into a Euclidean space of the same dimension.

\begin{lemma}
\label{sphericalcap}
For $x_0\in S^n$, $0<t<\frac{\pi}{2}$ and $0<\delta\leq\frac{t}{2}$, if $H^+$ is a hemisphere in $S^n$ whose distance from $x_0$ is $t-\delta$, then
$$
V_{S^n}\left(H^+\cap B_{S^n}\left(x_0,t\right)\right)\geq\frac{1}{2^n\cdot n^{\frac{n}{2}}}\cdot t^{\frac{n-1}{2}}\cdot\delta^{\frac{n+1}{2}}.
$$
\end{lemma}
\proof
The intersection $\partial_{S^n}H^+\cap B_{S^n}\left(x_0,t\right)$ is an $\left(n-1\right)$-ball in $S^n$ of some radius $a$. Let $s=t-\delta$, then
$$
\cos\left(s+\delta\right)=\cos\left(t\right)=\cos\left(s\right)\cdot\cos\left(a\right).
$$
We claim that $a\geq\sqrt{s\delta}$, which is equivalent to
\begin{equation}
\label{cosineq}
\cos\left(s\right)\cdot\cos\left(\sqrt{s\delta}\right)\geq\cos\left(s+\delta\right).
\end{equation}
This holds with equality if $\delta=0$, so it is sufficient to see
$$
\cos\left(s\right)\cdot\frac{\sin\left(\sqrt{s\delta}\right)}{\sqrt{s\delta}}\cdot\frac{s}{2}\leq\sin\left(s+\delta\right).
$$
after differentiation. Since $0<s<\frac{\pi}{2}$, we have $s<\tan\left(s\right)$, so
$$
\cos\left(s\right)\cdot\frac{\sin\left(\sqrt{s\delta}\right)}{\sqrt{s\delta}}\cdot\frac{s}{2}<\cos\left(s\right)\cdot\frac{s}{2}<\frac{\sin\left(s\right)}{2}.
$$
Here $\delta<\frac{\pi}{3}$ so $\frac{1}{2}<\cos\left(\delta\right)$, therefore
$$
\cos\left(s\right)\cdot\frac{\sin\left(\sqrt{s\delta}\right)}{\sqrt{s\delta}}\cdot\frac{s}{2}<\sin\left(s\right)\cdot\cos\left(\delta\right)<\sin\left(s+\delta\right),
$$
proving \eqref{cosineq}.

Let $w$ be the closest point of $H^+$ to $x_0$, and let $\Pi\colon S^n\to w^{\perp}$ be the orthogonal projection. Then $\Pi\left(\partial_{S^n}H^+\cap B_{S^n}\left(x_0,t\right)\right)$ is a Euclidean $\left(n-1\right)$-ball of radius $\sin\left(a\right)$. By the choice of $H^+$ and from $d_{S^n}\left(x_0,w\right)=t-\delta$, there is a point $\widetilde{w}\in H^+\cap\partial_{S^n}B_{S^n}\left(x_0,t\right)$ such that $d_{S^n}\left(w,\widetilde{w}\right)=\delta$ and $w\in\left[\widetilde{w},x_0\right]_{S^n}$. Hence, there is a cone $C\subset\Pi\left(H^+\cap B_{S^n}\left(x_0,t\right)\right)$ with base $\Pi\left(\partial_{S^n}H^+\cap B_{S^n}\left(x_0,t\right)\right)$ and height $\sin\left(\delta\right)$. Thus,
$$
V_{S^n}\left(H^+\cap B_{S^n}\left(x_0,t\right)\right)\geq V_{w^{\perp}}\left(\Pi\left(H^+\cap B_{S^n}\left(x_0,t\right)\right)\right)\geq V_{w^{\perp}}\left(C\right)=\frac{\kappa_{n-1}}{n}\cdot\left(\sin\left(a\right)\right)^{n-1}\cdot\sin\left(\delta\right).
$$
From \eqref{cosineq} we can imply
$$
\sin\left(a\right)\geq\sin\left(\sqrt{\frac{t\delta}{2}}\right),
$$
so from $\sqrt{\frac{t\cdot\delta}{2}}<\frac{\pi}{4}$ and $\delta<\frac{\pi}{3}$ we have
$$
V_{S^n}\left(H^+\cap B_{S^n}\left(x_0,t\right)\right)\geq\frac{\kappa_{n-1}}{n}\cdot\left(\frac{\sqrt{t\cdot\delta}}{2}\right)^{n-1}\cdot\frac{\delta}{2}=\frac{1}{2^n}\cdot t^{\frac{n-1}{2}}\cdot\delta^{\frac{n+1}{2}}.
$$
Finally, \eqref{kappanest} completes the proof.
\endproof

\section{The proof of Theorem~\ref{Isodiametricstab}}

We need some upper bound on the volume of balls in
$\mathcal{M}^n$ where $\mathcal{M}^n$ is either $\R^n$, $S^n$ or $H^n$.

\begin{lemma}
\label{kappanballestup}
If  $\mathcal{M}^n$ is either $\R^n$, $S^n$ or $H^n$, and $r>0$ where $r<\frac{\pi}2$ if $\mathcal{M}^n=S^n$, then
$$
V_{\mathcal{M}^n}\left(B_{\mathcal{M}^n}\left(z_0,r\right)\right)\leq
\left\{
\begin{array}{ll}
  \frac{2^{3n}}{n^{\frac{n+1}{2}}}\cdot r^n&\mbox{ \ if either $\mathcal{M}^n=\R^n$, or $\mathcal{M}^n=S^n$, or
  $\mathcal{M}^n=H^n$ and $r\leq 1$}\\[2ex]
 \frac{2^{\frac{3n+2}2}e^{(n-1)r}}{n^{\frac{n+1}{2}}} &\mbox{ \ if $\mathcal{M}^n=H^n$ and $r>0$.}
\end{array} \right. 
$$
\end{lemma}
\noindent{\bf Remark } The second bound in the hyperbolic case is worse than the first bound if $r\leq 1$.
\proof We set $r=D/2$. In the Euclidean case, it follows from  
$\Gamma\left(x+1\right)>\left(\frac{x}{e}\right)^x\cdot\sqrt{2\pi x}$ for $x\geq 1$ 
and from  $n\geq 2$ and $\frac1{\sqrt{\pi}}\cdot \left(\frac{e\pi}{2^4}\right)^{\frac{n}2}\leq \frac12$  that 
\begin{equation}
\label{kappanestupEuc}   
\kappa_n=\frac{\pi^{\frac{n}{2}}}{\Gamma\left(\frac{n}{2}+1\right)}<
\frac{\pi^{\frac{n}{2}}}{\left(\frac{n}{2e}\right)^{\frac{n}{2}}\cdot\sqrt{\pi n}} =\frac1{\sqrt{\pi}}\cdot \left(\frac{e\pi}{2^4}\right)^{\frac{n}2}\cdot \frac{2^{\frac{5n}2}}{n^{\frac{n+1}{2}}}<
\frac12 \cdot \frac{2^{\frac{5n}2}}{n^{\frac{n+1}{2}}}.
\end{equation}
In the spherical case, if $r\in(0,\frac{\pi}2)$, then
\begin{equation}
\label{kappanestupSn}
V_{S^n}\left(B_{S^n}\left(z_0,r\right)\right)
=n\kappa_n\int_0^r(\sin s)^{n-1}\,ds\leq
n\kappa_n\int_0^rs^{n-1}\,ds=r^n\kappa_n
 <\frac{2^{\frac{5n}2}}{n^{\frac{n+1}{2}}}\cdot r^n.
\end{equation}
Finally, in the hyperbolic 
case, if $r\in(0,1]$, then $\sinh s<\sqrt{2}\,s$ for  $s\in(0,1]$ yields
\begin{equation}
\label{kappanestupHn}
V_{H^n}\left(B_{H^n}\left(z_0,r\right)\right)
=n\kappa_n\int_0^r(\sinh s)^{n-1}\,ds\leq
n\kappa_n\int_0^r\sqrt{2}^{n-1}s^{n-1}\,ds
 <\frac{2^{3n}}{n^{\frac{n+1}{2}}}\cdot r^n.
\end{equation}
Therefore,  \eqref{kappanestupEuc}, \eqref{kappanestupSn} and \eqref{kappanestupHn} yield 
Lemma~\ref{kappanballestup} if either $\mathcal{M}^n=\R^n$, or $\mathcal{M}^n=S^n$, or
  $\mathcal{M}^n=H^n$ and $r\leq 1$.

Finally, we consider the case  $\mathcal{M}^n=H^n$ and $r\geq 1$. As a first step, we observe that
$\sinh s\leq e^s/2$ for $s\geq 0$; therefore,
we deduce using \eqref{kappanestupEuc} and $\frac{n}{n-1}\leq 2$ that
\begin{eqnarray*}
V_{H^n}\left(B_{H^n}\left(z_0,r\right)\right)
&=&n\kappa_n\int_0^r(\sinh s)^{n-1}\,ds\leq
n\kappa_n\int_0^r\frac{e^{(n-1)s}}{2^{n-1}}\,ds\\
&\leq&
n\kappa_n\cdot \frac{e^{(n-1)r}}{(n-1)2^{n-1}}\leq
\frac{2\cdot 2^{\frac{3n}2}e^{(n-1)r}}{n^{\frac{n+1}{2}}}
\leq
\frac{2^{\frac{3n+2}2}e^{(n-1)r}}{n^{\frac{n+1}{2}}},
\end{eqnarray*}
proving Lemma~\ref{kappanballestup} also if $\mathcal{M}^n=H^n$ and $r\geq 1$.
\endproof

Now we are ready to prove Theorem~\ref{Isodiametricstab}, which we are restating including the exact values of
$\varepsilon_{\mathcal{M}^n}\left(D\right)$.

\begin{theorem}
\label{Isodiametricstab0} 
 For $n\geq 2$ if $\mathcal{M}^n$ is either $\R^n$, $S^n$ or $H^n$, $D>0$ (where $D<\frac{\pi}{2}$ if $\mathcal{M}^n=S^n$) and $X\subset \mathcal{M}^n$ is measurable with ${\rm diam}X\leq D$ and
$$
V_{\mathcal{M}^n}\left(X\right)\geq (1-\varepsilon)
V_{\mathcal{M}^n}\left(B_{\mathcal{M}^n}\left(z_0,\frac{D}{2}\right)\right),
$$
for $\varepsilon\in\left[0,\varepsilon_{\mathcal{M}^n}\left(D\right)\right)$, then there exists a $c\in \mathcal{M}^n$ such that
$$
B\left(c,\mbox{$\frac{D}2$}-\gamma_{\mathcal{M}^n}\left(D\right)\cdot \varepsilon^{\frac2{3n+2}}\right)\subset {\rm conv}_{\mathcal{M}^n}X\subset B\left(c,\mbox{$\frac{D}2$}+\gamma_{\mathcal{M}^n}\left(D\right)\cdot\varepsilon^{\frac2{3n+2}}\right)
$$
where
$$
\gamma_{\mathcal{M}^n}\left(D\right)=
\left\{
\begin{array}{ll}
  e^{21}n\cdot D&\mbox{ if $\mathcal{M}^n=H^n$ and $D\leq 2$, or $\mathcal{M}^n=\R^n$, or $\mathcal{M}^n=S^n$;}\\[1ex]
  n\cdot e^{7D+8}
 &\mbox{ if $\mathcal{M}^n=H^n$ and $D\geq 1$,}
\end{array} \right. 
$$
$$
\varepsilon_{\mathcal{M}^n}\left(D\right)=
\left\{\begin{array}{lll}
e^{-28n}n^{-\frac{n}{2}}
 & \mbox{ if } & \mathcal{M}^n=\R^n,\mbox{ or }\mathcal{M}^n=S^n
 \mbox{ and } D\leq \frac{\pi}6,\\
 &&\mbox{ or }\mathcal{M}^n=H^n \mbox{ and } D\leq 2;\\
e^{-30n}n^{-\frac{n}{2}}
\left(\frac{\pi}2-D\right)^{3n+2}
 & \mbox{ if } & \mathcal{M}^n=S^n
 \mbox{ and } \frac{\pi}6\leq D< \frac{\pi}2;\\
 e^{-18D}n^{-\frac{n}{2}}
 & \mbox{ if } & \mathcal{M}^n=H^n
\mbox{ and } D\geq 2
\end{array}\right.
$$
In addition, 
$V_{\mathcal{M}^n}\left(\left({\rm conv}_{\mathcal{M}^n}X\right)\backslash X\right)\leq \varepsilon\cdot 
V_{\mathcal{M}^n}\left(B_{\mathcal{M}^n}\left(z_0,\frac{D}{2}\right)\right)$.
\end{theorem}

\proof
Let $\mathcal{M}^n$ be either $\R^n$,  $H^n$ or $S^n$, 
let $D>0$ where $D<\frac{\pi}2$  if $\mathcal{M}^n=S^n$, and let $X\subset \mathcal{M}^n$ be measurable satisfying ${\rm diam}\, X\leq D$ 
and $V(X)\geq (1-\varepsilon)V(B(z_0,\frac{D}2))$. 
We set
$$
\widetilde{\varepsilon}=\varepsilon\cdot V\left(B\left(z_0,\frac{D}{2}\right)\right)
\mbox{ \ and \ }
\widetilde{\varepsilon}_{\mathcal{M}^n}\left(D\right)
=\varepsilon_{\mathcal{M}^n}\left(D\right)\cdot V\left(B\left(z_0,\frac{D}{2}\right)\right),
$$
and hence 
$V(X)\geq V(B(z_0,\frac{D}2))-\widetilde{\varepsilon}$
and
$\widetilde{\varepsilon}<  \widetilde{\varepsilon}_{\mathcal{M}^n}\left(D\right)$.

We consider a $\widetilde{X}\subset \mathcal{M}^n$
that has maximal volume under the conditions $X\subset \widetilde{X}$  and
${\rm diam}\, \widetilde{X}\leq D$. In particular, $\widetilde{X}$ is 
a convex body of constant width $D$ according to Corollary~\ref{maximal-convex}.
We fix some $x_1,x_2\in \widetilde{X}$ satisfying that $d_{\mathcal{M}^n}(x_1,x_2)=D$,
let $y_0$
be the midpoint of $[x_1,x_2]_{\mathcal{M}^n}$, and let 
$$
\varrho=\left\{
\begin{array}{ll}
\frac{D}2&\mbox{ \ if $\mathcal{M}^n=\R^n$ or $\mathcal{M}^n=H^n$,}\\[1ex]
\min\left\{\frac{D}{2},\frac{\pi}{8}-\frac{D}{4}\right\}&\mbox{ \ if $\mathcal{M}^n=S^n$ and $D<\frac{\pi}{2}$}.
\end{array} \right.
$$
be the $\varrho$ in Proposition~\ref{Improve} that satisfies $D+2\varrho<\frac{\pi}2$ if $\mathcal{M}^n=S^n$.

We consider the parallel domain $\widetilde{X}^{(\varrho)}$, which is a convex body of constant width $D+2\varrho$ by
Lemma~\ref{completeparallel}. In addition, Proposition~\ref{error} yields that
\begin{equation}
\label{tildeXrhovol}
V_{\mathcal{M}^n}\left(\widetilde{X}^{\left(\varrho\right)}\right)\geq V_{\mathcal{M}^n}\left(B_{\mathcal{M}^n}\left(x_0,\frac{D}{2}+\varrho\right)\right)-
\widetilde{E}_{\mathcal{M}^n}\left(\widetilde{\varepsilon},D\right)
\end{equation}
where
$$
\widetilde{E}_{\mathcal{M}^n}\left(\widetilde{\varepsilon},D\right)= \left\{\begin{array}{lll}
    \widetilde{\varepsilon}\cdot 4^n & \mbox{ if } & \mathcal{M}^n=\R^n\mbox{ or }\mathcal{M}^n=S^n \\
    \widetilde{\varepsilon}\cdot 4^n\left(\cosh D\right)^{n} & \mbox{ if } & \mathcal{M}^n=H^n
\end{array}\right..
$$
As we verify it at the end of the proof of Theorem~\ref{Isodiametricstab0}, we have
\begin{equation}
\label{Egammaeta}
\widetilde{E}_{\mathcal{M}^n}\left(\widetilde{\varepsilon}_{\mathcal{M}^n}\left(D\right),D\right)< \widetilde{\gamma}_1\eta_0^{\frac{3n+2}2}
\end{equation}
where the constants $\widetilde{\gamma}_1$ and $\eta_0$
come from Proposition~\ref{Improve}; namely,
$$
\widetilde{\gamma}_1=
\left\{
\begin{array}{ll}
\frac{1}{2^{12n}n^n}\cdot
\frac{1}{D^{\frac{n+2}2}}&\mbox{ \ if $\mathcal{M}^n=\R^n$ or $\mathcal{M}^n=S^n$ and $D<\frac{\pi}2$;}\\[1ex]
\frac{1}{2^{8n}n^n}\cdot
\frac{1}{(\sinh{10D})^{\frac{n+2}2}}&\mbox{ \ if $\mathcal{M}^n=H^n$;}
\end{array} \right.
$$

We recall that two-point symmetrization preserves Lebesgue measure and does not increase diameter, and hence
$$
V_{\mathcal{M}^n}\left(\widetilde{X}^{\left(\varrho\right)}\right)=
V_{\mathcal{M}^n}\left(\tau_H\widetilde{X}^{\left(\varrho\right)}\right) \leq
V_{\mathcal{M}^n}\left(\mathrm{conv}_{\mathcal{M}^n}\tau_H\widetilde{X}^{\left(\varrho\right)}\right)\leq V_{\mathcal{M}^n}\left(B_{\mathcal{M}^n}\left(z_0,\frac{D}{2}+\varrho\right)\right)
$$ 
by the Isodiametric Inequality~\ref{Isodiametric}. 
Let
$$
\frac{D}2+\eta=\min\left\{\frac{D}2+\eta_0,\max_{z\in \widetilde{X}}d_{\mathcal{M}^n}(z,y_0)\right\},
$$
and hence Proposition~\ref{Improve} and \eqref{tildeXrhovol} yield the
existence of a hyperplane $H$ such that
\begin{eqnarray*}
\widetilde{\gamma}_1 \eta^{\frac{3n+2}2}&\leq&
V_{\mathcal{M}^n}\left(\mathrm{conv}_{\mathcal{M}^n}\tau_H\widetilde{X}^{\left(\varrho\right)}\right)-V_{\mathcal{M}^n}\left(\tau_H\widetilde{X}^{\left(\varrho\right)}\right)=V_{\mathcal{M}^n}\left(\mathrm{conv}_{\mathcal{M}^n}\tau_H\widetilde{X}^{\left(\varrho\right)}\right)-V_{\mathcal{M}^n}\left(\widetilde{X}^{\left(\varrho\right)}\right)\\
&&\leq V_{\mathcal{M}^n}\left(B_{\mathcal{M}^n}\left(z_0,\frac{D}{2}+\varrho\right)\right)-V_{\mathcal{M}^n}\left(\widetilde{X}^{\left(\varrho\right)}\right)\leq \widetilde{E}_{\mathcal{M}^n}\left(\widetilde{\varepsilon},D\right).
\end{eqnarray*}
It follows, using the condition \eqref{Egammaeta} and
$\sinh 10D\cdot(\cosh  D)^2\leq \sinh 11D\cdot\cosh  D\leq \sinh 12D$, that
$$
\frac{\eta}{\widetilde{\varepsilon}^{\frac2{3n+2}}}\leq 
\left\{
\begin{array}{ll}
\left(2^{12n}n^n4^n\cdot
D^{\frac{n+2}2}\right)^{\frac2{3n+2}}\leq 2^{10}nD^{\frac{n+2}{3n+2}}&\mbox{ \ if $\mathcal{M}^n=\R^n$ or $\mathcal{M}^n=S^n$;}\\[1ex]
\left(2^{8n}n^n
(\sinh{10D})^{\frac{n+2}2}\cdot 4^n(\cosh D)^n\right)^{\frac2{3n+2}}\leq
2^7n(\sinh{12 D})^{\frac{n+2}{3n+2}}&\mbox{ \ if $\mathcal{M}^n=H^n$.}
\end{array} \right. .
$$
We observe that $\frac{n+2}{3n+2}\leq\frac12$, and  $\sinh t\leq e^{22}t$ for $t\in[0,24]$; thus 
if $D\leq 2$, then
$$
(\sinh{12 D})^{\frac{n+2}{3n+2}}\leq 
\left\{
\begin{array}{ll}
 \left(e^{22} \cdot 12 \cdot D\right)^{\frac{n+2}{3n+2}}\leq
e^{13}\cdot D^{\frac{n+2}{3n+2}} &\mbox{ if $D\leq 2$;}\\[2ex]
e^{12D\cdot \frac{n+2}{3n+2}}\leq e^{6D}  &\mbox{ if  $D\geq 1$.}
\end{array} \right.
$$
Next, Lemma~\ref{kappanballestup} yields 
\begin{eqnarray}
\label{epstildeeps}
\frac{\widetilde{\varepsilon}^{\frac2{3n+2}}}{\varepsilon^{\frac2{3n+2}}}&=&V\left(B\left(z_0,\frac{D}2\right)\right)^{\frac2{3n+2}}\\
\nonumber
&\leq &
\left\{
\begin{array}{ll}
  \left(\frac{2^{3n}}{n^{\frac{n+1}{2}}}\cdot \frac{D^n}{2^n}\right)^{\frac2{3n+2}}
  \leq \frac{e\cdot D^{\frac{2n}{3n+2}}}{n^{\frac13}}&\mbox{ if $\mathcal{M}^n=\R^n$, or $\mathcal{M}^n=S^n$,}\\
  &\mbox{ or $\mathcal{M}^n=H^n$ and $D\leq 2$;}\\[2ex]
 \left(\frac{2^{\frac{3n+2}2}e^{\frac{(n-1)D}2}}{n^{\frac{n+1}{2}}}\right)^{\frac2{3n+2}} 
 \leq \frac{e^{\frac{D}3+1}}{n^{\frac13}}
 &\mbox{ if $\mathcal{M}^n=H^n$ and $D\geq 1$.}
\end{array} \right. 
\end{eqnarray}
Since $2^{7}<e^6$ and $2^{10}<e^7$, we deduce 
from the estimates above that
$$
R\left(\widetilde{X}\right)\leq \max_{z\in \widetilde{X}}d_{\mathcal{M}^n}(z,y_0)\leq
\frac{D}2+\widetilde{\gamma}_2\varepsilon^{\frac2{3n+2}}
$$
where 
$$
\widetilde{\gamma}_2= 
\left\{
\begin{array}{ll}
  e^{20}n\cdot D&\mbox{ if $\mathcal{M}^n=H^n$ and $D\leq 2$, or $\mathcal{M}^n=\R^n$, or $\mathcal{M}^n=S^n$;}\\[1ex]
  e^{7D+7}n
 &\mbox{ if $\mathcal{M}^n=H^n$ and $D\geq 1$.}
\end{array} \right. 
$$
We deduce from  $\varepsilon<\varepsilon_{\mathcal{M}^n}(D)$ 
(see the argument at the end of the proof of Theorem~\ref{Isodiametricstab0}) that
\begin{equation}
\label{tildegamma2cond}
\widetilde{\gamma}_2\varepsilon^{\frac2{3n+2}}\leq\frac{D}{8}.
\end{equation}
In turn, Lemma~\ref{r_plus_R} yields that
$$
r\left(\widetilde{X}\right)\geq 
\frac{D}2-\widetilde{\gamma}_2\varepsilon^{\frac2{3n+2}}
\geq \frac{3D}{8}.
$$
In particular, writing $c$ to denote the circumcenter of $\widetilde{X}$, we have
\begin{equation}
\label{CompletionAnnulus}    
B_{\mathcal{M}^n}\left(c,\frac{D}2-\widetilde{\gamma}_2\varepsilon^{\frac2{3n+2}}\right)
\subset \widetilde{X}\subset
B_{\mathcal{M}^n}\left(c,\frac{D}2+\widetilde{\gamma}_2\varepsilon^{\frac2{3n+2}}\right).
\end{equation}

Next, we write $\overline{X}={\rm conv}\,X$, and hence
\begin{equation}
\label{epsdiffchain}
V_{\mathcal{M}^n}\left(B_{\mathcal{M}^n}\left(c,\frac{D}2\right)\right)-\widetilde{\varepsilon}\leq 
V_{\mathcal{M}^n}\left(X\right)\leq 
V_{\mathcal{M}^n}\left(\overline{X}\right)\leq
V_{\mathcal{M}^n}\left(\widetilde{X}\right)\leq
V_{\mathcal{M}^n}\left(B_{\mathcal{M}^n}\left(c,\frac{D}2\right)\right).
\end{equation}
For any $x\in\partial_{\mathcal{M}^n}\overline{X}$,
writing $H^+_x$ to denote the closed ``supporting'' half-space
of $\mathcal{M}^n$ such that $x\in H^+_x$ and
$\overline{X}\cap {\rm int}\,H^+_x=\emptyset$, 
we deduce that
\begin{equation}
\label{capvolepsilon}
V_{\mathcal{M}^n}\left(B_{\mathcal{M}^n}\left(c,\frac{D}2-\widetilde{\gamma}_2\varepsilon^{\frac2{3n+2}}\right)\cap H^+_x\right)\leq \widetilde{\varepsilon}.
\end{equation}
We write $\delta_x$ to denote the depth of
$B_{\mathcal{M}^n}\left(c,\frac{D}2-\widetilde{\gamma}_2\widetilde{\varepsilon}^{\frac2{3n+2}}\right)\cap H^+_x$ for
$x\in\partial_{\mathcal{M}^n}\overline{X}$, and hence
combining \eqref{capvolepsilon} with Lemmas~\ref{Euclideancap}, \ref{hyperboliccap}
and \ref{sphericalcap} yield that
$$
\widetilde{\varepsilon}\geq 
\left\{
\begin{array}{ll}
 \frac{\left(\frac{D}4\right)^{\frac{n-1}{2}}}{n^{\frac{n}{2}}}\left(\min\left\{\delta_x,\frac{D}8\right\}\right)^{\frac{n+1}{2}}
 \geq \frac{D^{\frac{n-1}{2}}}{2^{n}\cdot n^{\frac{n}{2}}}\left(\min\left\{\delta_x,\frac{D}8\right\}\right)^{\frac{n+1}{2}}&\mbox{ \ if $\mathcal{M}^n=\R^n$,}\\[2ex]
 \frac{\left(\frac{D}4\right)^{\frac{n-1}{2}}}{2^n\cdot n^{\frac{n}{2}}}\left(\min\left\{\delta_x,\frac{D}8\right\}\right)^{\frac{n+1}{2}}
 \geq \frac{D^{\frac{n-1}{2}}}{2^{2n}\cdot n^{\frac{n}{2}}}\left(\min\left\{\delta_x,\frac{D}8\right\}\right)^{\frac{n+1}{2}}&\mbox{ \ if $\mathcal{M}^n=S^n$}\\[2ex]
 \frac{\left(\tanh\frac{D}4\right)^{\frac{n-1}{2}}}{n^{\frac{n}{2}}}\left(\min\left\{\delta_x,1,\frac{D}8\right\}\right)^{\frac{n+1}{2}}&\mbox{ \ if $\mathcal{M}^n=H^n$.}
\end{array} \right. 
$$
Since $\widetilde{\varepsilon}\leq \widetilde{\varepsilon}_{\mathcal{M}}(D)$, we deduce that
if $x\in\partial_{\mathcal{M}^n}\overline{X}$, then
$\delta_x<\frac{D}8$ if $\mathcal{M}^n=\R^n$ or $\mathcal{M}^n=S^n$, and 
$\delta_x<\min\left\{1,\frac{D}8\right\}$ if $\mathcal{M}^n=H^n$. Therefore,
if $x\in\partial_{\mathcal{M}^n}\overline{X}$, then
$$
\delta_x\leq \widetilde{\gamma}_3\widetilde{\varepsilon}^{\frac{2}{n+1}}
$$
where
$$
\widetilde{\gamma}_3=
\left\{
\begin{array}{ll}
  \frac{4n}{D^{\frac{n-1}{n+1}}}&\mbox{ \ if $\mathcal{M}^n=\R^n$,}\\[2ex]
 \frac{8n}{D^{\frac{n-1}{n+1}}}&\mbox{ \ if $\mathcal{M}^n=S^n$}\\[2ex]
 \frac{n}{\left(\tanh\frac{D}4\right)^{\frac{n-1}{n+1}}}&\mbox{ \ if $\mathcal{M}^n=H^n$.}
\end{array} \right. 
$$
It follows from Lemma~\ref{kappanballestup} that
$$
V\left(B\left(z_0,\frac{D}2\right)\right)^{\frac{2}{n+1}}\leq
\left\{
\begin{array}{ll}
  \frac{2^4}{n}\cdot D^{\frac{2n}{n+1}}&\mbox{ \ if either $\mathcal{M}^n=\R^n$, or $\mathcal{M}^n=S^n$, or
  $\mathcal{M}^n=H^n$ and $D\leq 2$}\\[2ex]
 \frac{e^{\frac{2D}3+2}}{n} &\mbox{ \ if $\mathcal{M}^n=H^n$ and $r>0$.}
\end{array} \right. 
$$

Now $\widetilde{\varepsilon}=\varepsilon\cdot V\left(B\left(z_0,\frac{D}{2}\right)\right)$, 
$\left(\tanh\frac{D}4\right)^{-\frac{n-1}{n+1}}< e$ for $D\geq 1$
and $\tanh t\geq \frac{t}2$ for $t\in[0,1]$ imply
$$
\delta_x\leq \widetilde{\gamma}_4\varepsilon^{\frac{2}{n+1}}
$$
for any $x\in\partial_{\mathcal{M}^n}\overline{X}$ where
$$
\widetilde{\gamma}_4=
\left\{
\begin{array}{ll}
  2^7 D&\mbox{ \ if either $\mathcal{M}^n=\R^n$, or $\mathcal{M}^n=S^n$, or
  $\mathcal{M}^n=H^n$ and $D\leq 2$}\\[2ex]
 e^{D+3}&\mbox{ \ if $\mathcal{M}^n=H^n$ and $D\geq 1$.}
\end{array} \right. 
$$
In turn, we conclude from \eqref{CompletionAnnulus} that   
$$
B_{\mathcal{M}^n}\left(c,\frac{D}2-e\cdot\widetilde{\gamma}_2\varepsilon^{\frac2{3n+2}}\right)
\subset \overline{X}\subset
B_{\mathcal{M}^n}\left(c,\frac{D}2+\widetilde{\gamma}_2\varepsilon^{\frac2{3n+2}}\right).
$$

Finally, \eqref{epsdiffchain} yields that
$V_{\mathcal{M}^n}\left(\overline{X}\backslash X\right)\leq \varepsilon
V_{\mathcal{M}^n}\left(B_{\mathcal{M}^n}\left(z_0,\frac{D}{2}\right)\right)$.

We still need to verify that our choice of 
$\varepsilon_{\mathcal{M}^n}\left(D\right)$ works; in other words,
it satisfies \eqref{Egammaeta} and \eqref{tildegamma2cond}.
We prove \eqref{Egammaeta} on a case-by-case basis
using \eqref{epstildeeps} and the values of
$\widetilde{\gamma}_1$ and $\eta_0$
in Proposition~\ref{Improve}.
It follows from $\widetilde{\varepsilon}=
\varepsilon\cdot V\left(B\left(z_0,\frac{D}{2}\right)\right)$,
$\cosh D<4$ if $0<D\leq 2$, $\cosh D<e^D$ if $D\geq 2$
and Lemma~\ref{kappanballestup} that
\begin{eqnarray*}
\widetilde{E}_{\mathcal{M}^n}\left(\widetilde{\varepsilon},D\right)^{\frac2{3n+2}}&=& \left\{\begin{array}{lll}
\varepsilon^{\frac2{3n+2}}\cdot V\left(B\left(z_0,\frac{D}{2}\right)\right)^{\frac2{3n+2}}\cdot 4^{\frac{2n}{3n+2}} & \mbox{ if } & \mathcal{M}^n=\R^n\mbox{ or }\mathcal{M}^n=S^n \\
\varepsilon^{\frac2{3n+2}}\cdot V\left(B\left(z_0,\frac{D}{2}\right)\right)^{\frac2{3n+2}}\cdot 4^{\frac{2n}{3n+2}}
\left(\cosh D\right)^{\frac{2n}{3n+2}} & \mbox{ if } & \mathcal{M}^n=H^n
\end{array}\right.\\
&\leq&
\left\{\begin{array}{lll}
\varepsilon^{\frac2{3n+2}}\cdot D^{\frac{2n}{3n+2}}\cdot
\frac{2^4}{n^{\frac{n+1}{3n+2}}} & \mbox{ if } & \mathcal{M}^n=\R^n\mbox{ or }\mathcal{M}^n=S^n
\mbox{ or }\mathcal{M}^n=H^n \mbox{ and } D\leq 2\\
\varepsilon^{\frac2{3n+2}}\cdot  e^D\cdot \frac{2^{\frac73}}{n^{\frac{n+1}{3n+2}}}
 & \mbox{ if } & \mathcal{M}^n=H^n
\mbox{ and } D\geq 2
\end{array}\right.
\end{eqnarray*}
We deduce from $2^4<e^3$, $2^{\frac73}<e^2$, $2^{12n\cdot \frac{2}{3n+2}}<2^8<e^{6}$, 
$2^{8n \cdot \frac{2}{3n+2}}<2^{\frac{16}3}<e^{4}$ and
$$
(\sinh{10D})^{\frac{n+2}2\cdot \frac{2}{3n+2}}<
\left\{\begin{array}{lll}
\left(e^{20} D\right)^{\frac{n+2}{3n+2}}<e^{10}D^{\frac{n+2}{3n+2}}& \mbox{ if } &  0<D\leq 2\\
 e^{5D}
 & \mbox{ if } &  D\geq 2
\end{array}\right.
$$
the estimate
$$
\frac{\widetilde{E}_{\mathcal{M}^n}\left(\widetilde{\varepsilon},D\right)^{\frac2{3n+2}}}
{ \widetilde{\gamma}_1^{\frac2{3n+2}}}\leq
\left\{\begin{array}{lll}
\varepsilon^{\frac2{3n+2}}\cdot D\cdot e^{13}n^{\frac{n}{3n+2}}
 & \mbox{ if } & \mathcal{M}^n=\R^n\mbox{ or }\mathcal{M}^n=S^n
\mbox{ or }\mathcal{M}^n=H^n \mbox{ and } D\leq 2\\
\varepsilon^{\frac2{3n+2}}\cdot  e^{6D+6}n^{\frac{n}{3n+2}}
 & \mbox{ if } & \mathcal{M}^n=H^n
\mbox{ and } D\geq 2.
\end{array}\right.
$$
Therefore, using $3n+2\leq 4n$, we may choose
$$
\varepsilon_{\mathcal{M}^n}\left(D\right)=
\left\{\begin{array}{lll}
e^{-28n}n^{-\frac{n}{2}}
 & \mbox{ if } & \mathcal{M}^n=\R^n,\mbox{ or }\mathcal{M}^n=S^n
 \mbox{ and } D\leq \frac{\pi}6,\\
 &&\mbox{ or }\mathcal{M}^n=H^n \mbox{ and } D\leq 2;\\
e^{-30n}n^{-\frac{n}{2}}
\left(\frac{\pi}2-D\right)^{3n+2}
 & \mbox{ if } & \mathcal{M}^n=S^n
 \mbox{ and } \frac{\pi}6\leq D< \frac{\pi}2;\\
 e^{-18D}n^{-\frac{n}{2}}
 & \mbox{ if } & \mathcal{M}^n=H^n
\mbox{ and } D\geq 2
\end{array}\right.
$$
in order to ensure that \eqref{Egammaeta} holds.

Luckily, this choice of $\varepsilon_{\mathcal{M}^n}\left(D\right)$ also satisfies \eqref{tildegamma2cond}.
\endproof

\section{The Euclidean case revisited}
\label{secEuclid}

We write $B^n$ to denote the Euclidean unit ball centered at the origin,
and $X\Delta Y$ to denote the symmetric difference of the sets $X$ and $Y$. The goal of this section is to prove the following theorem.

\begin{theorem}
\label{IsodiametricstabEuclid}
For $n\geq 2$, if $X\subset \R^n$ is measurable with ${\rm diam}X\leq D$ for $D>0$ and
$$
V\left(X\right)\geq (1-\varepsilon)
V\left(\mbox{$\frac{D}{2}$}\,B^n\right)
$$
for $\varepsilon\in\left[0,\frac12\right]$, then there exists a $z\in \R^n$  such that
\begin{equation}
\label{IsodiametricstabEuclidvol}    
V\left(X\Delta \left(z+\mbox{$\frac{D}{2}$}\,B^n\right)\right)\leq \gamma_0\sqrt{\varepsilon}\cdot V\left(\mbox{$\frac{D}{2}$}\,B^n\right)
\end{equation}
for $\gamma_0=cn^{\frac52}(\log n)^5$, and assuming $\varepsilon<c^{-n}$, we have
\begin{equation}
\label{IsodiametricstabEuclidHaus} 
z+\left(1-c\varepsilon^{\frac1{n+1}}\right)\mbox{$\frac{D}{2}$}\,B^n\subset {\rm conv}X\subset z+\left(1-c\varepsilon^{\frac1{n+1}}\right)\mbox{$\frac{D}{2}$}\,B^n
\end{equation}
where  $V\left(\left({\rm conv}X\right)\backslash X\right)\leq \varepsilon\cdot 
V\left(B_{\mathcal{M}^n}\left(z_0,\frac{D}{2}\right)\right)$ and $c>1$ is an absolute constant.
\end{theorem}

The first stability forms of the Brunn-Minkowski inequality were due to Minkowski himself
(see Groemer \cite{Gro93}).
If the distance of convex bodies $K$ and $C$ in $\R^n$ is measured in terms of the Hausdorff distance, then
Diskant \cite{Dis73} and Groemer \cite{Gro88} provided close to be optimal stability versions
(see Groemer \cite{Gro93}). However, the natural distance is in terms volume of the symmetric difference, and the optimal result is due to Figalli, Maggi, Pratelli \cite{FMP09,FMP10} (see Kolesnikov, Milman \cite{KoM22} Section 12 and Klartag, Lehec \cite{KlL} for an improvement of the factor $\gamma^*$ involved in Theorem~\ref{BMstabFMP}).

For convex bodies $K$ and $C$ in $\R^n$, we define the ``homothetic distance'' $A(K,C)$
of convex bodies $K$ and $C$ to be
$$
A(K,C)=\min\left\{V\left(\alpha K\Delta (x+\beta C)\right):\,x\in\R^n\right\}
$$
where $\alpha=V(K)^{\frac{-1}n}$ and
$\beta=V(C)^{\frac{-1}n}$, and $K\Delta Q$ stands for the symmetric difference of $K$ and $Q$.
In addition, let $\sigma(K,C)=\max\left\{\frac{V(C)}{V(K)},\frac{V(K)}{V(C)}\right\}$.
Now Figalli, Maggi, Pratelli \cite{FMP09,FMP10} proved Theorem~\ref{BMstabFMP} up to the factor $\gamma^*$
depending on $n$. The original factor $\gamma^*$ of \cite{FMP10} was improved by Kolesnikov, Milman \cite{KoM22}, 
and Section 12.2 of \cite{KoM22} used the upper bound $c_1\sqrt[4]{n}$ on the Cheeger constant that was available that time for an absolute constant $c_1>0$. However, since then Klartag, Lehec \cite{KlL} proved the upper bound $c_2(\log n)^5$ on the Cheeger constant for an absolute constant $c_2>0$ (note that the Cheeger constant in $\R^n$ is upper bounded by an absolute constant according to the celebrated Kannan-Lov\'asz-Simonovits conjecture).

\begin{theorem}[Figalli-Maggi-Pratelli, Kolesnikov-Milman, Klartag-Lehec]
\label{BMstabFMP}
For convex bodies $K$ and $C$ in $\R^n$, and $\gamma^*=c_0n^{-5}(\log n)^{-10}$ where
$c\in(0,1]$ is an absolute constant, we have
$$
V(K+C)^{\frac1n}\geq (V(K)^{\frac1n}+V(C)^{\frac1n})
\left[1+\frac{\gamma^*}{\sigma(K,C)^{\frac1n}}\cdot A(K,C)^2\right].
$$
\end{theorem}

Here the exponent $2$ of $A(K,C)^2$ is optimal ({\it cf.} Figalli, Maggi, Pratelli \cite{FMP10}).  
We note that 
Diskant \cite{Dis73} and  Groemer \cite{Gro88} verified stability versions of the Brunn-Minkowski inequality in terms of the Hausdorff distance
(see also Groemer \cite{Gro93}).

If $C$ is a convex body in $\R^n$, then its support function is
$$
h_C(u)=\max_{x\in C}\langle x,u\rangle\mbox{ \ for \ }u\in\R^n,
$$
and hence ${\rm diam}\,C=\max_{u\in S^{n-1}}(h_C(u)+h_C(-u))$. Since $h_{C+M}=h_C+h_M$ for another convex body $M$,
it follows that the origin symmetric convex body $\frac12\,(C-C)$ satisfies
\begin{equation}
\label{C-Cball}    
\frac12\,(C-C)\subset  \frac{{\rm diam}\,C}2\,B^n.
\end{equation}

\noindent{\bf Proof of Theorem~\ref{IsodiametricstabEuclid}: } We may assume that $D=2$.
We consider the convex body $K={\rm cl}\,{\rm conv}\,X$, 
 and hence ${\rm diam}\,K\leq 2$ and $V(K)\geq (1-\varepsilon)\kappa_n$ for $V(B^n)=\kappa_n$. It follows from the Isodiametric Inequality Theorem~\ref{Isodiametric} that
$V(K)\leq \kappa_n$. We may translate $K$ in a way such that
$$
A(K,-K)=V(K)^{-1}\cdot V\left(K\Delta (-K)\right)
\geq \kappa_n^{-1}\cdot V\left(K\Delta (-K)\right).
$$
As $1-\varepsilon\geq (1+2\varepsilon)^{-1}$ follows from $\varepsilon\leq \frac12$, we deduce from 
\eqref{C-Cball}, $\sigma(K,-K)=1$ and 
Theorem~\ref{BMstabFMP} that
$$
\kappa_n\geq V\left(\frac12\,(K-K)\right)\geq V(K)\left(1+\gamma^*\cdot \frac{V\left(K\Delta (-K)\right)}{\kappa_n^{2}}\right)^2
\geq \frac{\kappa_n}{1+2\varepsilon}\left(1+\gamma^*\cdot \frac{V\left(K\Delta (-K)\right)}{\kappa_n^{2}}\right)^2;
$$
therefore,
\begin{equation}
\label{symmdiffK-K}    
V\left(K\Delta (-K)\right)\leq \sqrt{\frac{2\varepsilon}{\gamma^*}}  \cdot \kappa_n.  
\end{equation}
In particular, \eqref{C-Cball} and \eqref{symmdiffK-K} yield that $K_0=K\cap (-K)\subset\frac12(K-K)$ satisfies
\begin{equation}
\label{symmdiffK0}    
 \left(1-\sqrt{\frac{\varepsilon}{\gamma^*}}\right)  \cdot \kappa_n\leq V(K)-\frac12\,V\left(K\Delta (-K)\right)
 =V(K_0)\leq \kappa_n,  
\end{equation}
and hence
$$
V\left(K\Delta B^n\right)\leq
V\left(K\Delta K_0\right)+V\left(K_0\Delta B^n\right)
\leq 2\sqrt{\frac{\varepsilon}{\gamma^*}},
$$
proving \eqref{IsodiametricstabEuclidvol}.

Concerning the estimate for the Hausdorff distance of $K$ from $B^n$, let $\delta$ be the maximal depth of a cap of $B^n$ not overlapping with $K_0\subset B^n$. It follows from \eqref{symmdiffK0} and Lemma~\ref{Euclideancap} that
$$
\sqrt{\frac{\varepsilon}{\gamma^*}}  \cdot \kappa_n\geq \frac{2\kappa_{n-1}}{n+1}
\cdot\delta^{\frac{n+1}{2}},
$$
and hence
$$
\delta\leq \left(\frac{\varepsilon}{\gamma^*}\right)^{\frac{1}{n+1}}\cdot
\left(\frac{(n+1)\kappa_n}{2\kappa_{n-1}}\right)^{\frac{2}{n+1}},
$$
proving \eqref{IsodiametricstabEuclidHaus}.
\endproof

\noindent{\bf Acknowledgement: } 
K\'aroly J. B\"or\"oczky and \'Ad\'am Sagmeister are supported by NKFIH project  
 K 132002. We are grateful for the referee's comments improving our paper.

\end{document}